\documentclass[oneside, reqno]{amsart}
\pdfoutput=1
\usepackage{subfiles}
\usepackage[utf8]{inputenc}
\usepackage[T1]{fontenc}
\usepackage[dvipsnames]{xcolor}
\usepackage{ifdraft}
\usepackage{amssymb}
\usepackage{amsmath}
\usepackage{amsthm}
\usepackage{caption}
\usepackage{microtype}
\usepackage{mathtools}
\usepackage{tikz}
\usepackage[bookmarks, bookmarksdepth=subsection, colorlinks, allcolors=Blue, final]{hyperref}
\usepackage{tikz-cd}
\usepackage{lipsum}
\usepackage[foot]{amsaddr}

\ifdraft{
    \usepackage{showlabels}
}{}

\usepackage[colorinlistoftodos, obeyDraft, textsize=tiny]{todonotes}

\usepackage[shortlabels]{enumitem}
\setenumerate[0]{label=(\roman*), font=\normalfont}

\usepackage[capitalize, noabbrev, nameinlink]{cleveref}
\RemoveFromHook{label}[firstaid/cleveref]
\Crefname{equation}{}{}
\Crefname{enumi}{}{}
\Crefname{axiom}{Axiom}{Axioms}
\Crefname{diagram}{Diagram}{Diagrams}
\crefformat{footnote}{#2\footnotemark[#1]#3}

\usepackage[
    backend=biber,
    style=numeric-comp,
    giveninits=true,
    maxbibnames=10,
    minbibnames=10,
    mincrossrefs=10,
    bibencoding=utf8]{biblatex} 

\AtBeginBibliography{
    \raggedright
}

\ifSubfilesClassLoaded{
    \AtEndDocument{%
        \printbibliography
    }
}{}

\makeatletter
\renewcommand{\@mkboth}[2]{}
\makeatother

\usetikzlibrary{decorations.markings}
\usetikzlibrary{cd}

\tikzcdset{
    diagrams={ampersand replacement=\&},
    postmark/.style={
        postaction={
            decorate,
            decoration={
                markings,
                mark={#1}
            }
        }
    },
    display style/.style={
        cells={font=\everymath\expandafter{\the\everymath\displaystyle}},
    },
    dagger arrow/.style={
        postmark=at position 0.5 with {
            \draw[-] (#1,1.5pt) -- (#1,-1.5pt);
        },
        labels={
            inner sep=0.75ex,
        }
    },
    mono/.style={tail},
    epi/.style={two heads},
    dagger mono/.style={dagger arrow={0pt}, mono},
    dagger epi/.style={dagger arrow={-1pt}, epi},
    left adjoint/.style={
        phantom,
        "\dashv"{font=\scriptsize},
        sloped,
        rotate={90}
    },
    iso/.style={"\cong"{sloped}}
}
\input{preamble/widebar}
\theoremstyle{plain}
\newtheorem{theorem}                {Theorem}     [section]

\newtheorem{lemma}        [theorem] {Lemma}
\newtheorem{proposition}  [theorem] {Proposition}
\newtheorem{corollary}    [theorem] {Corollary}
\newtheorem*{theorem*}              {Theorem}

\theoremstyle{definition}

\newtheorem{definition}   [theorem] {Definition}
\newtheorem*{definition*}           {Definition}
\newtheorem{construction} [theorem] {Construction}

\theoremstyle{remark}
\newtheorem{remark}       [theorem] {Remark}

\newtheorem*{remark*}               {Remark}

\renewcommand{\leq}{\leqslant}
\renewcommand{\nleq}{\nleqslant}
\renewcommand{\geq}{\geqslant}

\newcommand{\sgt}{\succ}

\newcommand{\nqquad}{\mkern-36mu}

\newcommand{\blank}{\,\underline{\phantom{n}}\,}

\DeclarePairedDelimiter{\paren}{(}{)}
\DeclarePairedDelimiter{\abs}{|}{|}
\DeclarePairedDelimiter{\norm}{\|}{\|}
\DeclarePairedDelimiterX{\innerProd}[2]{\langle}{\rangle}{#1 \delimsize\vert\mathopen{} #2}
\DeclarePairedDelimiterX{\setb}[2]{\lbrace}{\rbrace}{#1 \,\delimsize\vert\,\mathopen{} #2}
\DeclarePairedDelimiter{\set}{\lbrace}{\rbrace}

\DeclareMathOperator*{\olim}{olim}

\DeclareMathOperator*{\unlim}{unlim}

\newcommand{\Pos}[1]{{#1}_{\geq 0}}
\newcommand{\SelfAdj}[1]{{#1}_{\mathrm{H}}}

\newcommand{\inv}{\ast}

\newcommand{\field}[1]{\mathbb{#1}}
\newcommand{\Nats}{\field{N}}

\newcommand{\Rats}{\field{Q}}

\newcommand{\Reals}{\field{R}}
\newcommand{\PosReals}{\Pos{\Reals}}
\newcommand{\Comps}{\field{C}}
\newcommand{\Quats}{\field{H}}

\newcommand{\cat}[1]{\mathbf{#1}}

\newcommand{\ACon}[1]{{{#1}_{\leq 1}}}
\newcommand{\AIsom}[1]{{{#1}_{1}}}

\newcommand{\Hilb}{\cat{Hilb}}
\newcommand{\Mod}{\cat{Mod}}
\newcommand{\URep}{\cat{URep}}

\newcommand{\Set}{\cat{Set}}
\newcommand{\Vect}{\cat{Vect}}
\newcommand{\C}{\cat{C}}
\newcommand{\J}{\cat{J}}
\newcommand{\D}{\cat{D}}

\newcommand{\pair}[2]{\begin{bsmallmatrix} #1 \\ #2 \end{bsmallmatrix}}
\newcommand{\pairBig}[2]{\begin{bmatrix} #1 \\ #2 \end{bmatrix}}

\newcommand{\He}{\operatorname{He}}

\newcommand{\Sk}{\operatorname{Sk}}

\newcommand*{\CommaCat}[3][]{#2 / #3}

\newcommand{\FinSubs}[1]{\mathcal{F}\!#1}

\makeatletter
\DeclareRobustCommand\vdots{%
  \mathpalette\@vdots{}%
}
\newcommand*{\@vdots}[2]{%
  \sbox0{$#1\cdotp\cdotp\cdotp\m@th$}%
  \sbox2{$#1.\m@th$}%
  \vbox{%
    \dimen@=\wd0 %
    \advance\dimen@ -3\ht2 %
    \kern.5\dimen@
    \dimen@=\wd2 %
    \advance\dimen@ -\ht2 %
    \dimen2=\wd0 %
    \advance\dimen2 -\dimen@
    \vbox to \dimen2{%
      \offinterlineskip
      \copy2 \vfill\copy2 \vfill\copy2 %
    }%
  }%
}
\DeclareRobustCommand\ddots{%
  \mathinner{%
    \mathpalette\@ddots{}%
    \mkern\thinmuskip
  }%
}
\newcommand*{\@ddots}[2]{%
  \sbox0{$#1\cdotp\cdotp\cdotp\m@th$}%
  \sbox2{$#1.\m@th$}%
  \vbox{%
    \dimen@=\wd0 %
    \advance\dimen@ -3\ht2 %
    \kern.5\dimen@
    \dimen@=\wd2 %
    \advance\dimen@ -\ht2 %
    \dimen2=\wd0 %
    \advance\dimen2 -\dimen@
    \vbox to \dimen2{%
      \offinterlineskip
      \hbox{$#1\mathpunct{.}\m@th$}%
      \vfill
      \hbox{$#1\mathpunct{\kern\wd2}\mathpunct{.}\m@th$}%
      \vfill
      \hbox{$#1\mathpunct{\kern\wd2}\mathpunct{\kern\wd2}\mathpunct{.}\m@th$}%
    }%
  }%
}

\newcommand{\bsmallmat}[1]{\begin{bsmallmatrix}#1\end{bsmallmatrix}}
\makeatletter
\newcommand{\fps@diagram}{tbp}
\newcounter{diagram}
\def\ftype@diagram{1}
\def\ext@diagram{lof}
\def\fnum@diagram{\diagramname\ \thediagram}
\def\diagramname{Diagram}
\newenvironment{diagram}{%
\@float{diagram}%
}{%
\end@float
}
\newenvironment{diagram*}{%
\@dblfloat{diagram}%
}{%
\end@dblfloat
}
\makeatother

\usepackage{bm}

\reversemarginpar

\begin{document}
\title[Hilbert \(*\)-categories]{Hilbert \texorpdfstring{\(\bm{*}\)}{*}-categories\\[2ex]\footnotesize Where limits in analysis and\\category theory meet}

\author{Matthew {Di Meglio}}
\author{Chris Heunen}
\address{The University of Edinburgh, 10 Crichton St, Edinburgh, United Kingdom}
\email{m.dimeglio@ed.ac.uk}
\email{chris.heunen@ed.ac.uk}
\date{\today}
\subjclass[2020]{18M40, 46B15, 46L08, 46M15, 46M40, 06F25}
\keywords{Hilbert space, category, directed colimit, monotone completeness}

\begin{abstract}
    This article introduces \textit{Hilbert \(\inv\)\nobreakdash-categories}---an abstraction of categories having algebraic and analytic properties similar to those of the categories of real, complex, and quaternionic Hilbert spaces and bounded linear maps. Other examples include categories of Hilbert W*\nobreakdash-modules and of unitary representations of groupoids. Hilbert \(\inv\)\nobreakdash-categories are ``analytically'' complete in two ways: every bounded increasing sequence of Hermitian endomorphisms has a supremum, and every suitably bounded orthogonal family of parallel morphisms is summable. These ``analytic'' completeness properties are not assumed outright; rather, they are derived, respectively, from two new universal constructions: codirected \(\ell^2\)\nobreakdash-limits of contractions and \(\ell^2\)\nobreakdash-products. In turn, these are built from directed colimits in the wide subcategory of isometries.
\end{abstract}

\maketitle

\setcounter{tocdepth}{1}
\tableofcontents

\section{Introduction}

\textit{C*-categories}~\cite{ghez:w-categories} are a long-established abstract framework for the theory of Hilbert spaces originating in operator algebra. They are multi-object C*\nobreakdash-algebras, that is, they are \textit{\(\inv\)\nobreakdash-categories} (i.e., categories equipped with a functorial and involutive choice of morphism \(f^\inv \colon Y \to X\) for each morphism \(f \colon X \to Y\))\footnote{For consistency with operator algebra, we use the terminology of \(\inv\)\nobreakdash-categories~\cite{ghez:w-categories} instead of dagger categories~\cite{heunenvicary:quantum}. For us, \(\inv\)\nobreakdash-categories are not, \textit{a priori}, enriched in complex vector spaces.} with an involution-compatible enrichment in the category of complex Banach spaces. The involution abstracts adjoints of bounded linear maps, while the Banach space structure of each hom-set abstracts fundamental properties of the operator norm. Notably absent from the definition of C*\nobreakdash-category are universal constructions---the bread and butter of category theory. The theory of C*\nobreakdash-categories is closer in style to operator algebra than to category theory, with many proofs of results about C*\nobreakdash-algebras working \textit{mutatis mutandis} for C*\nobreakdash-categories.

\textit{Pre-Hilbert \(*\)-categories}~\cite{di-meglio:r-star-cats} are a nascent abstract framework for the theory of Hilbert spaces built on ideas from categorical quantum mechanics~\cite{abramsky:categorical-semantics-quantum,heunen:quantum-logic-dagger-order,vicary:completeness--categories-complex,heunen:hilb}. Like C*\nobreakdash-categories, they are \(\inv\)\nobreakdash-categories. Unlike C*\nobreakdash-categories, this is the only assumed \textit{structure}, and all of the assumed \textit{properties} are about universal constructions. The pre\nobreakdash-Hilbert \(\inv\)\nobreakdash-category axioms (see \cref{s:recap}) closely parallel the axioms for abelian categories~(see, e.g.,~\cite[Definition~1.4.1]{borceux:handbook-categorical-algebra-2}). They endow pre\nobreakdash-Hilbert \(\inv\)\nobreakdash-categories with canonical structure that for C*\nobreakdash-categories is assumed from the outset, such as addition of morphisms.

Missing from the theory of pre\nobreakdash-Hilbert \(\inv\)\nobreakdash-categories is a similar universal treatment of analytic completeness. We address this gap with two new universal constructions for pre\nobreakdash-Hilbert \(\inv\)\nobreakdash-categories: \textit{\(\ell^2\)\nobreakdash-limits} of codirected diagrams of contractions, and \textit{\(\ell^2\)\nobreakdash-products} (which are \(\ell^2\)\nobreakdash-limits of finite orthonormal products). We call pre\nobreakdash-Hilbert \(\inv\)\nobreakdash-categories that have all codirected \(\ell^2\)\nobreakdash-limits (and thus all \(\ell^2\)\nobreakdash-products) \textit{Hilbert \(\inv\)\nobreakdash-categories} (see \cref{d:m-category}). This extra axiom makes them the Hilbert analogue of AB3 abelian categories~\cite{grothendieck:quelques-points-dalgebre} (see also~\cite{popescu:abelian-categories}). In our main examples of Hilbert \(\inv\)\nobreakdash-categories (see \cref{s:hilbert-spaces} and \cref{s:examples}), the objects---Hilbert W*\nobreakdash-modules and unitary representations---are in some sense analytically complete. This is not a coincidence. We will see that all Hilbert \(\inv\)\nobreakdash-categories exhibit two ``analytic'' completeness properties: \textit{monotone completeness} and \textit{orthogonal completeness}.

\subsection{Codirected \texorpdfstring{\(\ell^2\)}{l2}-limits and monotone completeness}
\label{s:intro-monotone-complete}

Consider the category \(\Hilb\) of complex Hilbert spaces and bounded linear maps. The Hilbert space
\[
    \ell^2(\Nats) = \setb[\Big]{x \in \Comps^\Nats}{\sum_{n=1}^\infty \norm{x_n}^2 < \infty} \text,
\]
of square-summable sequences looks like a limit in \(\Hilb\) of the codirected diagram
\[
    \begin{tikzcd}[cramped, column sep={0em}]
        \Comps
        \&[4em]
        \Comps^2
            \arrow[l, "{\bsmallmat{
                1 & 0
            }}" swap]
        \&[5em]
        \Comps^3
            \arrow[l, "{\bsmallmat{
                1 & 0 & 0 \\
                0 & 1 & 0
            }}" swap]
        \&[6em]
        \Comps^4
            \arrow[l, "{\bsmallmat{
                1 & 0 & 0 & 0\\
                0 & 1 & 0 & 0\\
                0 & 0 & 1 & 0
            }}" swap]
        \&[2em]
        \cdots
            \arrow[l]
    \end{tikzcd}
\]
of coisometries (orthogonal projections onto closed subspaces), but it is not. While every element \(x \in \ell^2(\Nats)\) induces a cone
\[
    \begin{tikzcd}[cramped, column sep={0em}, row sep={5em}]
        \&[4em]\&[5em]\&[6em]\&[2em]
        \Comps 
            \arrow[dllll, "{x_1}"{swap, pos=0.95, inner sep={1pt}}, in=60, out=-180, looseness=0.6]
            \arrow[dlll, "{\bsmallmat{x_1 \\ x_2}}" {swap, pos=0.96, inner sep={1pt}}, in=70, out=-170, looseness=0.7]
            \arrow[dll, "{\bsmallmat{x_1 \\ x_2 \\ x_3}}"{swap, pos=0.95, inner sep={1pt}}, in=80, out=-160, looseness=0.8]
            \arrow[dl, "{\bsmallmat{x_1 \\ x_2 \\ x_3 \\ x_4}}"{swap, pos=0.92, inner sep={1pt}}, in=90, out=-150, looseness=0.9]
        \\
        \Comps
        \&
        \Comps^2
            \arrow[l, "{\bsmallmat{
                1 & 0
            }}"]
        \&
        \Comps^3
            \arrow[l, "{\bsmallmat{
                1 & 0 & 0 \\
                0 & 1 & 0
            }}"]
        \&
        \Comps^4
            \arrow[l, "{\bsmallmat{
                1 & 0 & 0 & 0\\
                0 & 1 & 0 & 0\\
                0 & 0 & 1 & 0
            }}"]
        \&
        \cdots
            \arrow[l]
    \end{tikzcd}
\]
on the diagram, not every such cone comes from such an element. The ones that do are ``bounded above'' in the sense that there exists \(b \in \PosReals\) such that
\[
    \norm*{\begin{bmatrix} x_1 \\ \vdots \\ x_n \end{bmatrix}}^2
    = \sum_{k=1}^n \norm{x_k}^2 \leq b
\]
for each \(n \in \Nats\). Also
\[
    \norm{x}^2
    = \sum_{n=1}^\infty \norm{x_n}^2
    = \sup_{n \in \Nats} {\norm*{\begin{bmatrix} x_1 \\ \vdots \\ x_n \end{bmatrix}}^2}
    = \lim_{n \to \infty} {\norm*{\begin{bmatrix} x_1 \\ \vdots \\ x_n \end{bmatrix}}^2}
    \text.
\]

The actual universal property of \(\ell^2(\Nats)\) in \(\Hilb\) generalises this situation: the cone formed by the truncation maps \(\ell^2(\Nats) \to \Comps^n\) is terminal among the cones that are suitably ``bounded above'', and the comparison morphism from such a cone to the terminal one encodes a supremum similar to the one above. The new notion of \emph{\(\ell^2\)\nobreakdash-limit} (\cref{d:l2limit}) in pre\nobreakdash-Hilbert \(\inv\)\nobreakdash-categories is an abstraction of this universal property.

Every increasing sequence \((a_n)_{n \in \Nats}\) of strictly positive real numbers forms a cone
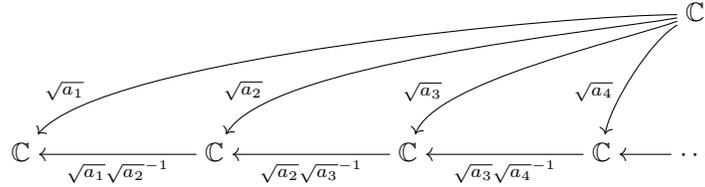
\[
    \begin{tikzcd}[cramped, column sep={0em}, row sep={4em}]
        \&[6em]\&[6em]\&[6em]\&[2em]
        \Comps 
            \arrow[dllll, "\sqrt{a_1}"{swap, pos=0.86, inner sep={1pt}}, in=50, out=-180, looseness=0.4]
            \arrow[dlll, "\sqrt{a_2}" {swap, pos=0.84, inner sep={1pt}}, in=60, out=-170, looseness=0.45]
            \arrow[dll, "\sqrt{a_3}"{swap, pos=0.8, inner sep={1pt}}, in=70, out=-160, looseness=0.5]
            \arrow[dl, "\sqrt{a_4}"{swap, pos=0.73, inner sep={1pt}}, in=80, out=-150, looseness=0.6]
        \\
        \Comps
        \&
        \Comps
            \arrow[l, "\sqrt{a_1}\sqrt{a_2}^{-1}"]
        \&
        \Comps
            \arrow[l, "\sqrt{a_2}\sqrt{a_3}^{-1}"]
        \&
        \Comps
            \arrow[l, "\sqrt{a_3}\sqrt{a_4}^{-1}"]
        \&
        \cdots
            \arrow[l]
    \end{tikzcd}
\]
on a codirected diagram of contractions. If the sequence is bounded above, then it converges to its supremum, and this supremum can be extracted from the comparison morphism from the cone to the \(\ell^2\)\nobreakdash-limit of the diagram. More generally, if an increasing net of strictly positive endomorphisms on a Hilbert space is bounded above, then it converges in the strong operator topology to its supremum, and this supremum can be recovered from an \(\ell^2\)\nobreakdash-limit in a similar way.

Inspired by the recent characterisation of the category of finite-dimensional Hilbert spaces and contractions~\cite{di-meglio:fcon}, we abstract this construction to the setting of Hilbert \(\inv\)\nobreakdash-categories, yielding the following result (\cref{p:monotone-completeness}).

\begin{theorem*}
In a Hilbert \(\inv\)\nobreakdash-category, every increasing net of Hermitian endomorphisms with an upper bound has a supremum.
\end{theorem*}

This property is called \textit{monotone completeness}. It characterises the field \(\Reals\) of real numbers among partially ordered fields~\cite{demarr:partially-ordered-fields}. It is also important in operator algebra~\cite{saitowright:monotone-complete-generic}, featuring in characterisations of W*\nobreakdash-algebras~\cite[Definition~1]{kadison:operator-algebras-faithful} and AW*\nobreakdash-algebras~\cite[Theorem~1.5]{saito:defining-aw-algebras-rickart}.

\subsection{\texorpdfstring{\(\ell^2\)}{l2}-products and orthogonal completeness}

Rather than thinking of the Hilbert space \(\ell^2(\Nats)\) as a codirected \(\ell^2\)\nobreakdash-limit, it is simpler to think of it as being like a product of a copy of \(\Comps\) for each natural number. Of course, it is not actually a product. The (wide) span formed by the coordinate maps \(\ell^2(\Nats) \to \Comps\) is terminal, not among \textit{all} spans, but only among the ones that are suitably ``bounded above''. The notion of \(\ell^2\)\nobreakdash-product (see~\cref{d:l2biproduct}) is an abstraction of this universal property to the setting of pre\nobreakdash-Hilbert \(\inv\)\nobreakdash-categories. It is the pre\nobreakdash-Hilbert \(\inv\)\nobreakdash-category analogue of the C*\nobreakdash-category notion of \textit{infinite direct sum}~\cite{fritz:universal-property-infinite}.

Although \(\ell^2\)\nobreakdash-products are really just \(\ell^2\)\nobreakdash-limits of finite orthonormal products, they deserve their own independent treatment. For example, there is a purely equational characterisation of \(\ell^2\)\nobreakdash-products in Hilbert \(\inv\)\nobreakdash-categories (see \cref{p:l2-biproduct-eqs,p:l2-biproduct-eqs-conv}) that is similar to the purely equational characterisation of finite biproducts in categories enriched in commutative monoids. There is no reason to expect arbitrary codirected \(\ell^2\)\nobreakdash-limits to have such an equational characterisation.

Just as finite biproducts encode finite sums of vectors, infinite \(\ell^2\)\nobreakdash-products encode certain infinite sums of vectors. Consider an orthogonal family \((y_j)_{j \in J}\) of vectors in a Hilbert space~\(X\). For each \(j \in J\), let \(Y_j\) denote the linear span of the vector \(y_j\). As the subspaces \(Y_j\) of \(X\) are closed and pairwise orthogonal, their \(\ell^2\)\nobreakdash-product \(Y\) is again a closed subspace of \(X\). Now, if the finite sums of the vectors are ``bounded above'' in the sense that there exists \(b \in \PosReals\) such that
\[
    \norm[\Big]{\sum_{j \in F} y_j}^2 = \sum_{j \in F} \norm{y_j}^2 \leq b
\]
for all finite \(F \subseteq J\), then the family of vectors is summable~\cite[Theorem~8.2]{halmos:introduction-hilbert}, and its sum---an element of \(Y\)---can be recovered from the universal property of \(\ell^2\)\nobreakdash-products.

To abstract this construction to arbitrary pre\nobreakdash-Hilbert \(\inv\)\nobreakdash-categories (whose assumed structure does not include a norm on each hom-set\footnote{Actually, each hom-set of a pre\nobreakdash-Hilbert \(\inv\)\nobreakdash-category has a canonical \textit{extended seminorm} (like a norm but possibly infinite or zero on some non-zero elements)~\cite{cimpric:2009:quadratic-module} (see also \cite{schotz:equivalence-order-algebraic}). This seminorm does not feature in the present article, but it will be important in future work.}), we use the following new characterisation of sums of orthogonal families (see \cref{d:order-sum,d:l2-family}).

\begin{definition*}
An orthogonal family \((x_j)_{j \in J}\) of elements of an inner product module over an ordered \(\inv\)\nobreakdash-ring
\begin{itemize}
    \item is an \textit{\(\ell^2\)\nobreakdash-family} if there exists \(b \geq 0\) such that
    \[
        \sum_{j \in F} \innerProd{x_j}{x_j} \leq b
    \]
    for all finite \(F \subseteq J\), and
    \item has \textit{order sum} \(s\) if
    \[\innerProd{s}{x_j} = \innerProd{x_j}{x_j}\]
    for each \(j \in J\), and
    \[\innerProd{s}{s} = \sup_{F} \sum_{j \in F} \innerProd{x_j}{x_j}\]
    where \(F\) ranges over the finite subsets of \(J\).
\end{itemize}
\end{definition*}

In inner product spaces over \(\Reals\), \(\Comps\) or \(\Quats\), an orthogonal family is norm summable exactly when it is order summable, in which case the two sums coincide (see \cref{p:hilb:orthogonallycomplete}). Since order sum is a new concept, we take time in \cref{s:sums} to prove that it behaves sensibly.

Using the notion of order sum, we abstract the link between sums of orthogonal vectors and \(\ell^2\)\nobreakdash-products of subspaces into the following theorem (\cref{p:orthogonal-completeness}).

\begin{theorem*}
In a pre\nobreakdash-Hilbert \(\inv\)\nobreakdash-category that has all \(\ell^2\)\nobreakdash-products (e.g., a Hilbert \(\inv\)\nobreakdash-category), every \(\ell^2\)\nobreakdash-family of parallel morphisms has an order sum.
\end{theorem*}

We call this property \textit{orthogonal completeness} (see \cref{d:l2-family}). It is similar to \textit{orthogonal completeness} of commutative semisimple rings~\cite[Definition~4]{abian:direct-product}, \textit{addability of partial isometries} in Baer \(\inv\)\nobreakdash-rings~\cite[Section~11]{berberian:baer-star-rings}, and one of the defining properties of AW*\nobreakdash-modules~\cite{kaplansky:modules-operator-algebras}. To justify why orthogonal completeness should be thought of as an ``analytic'' completeness property, we show that an inner product space is orthogonally complete if and only if it is a Hilbert space (\cref{p:hilb:orthogonallycomplete}).

\subsection{Alternatives to codirected \texorpdfstring{\(\ell^2\)}{l2}-limits}

The idea to use a variant of directed colimits to encode analytic completeness comes from the recent characterisations of categories of Hilbert spaces. The characterisation of the category of Hilbert spaces and bounded linear maps assumes that the wide subcategory of \textit{isometries} has directed colimits~\cite{heunen:hilb,tobin-lack:hilb}. The characterisation of the category of Hilbert spaces and \textit{contractions} assumes that the category itself has directed colimits~\cite{heunen:con}.

At this point, the reader may wonder why we did not merely adopt one of these alternative axioms about codirected limits, and then generalise the proof of analytic completeness given in the cited works. Unfortunately, the approach in the cited works---using Solèr's theorem~\cite{soler:characterization-hilbert-spaces,prestel:solers-characterization-hilbert} to deduce that the \(\inv\)\nobreakdash-ring of endomorphisms of a simple separating object is isomorphic to \(\Reals\), \(\Comps\) or \(\Quats\)---has no hope of explaining the analytic completeness of the endomorphism rings of arbitrary objects because they need not be division rings.

The question still remains: what is the link between our axiom about codirected \(\ell^2\)\nobreakdash-limits and these alternative axioms? The following answer is given in \cref{p:dir-colim-equivalence}.

\begin{theorem*}
The following statements about a pre\nobreakdash-Hilbert \(\inv\)\nobreakdash-category are equivalent:
\begin{enumerate}
    \item it is a Hilbert \(\inv\)\nobreakdash-category,
    \item its wide subcategory of contractions has directed colimits, and
    \item its wide subcategory of isometries has directed colimits.
\end{enumerate}
\end{theorem*}

This theorem affords us the luxury of verifying the simpler and weaker axiom about isometries when checking that particular categories are Hilbert \(\inv\)\nobreakdash-categories, while also using the stronger (but more complicated) axiom about codirected \(\ell^2\)\nobreakdash-limits when developing the theory of Hilbert \(\inv\)\nobreakdash-categories. Furthermore, by virtue of connecting codirected \(\ell^2\)\nobreakdash-limits to (ordinary) directed colimits, this theorem should reassure readers that \textit{codirected \(\ell^2\)\nobreakdash-limit} is a fine category-theoretic notion.

\subsection{Next steps}
\label{s:futurework}

While monotone completeness is already important in operator algebra (see the end of \cref{s:intro-monotone-complete}), it remains secondary to norm completeness, which is part of the definition of C*\nobreakdash-algebra. For order-first abstractions of operator algebras~\cite{schotz:equivalence-order-algebraic,handelman:rings-involution-partially,albeverio:partially-ordered-involutory}, monotone completeness can instead serve as the primary analytic completeness assumption. A forthcoming article will develop this idea.

Tentatively, an \textit{M*\nobreakdash-ring} is a symmetric ordered \(\inv\)\nobreakdash-ring that is monotone complete and orthogonally complete, and a \textit{Hilbert M*\nobreakdash-module} over an M*\nobreakdash-ring \(R\) is an inner product \(R\)\nobreakdash-module that is orthogonally complete. Not only are all W*\nobreakdash-algebras and commutative AW*\nobreakdash-algebras examples of M*\nobreakdash-rings, but so are some \emph{unbounded} operator algebras, such as the Arens algebra \(L^\omega = \bigcap_{p \geq 1} L^p\)~\cite{arens:space-lomega-convex}. Hilbert M*\nobreakdash-modules are an unbounded and non-commutative generalisation of \textit{Kaplansky--Hilbert modules}, which were originally called \textit{AW*\nobreakdash-modules}~\cite{kaplansky:modules-operator-algebras}. The Hilbert M*\nobreakdash-modules over a given M*\nobreakdash-ring \(R\) form a Hilbert \(\inv\)\nobreakdash-category \(\Hilb_R\). Every Hilbert \(\inv\)\nobreakdash-category is conjectured to be of this form: if a Hilbert \(\inv\)\nobreakdash-category~\(\C\) has a separating object \(A\), then the hom\nobreakdash-functor \(\C(A, \blank) \colon \C \to \Set\) ought to factor through the forgetful functor \(\Hilb_{\C(A, A)} \to \Set\) via an equivalence of Hilbert \(\inv\)\nobreakdash-categories. This conjecture is almost a Hilbert analogue of the Gabriel--Popescu theorem~\cite{gabriel-popescu} (see also~\cite{popescu:abelian-categories}).

\subsection{Acknowledgements}
Many thanks go to Steve Lack, Tom Leinster, Maria Pia Solèr and Luuk Stehouwer for their helpful feedback.

\section{Recap of pre-Hilbert \texorpdfstring{\(*\)}{*}-categories}
\label{s:recap}

This section recalls just enough about pre\nobreakdash-Hilbert \(\inv\)\nobreakdash-categories 
to define Hilbert \(\inv\)\nobreakdash-categories, which we will do in the next section. Other concepts and results about pre\nobreakdash-Hilbert \(\inv\)\nobreakdash-categories will be recalled as needed. Consult the original article on pre\nobreakdash-Hilbert \(\inv\)\nobreakdash-categories~\cite{di-meglio:r-star-cats} for more details, including proofs.

\subsection{Elementary concepts}
\label{s:elementary-concepts}
A \textit{\(\inv\)\nobreakdash-category} is a category equipped with a choice of morphism \(f^\inv \colon Y \to X\) for each morphism \(f \colon X \to Y\), such that
\[1^\inv = 1,\qquad (gf)^\inv = f^\inv g^\inv, \qquad\text{and}\qquad (f^{\inv})^\inv = f.\]
The operation \(f \mapsto f^\inv\) is called the \textit{involution}. It generalises adjoints of bounded linear maps between Hilbert spaces.

A morphism \(f\) in a \(\inv\)\nobreakdash-category is \textit{isometric} if \(f^\inv f = 1\), \textit{coisometric} if \(ff^\inv = 1\), and \textit{unitary} if it is both isometric and coisometric. An \textit{isometry} is an isometric morphism and a \textit{coisometry} is a coisometric morphism. A \emph{partial isometry} is a morphism \(f\) satisfying \(f=ff^\inv f\).

In a \(\inv\)\nobreakdash-category with a zero object, a binary product
\[
    \begin{tikzcd}
        X_1
            \arrow[from=r, "r_1" swap]
            \&
        X
            \&
        X_2
            \arrow[from=l, "r_2"]
    \end{tikzcd}
\]
is \textit{orthonormal} if \(r_1\) and \(r_2\) are coisometries and \(r_1{r_2}^\inv = 0\), in which case
\[
    \begin{tikzcd}
        X_1
            \arrow[r, "{r_1}^\inv" swap, shift right]
            \arrow[from=r, "r_1" swap, shift right]
            \&
        X
            \&
        X_2
            \arrow[l, "{r_2}^\inv", shift left]
            \arrow[from=l, "r_2", shift left]
    \end{tikzcd}
\]
is a biproduct. Orthonormal products of arbitrary finite families of objects are defined similarly, and also form biproducts. Orthonormal products generalise orthogonal direct sums of Hilbert spaces. We write
\[
    \begin{tikzcd}
        X_1
            \arrow[from=r, "p_1" swap]
            \&
        X_1 \oplus X_2
            \&
        X_2
            \arrow[from=l, "p_2"]
    \end{tikzcd}
\]
for a chosen orthonormal product of \(X_1\) and \(X_2\), and let \(i_1 = {p_1}^\inv\) and \(i_2 = {p_2}^\inv\), so that
\[
    \begin{tikzcd}
        X_1
            \arrow[r, "i_1", shift left]
            \arrow[from=r, "p_1", shift left]
            \&
        X_1 \oplus X_2
            \&
        X_2
            \arrow[l, "i_2" swap, shift right]
            \arrow[from=l, "p_2" swap, shift right]
    \end{tikzcd}
\]
is a chosen biproduct of \(X_1\) and \(X_2\). We use similar notation for chosen orthonormal products of arbitrary finite families of objects. 

We will also make use of the matrix calculus for biproducts. Given morphisms \(f_{jk} \colon X_j \to Y_k\) for all \(j, k \in \set{1, 2}\), we write the matrix
\[\begin{bmatrix} f_{11} & f_{21} \\ f_{12} & f_{22} \end{bmatrix}\]
for the unique morphism \(f \colon X_1 \oplus X_2 \to Y_1 \oplus Y_2\) such that \[p_k f i_j = f_{jk}\] for all \(j, k \in \set{1, 2}\); matrices of other dimensions are defined similarly. Composition of morphisms between biproducts is by matrix multiplication, and applying the involution to a morphism between orthonormal biproducts is the same as taking the transpose of its matrix and then applying the involution to each of its entries.

An \emph{pre\nobreakdash-Hilbert \(\inv\)\nobreakdash-category} is a \(\inv\)\nobreakdash-category with a zero object, binary orthonormal products and isometric kernels, in which every morphism of the form \(\pair{1}{1} \colon X \to X \oplus X\) is a kernel~\cite[Definition~2.4]{di-meglio:r-star-cats}. Given the other axioms, the last axiom is equivalent to the existence of an enrichment in abelian groups~\cite[Proposition~3.4]{di-meglio:r-star-cats}. Pre\nobreakdash-Hilbert \(\inv\)\nobreakdash-categories are to the category of Hilbert spaces and bounded linear maps what abelian categories are to the category of abelian groups. They are an abstraction of categories exhibiting algebraic features of the theory of Hilbert spaces, including aspects of orthogonal complements, positivity and contractions.

The prototypical pre\nobreakdash-Hilbert \(\inv\)\nobreakdash-categories are the \(\inv\)\nobreakdash-categories \(\Hilb_\Reals\), \(\Hilb_\Quats\) and \(\Hilb = \Hilb_\Comps\) of real, quaternionic, and complex Hilbert spaces and bounded linear maps. Generalising \(\Hilb_\Comps\), the following two families of examples are also of interest:
\begin{itemize}
    \item the \(\inv\)\nobreakdash-category \(\URep_G\) of unitary representations of a groupoid \(G\), and
    \item the \(\inv\)\nobreakdash-category \(\Hilb_A\) of Hilbert W*\nobreakdash-modules over a complex W*\nobreakdash-algebra \(A\).
\end{itemize}
Readers unfamiliar with these categories should not fret: discussion of them is, for the most part, confined to \cref{s:examples}, where they are described in more detail.

\subsection{Order, inner products and contractions}

A \emph{\(\inv\)\nobreakdash-ring} is a ring \(R\) equipped with a function \(\blank^\inv \colon R \to R\), called the \textit{involution}, such that 
\[r^{\inv\inv}=r, \qquad 1^\inv=1, \qquad (rs)^\inv = s^\inv r^\inv, \qquad \text{and}\qquad (r+s)^\inv=r^\inv+s^\inv\]
for all \(r,s \in R\). A \(\inv\)\nobreakdash-ring \(R\) is \textit{anisotropic} if \(r^\inv r = 0\) implies \(r = 0\) for all \(r \in R\). An element \(a\) of a \(\inv\)\nobreakdash-ring is \textit{Hermitian} if \(a^\inv = a\). The set of Hermitian elements of a \(\inv\)\nobreakdash-ring \(R\) is denoted \(\SelfAdj{R}\).

An \textit{ordered \(\inv\)\nobreakdash-ring}~\cite[Definition~5.4]{di-meglio:r-star-cats} is an anisotropic \(\inv\)\nobreakdash-ring \(R\) equipped with a partial order \(\leq\) on \(\SelfAdj{R}\) such that
\[
    0 \leq 1,
    \qquad\qquad\text{and}\qquad\qquad
    a \leq b \;\; \implies \;\; r^\inv ar + c \leq r^\inv br + c
\]
for all \(a, b, c \in \SelfAdj{R}\) and all \(r \in R\). Ordered \(\inv\)\nobreakdash-rings are the \(\inv\)\nobreakdash-ring analogue of \textit{ordered \(\inv\)\nobreakdash-algebras}~\cite{schotz:equivalence-order-algebraic}.

An element \(r\) of an ordered \(\inv\)\nobreakdash-ring is \textit{positive} if \(r \geq 0\), and is \textit{strictly positive}, written \(r \sgt 0\), if it is both positive and invertible. An ordered \(\inv\)\nobreakdash-ring is \textit{symmetric} if \(a \geq 1\) implies \(a \sgt 0\), and is \textit{inverse closed} if \(a \geq b\) and \(b \sgt 0\) implies \(a \sgt 0\). Clearly, every inverse-closed ordered \(\inv\)\nobreakdash-ring is symmetric. If an ordered \(\inv\)\nobreakdash-ring is symmetric, then~\cite[Proposition~6.6]{di-meglio:r-star-cats} it is canonically an ordered \(\inv\)\nobreakdash-algebra over~\(\Rats\) (i.e., it admits a unique ordered \(\inv\)\nobreakdash-ring homomorphism from \(\Rats\)).

Every endohom-set \(\C(A, A)\) of a pre\nobreakdash-Hilbert \(\inv\)\nobreakdash-category \(\C\) is canonically an ordered \(\inv\)\nobreakdash-ring~\cite[Proposition~5.11]{di-meglio:r-star-cats}, with \(f \leq g\) if and only if \(g-f=h^\inv h\) for some morphism \(h \colon A \to X\). Also, these ordered \(\inv\)\nobreakdash-rings are inverse closed~\cite[Proposition~6.5]{di-meglio:r-star-cats}, from which it follows that all pre\nobreakdash-Hilbert \(\inv\)\nobreakdash-categories are canonically enriched in the category of rational vector spaces~\cite[Proposition~6.7]{di-meglio:r-star-cats}.

Let \(R\) be an ordered \(\inv\)\nobreakdash-ring. An \textit{inner product \(R\)\nobreakdash-module}~\cite[Definition~5.5]{di-meglio:r-star-cats} is a right \(R\)\nobreakdash-module \(X\) equipped with a function \(\innerProd{\blank}{\blank} \colon X \times X \to R\), called the \textit{inner product}, that has the following properties:
\begin{itemize}
    \item \textit{linearity}: \(\innerProd{x}{yr + z} = \innerProd{x}{y}r + \innerProd{x}{z}\) for all \(x, y, z \in X\) and all \(r \in R\),
    \item \textit{Hermitianness}: \(\innerProd{x}{y}^\inv = \innerProd{y}{x}\) for all \(x, y \in X\),
    \item \textit{positivity}: \(\innerProd{x}{x} \geq 0\) for all \(x \in X\), and
    \item \textit{anisotropy}: \(\innerProd{x}{x} = 0\) implies \(x = 0\) for all \(x \in X\).
\end{itemize}
Although the positive square root \(\abs{x} = \sqrt{\innerProd{x}{x}}\) may not exist, we often write \(\abs{x}^2\) as shorthand for \(\innerProd{x}{x}\).

\begin{remark}
\label{p:weak-tri-ineq}
For all elements \(x\) and \(y\) of an inner product module,
\[
    \abs{x + y}^2 \leq 2 \abs{x}^2 + 2 \abs{y}^2
    \qquad\text{and}\qquad \abs{x - y}^2 \leq 2 \abs{x}^2 + 2 \abs{y}^2
\]
because \(\abs{x + y}^2 + \abs{x - y}^2 = 2 \abs{x}^2 + 2 \abs{y}^2\).
\end{remark}

\begin{proposition}[Generalised Cauchy--Schwarz]
\label{p:cauchy-schwarz-inverse}
For all elements \(x\) and \(y\) of an inner product module, and all invertible \(a \geq \abs{x}^2\) in its ordered \(\inv\)\nobreakdash-ring of scalars,
\[
    \innerProd{y}{x}a^{-1}\innerProd{x}{y} \leq \abs{y}^2.
\]
\end{proposition}

\begin{proof}
Observe that
\begin{align*}
    0
    &\leq \abs[\big]{y - xa^{-1} \innerProd{x}{y}}^2\\
    &= \abs{y}^2 - 2 \innerProd{y}{x}a^{-1} \innerProd{x}{y} + \innerProd{y}{x}a^{-1}\abs{x}^2a^{-1}\innerProd{x}{y}\\
    &\leq \abs{y}^2 - 2 \innerProd{y}{x}a^{-1} \innerProd{x}{y} + \innerProd{y}{x}a^{-1}aa^{-1}\innerProd{x}{y}\\
    &= \abs{y}^2 - \innerProd{y}{x}a^{-1} \innerProd{x}{y}. \qedhere
\end{align*}
\end{proof}

For all objects \(A\) and \(X\) of a pre\nobreakdash-Hilbert \(\inv\)\nobreakdash-category \(\C\), the hom-set \(\C(A,X)\) is an inner product \(\C(A, A)\)\nobreakdash-module~\cite[Proposition~5.11]{di-meglio:r-star-cats} with \(\innerProd{f}{g} = f^\inv g\).

A morphism \(f\) in a pre\nobreakdash-Hilbert \(\inv\)\nobreakdash-category is~\cite[Definition~8.1]{di-meglio:r-star-cats} a \textit{contraction} if \(f^\inv f \leq 1\). All isometries and coisometries are contractions, and contractions are closed under composition and involution~\cite[Proposition~8.3]{di-meglio:r-star-cats}. For each pre\nobreakdash-Hilbert \(\inv\)\nobreakdash-category~\(\C\), the wide \(\inv\)\nobreakdash-subcategory of contractions of \(\C\) is denoted \(\C_{\leq 1}\).

\section{\texorpdfstring{\(\ell^2\)}{l2}-limits and monotone completeness}\label{s:monotonecompleteness}

The notion of \emph{Hilbert \(\inv\)\nobreakdash-category} is introduced in \cref{s:l2limits}. It is defined in terms of \textit{\(\ell^2\)\nobreakdash-limits} of codirected diagrams of contractions: a new universal construction for pre\nobreakdash-Hilbert \(\inv\)\nobreakdash-categories (also introduced in \cref{s:l2limits}) that generalises the properties of \(\ell^2(\Nats)\) highlighted in \cref{s:intro-monotone-complete}. In particular, \(\ell^2\)\nobreakdash-limits are terminal among the cones for which certain associated increasing nets are bounded above. Moreover, the mediating morphism from such a cone encodes the supremum of the associated net. From this property of the mediating morphisms, we reverse-engineer an ``analytic'' completeness property of Hilbert \(\inv\)\nobreakdash-categories: \emph{monotone completeness}. As suprema of increasing nets feature in the definitions of both \(\ell^2\)\nobreakdash-limit and monotone completeness, their basic properties are reviewed first, in \cref{s:sup-properties}. Finally, to illustrate the new concepts of \(\ell^2\)\nobreakdash-limit and Hilbert \(\inv\)\nobreakdash-category, we prove, in \cref{s:hilbert-spaces}, that \(\Hilb_{\Reals}\), \(\Hilb_{\Comps}\) and \(\Hilb_{\Quats}\) are Hilbert \(\inv\)\nobreakdash-categories; additional examples of Hilbert \(\inv\)\nobreakdash-categories are established in \cref{s:examples}.

\subsection{Suprema of increasing nets}
\label{s:sup-properties}

A \textit{directed set} is a non-empty set equipped with a preorder \(\leq\) in which every pair of elements has an upper bound. For each element \(i\) of a directed set \(J\), the subset
\[
    J_{\geq i} = \setb[\big]{j \in J}{j \geq i}
\]
of \(J\) is a directed set when equipped with the restriction of the preorder of \(J\). 

An \textit{increasing net} in an ordered \(\inv\)\nobreakdash-ring is a family of Hermitian elements \((a_j)_{j \in J}\) indexed by a directed set \(J\) such that \(a_k \geq a_j\) for all \(j \leq k\) in \(J\). An \textit{upper bound} of such an increasing net is an element \(b\) such that \(b \geq a_j\) for all \(j \in J\). The \textit{supremum} (it is unique) of such an increasing net is an upper bound \(a\) such that every upper bound \(b\) satisfies \(b \geq a\). The notions of \textit{decreasing net}, \textit{lower bound} and \textit{infimum} are defined dually. An ordered \(\inv\)\nobreakdash-ring is \textit{monotone complete} if every increasing net with an upper bound has a supremum, or, equivalently, if every decreasing net with a lower bound has an infimum. As usual, we write
\[
    \sup_{j \in J}a_j \qquad\text{and}\qquad \inf_{j \in J}b_j
\]
for the supremum of \((a_j)_{j \in J}\) and the infimum of \((b_j)_{j \in J}\), and
\[
    \sup_{j \in J} a_j < \infty \qquad\text{or}\qquad \inf_{j \in J} b_j > -\infty,
\]
to mean that the supremum or infimum exists.

Next, we list some elementary properties of suprema and infima of monotone nets in ordered \(\inv\)\nobreakdash-rings that will be used throughout the rest of the article without further comment. The dual statements (the ones obtained by swapping increasing with decreasing and suprema with infima) are of course also true. The proofs are routine (see, e.g.,~\cite[Section~2.1]{saitowright:monotone-complete-generic}) and are therefore omitted.

First, the supremum of an increasing net depends only on its ``tail''. More precisely, for each increasing net \((a_j)_{j \in J}\) and each \(i \in J\), the net \((a_j)_{j \in J}\) has a supremum if and only if the net \((a_j)_{j \in J_{\geq i}}\) has a supremum, in which case
\[
    \sup_{j \in J}a_j = \sup_{j \in J_{\geq i}} a_j.
\]


Second, suprema and infima are compatible with addition and, in a weaker sense, multiplication. Indeed, for all increasing nets \((a_j)_{j \in J}\) and \((b_j)_{j \in J}\) that have suprema,
\begin{enumerate}
    \item \(\sup_j {(a_j + b_j)} = \sup_j a_j + \sup_j b_j\),
    \item \(\sup_j qa_j = q \sup_j a_j\) for all \(q \in \Rats_{\geq 0}\),
    \item \(\inf_j (-a_j) = -\sup_j a_j\),
    \item 
    \label{i:sup-conj-iso}
    \(\sup_j r a_j r^\inv = r (\sup_j a_j) r^\inv\) for all invertible elements \(r\), and
    \item \(\sup_j a_j \leq \sup_j b_j\) if \(a_j \leq b_j\) for all \(j \in J\).
\end{enumerate}
In fact, if the ambient ordered \(\inv\)\nobreakdash-ring is sufficiently well behaved, then property \cref*{i:sup-conj-iso} holds even when \(r\) is not invertible.


\begin{proposition}
\label{p:sup-mult}
Let \(R\) be an ordered \(\inv\)\nobreakdash-ring that is inverse closed and monotone complete. For all \(r \in R\) and all increasing nets \((a_j)_{j \in J}\) in \(R\) with an upper bound,
\[
    \sup_{j \in J} ra_jr^\inv = r(\sup_{j \in J} a_j)r^\inv.
\]
\end{proposition}

A symmetric ordered \(\inv\)\nobreakdash-ring \(R\) is \textit{Archimedean} if the decreasing sequence \((4^{-n})_{n \in \Nats}\) in \(R\) has an infimum. The reader should be reassured by the fact that \(R\) is Archimedean if and only if every decreasing sequence of rationals that has an infimum in \(\Rats\) also has an infimum in \(R\); we do not prove this fact as we do not need it. Clearly, if \(R\) is monotone complete, then it is Archimedean.

\begin{lemma}
\label{l:arch}
In an Archimedean symmetric ordered \(\inv\)\nobreakdash-ring,
\[
    \inf_{n \in \Nats} 4^{-n} = 0.
\]
\end{lemma}

\begin{proof}
We have
\(
    \displaystyle
    4 \inf_{n \in \Nats} 4^{-n} = \inf_{n \in \Nats} 4^{-n + 1} = \inf_{n \in \Nats} 4^{-n}.
\)
Hence
\(
    \displaystyle
    3 \inf_{n \in \Nats} 4^{-n} = 0.
\)
\end{proof}

Recall~\cite[Proposition~6.2]{di-meglio:r-star-cats} that if \(a \sgt 0\) and \(b \sgt 0\), then \(a \geq rbr^\inv\) if and only if \(b^{-1} \geq r^\inv a^{-1} r\). In particular, inversion is antimonotonic.
We will use these facts repeatedly in the argument below.

\begin{lemma}
\label{l:sup-mult-inf}
Let \(R\) be an inverse-closed Archimedean ordered \(\inv\)\nobreakdash-ring, and let \((b_j)_{j \in J}\) be a decreasing net in \(R\) with an infimum. If \(\inf_j b_j \sgt 0\), then, for all \(r \in R\),
\[
    \sup_{j \in J} r{b_j}^{-1}r^\inv = r(\inf_{j \in J} b_j)^{-1} r^\inv.
\]
\end{lemma}

\begin{proof}
By inverse closure, \(b_k \sgt 0\) for all \(k \in J\). As inversion is antimonotonic, the net \((r {b_j}^{-1} r^\inv)_{j \in J}\) is increasing and bounded above by \(r (\inf_j b_j)^{-1} r^\inv\). Let \(c\) be another upper bound. For each \(n \in \Nats\), the element \(c + 4^{-n}\) is invertible by inverse closure. For each \(n \in \Nats\) and each \(k \in J\),
\begin{align*}
    r{b_k}^{-1} r^\inv \leq c + 4^{-n},
    \qquad&\text{so}\qquad
    r^\inv (c + 4^{-n})^{-1} r \leq b_k.
    \\\intertext{Hence, for each \(n \in \Nats\),}
    r^\inv (c + 4^{-n})^{-1} r \leq \inf_{j \in J} b_j
    \qquad&\text{so}\qquad
    r(\inf_{j \in J} b_j)^{-1}r^\inv \leq c + 4^{-n}.
\end{align*}
By \cref{l:arch}, it follows that
\[
    r(\inf_{j \in J} b_j)^{-1}r^\inv \leq \inf_{n \in \Nats}(c + 4^{-n}) = c. \qedhere
\]
\end{proof}

\begin{proof}[Proof of \cref{p:sup-mult}]
By assumption, the net \((a_j)_{j \in J}\) has an upper bound \(a\). As \(J\) is a directed set, there exists \(k \in J\). For all \(j \in J_{\geq k}\), as \(a_j - a_k + 1 \geq 1\), also \(a_j - a_k + 1 \sgt 0\) by symmetry; let \(b_j = (a_j - a_k + 1)^{-1}\). Similarly let \(b = (a - a_k + 1)^{-1}\). As inversion is antimonotonic, the net \((b_j)_{j \in J_{\geq k}}\) is decreasing and bounded below by~\(b\). Hence, by monotone completeness, it has an infimum. Also \(\inf_{j \in J_{\geq k}} b_j \geq b \sgt 0\), so \(\inf_{j \in J_{\geq k}} b_j \sgt 0\) by inverse closure. By \cref{l:sup-mult-inf},
\[
    \sup_{j \in J_{\geq k}} {b_j}^{-1} = (\inf_{j \in J_{\geq k}} b_j)^{-1}
    \qquad\text{and}\qquad
    \sup_{j \in J_{\geq k}} r{b_j}^{-1}r^\inv = r(\inf_{j \in J_{\geq k}} b_j)^{-1}r^\inv.
\]
Hence
\begin{align*}
    \sup_{j \in J} r{a_j}r^\inv
    &= \sup_{j \in J_{\geq k}} r{a_j}r^\inv
    = \sup_{j \in J_{\geq k}} r{b_j}^{-1}r^\inv + ra_kr^\inv - rr^\inv
    \\&= r\paren[\big]{\sup_{j \in J_{\geq k}} {b_j}^{-1}}r^\inv + ra_kr^\inv - rr^\inv
    = r (\sup_{j \in J_{\geq k}} a_j) r^\inv
    = r (\sup_{j \in J}a_j) r^\inv \qedhere
\end{align*}
\end{proof}

\subsection{Codirected \texorpdfstring{\(\ell^2\)}{l2}-limits and Hilbert \texorpdfstring{\(*\)}{*}-categories}
\label{s:l2limits}

Every partially ordered set can be viewed as a category: the objects of the category are the elements of the set, and there is a morphism \(a \to b\) whenever \(a \geq b\). From this perspective, the category-theoretic generalisations of increasing net, upper bound, and supremum are, respectively, \textit{codirected diagram}, \textit{cone} and \textit{limit}, which we now define.

A \textit{codirected diagram} \((f_{j \leq k} \colon X_k \to X_j)_{j, k \in J}\) in a category consists of a family of objects \((X_j)_{j \in J}\) indexed by a directed set \(J\), and chosen morphisms \(f_{j \leq k} \colon X_k \to X_j\) for all \(j \leq k\) in \(J\), such that \(f_{j \leq j} = 1\) and \(f_{j \leq k}f_{k \leq l} = f_{j \leq l}\) for all \(j \leq k \leq l\) in \(J\). A \textit{cone} on a codirected diagram \((f_{j \leq k} \colon X_k \to X_j)_{j, k \in J}\) is a family of morphisms \((g_j \colon Y \to X_j)_{j \in J}\) such that \(f_{j \leq k}g_k = g_j\) for all \(j \leq k\) in \(J\). A \textit{limit} of a codirected diagram \((f_{j \leq k} \colon X_k \to X_j)_{j, k \in J}\) is a cone \((f_j \colon X \to X_j)_{j \in J}\) such that, for all cones \((g_j \colon Y \to X_j)_{j \in J}\), there is a unique morphism \(g \colon Y \to X\) such that \(f_jg = g_j\). Up to isomorphism, every codirected diagram has at most one limit. Limits of codirected diagrams are collectively referred to as \textit{codirected limits}. The terms \textit{directed diagram}, \textit{cocone}, \textit{colimit} and \textit{directed colimits} are defined dually.

Consider now a codirected diagram \((f_{j \leq k} \colon X_k \to X_j)_{j,k \in J}\) of contractions in a pre\nobreakdash-Hilbert \(\inv\)\nobreakdash-category \(\C\). If \((g_j \colon Y \to X_j)_{j \in J}\) is a cone on this diagram, then \(({g_j}^\inv g_j)_{j \in J}\) is an increasing net in the ordered \(\inv\)\nobreakdash-ring \(\C(Y,Y)\). Indeed,
\[{g_j}^\inv g_j = {g_k}^\inv {f_{j \leq k}}^\inv f_{j \leq k} g_k \leq {g_k}^\inv g_k\]
for all \(j \leq k\) in \(J\), because the morphism \(f_{j \leq k}\) is a contraction.

\begin{definition}
\label{d:l2limit}
Let \((f_{j \leq k} \colon X_k \to X_j)_{j,k \in J}\) be a codirected diagram of contractions in a pre\nobreakdash-Hilbert \(\inv\)\nobreakdash-category \(\C\).
\begin{itemize}
    \item An \textit{\(\ell^2\)\nobreakdash-cone} on the diagram is a cone \((g_j \colon Y \to X_j)_{j \in J}\) on the diagram such that the increasing net \(({g_j}^\inv g_j)_{j \in J}\) in \(\C(Y,Y)\) is bounded above.
    \item An \textit{\(\ell^2\)\nobreakdash-limit} of the diagram is an \(\ell^2\)\nobreakdash-cone \((f_j \colon X \to X_j)_{j \in J}\) on the diagram such that for all \(\ell^2\)\nobreakdash-cones \((g_j \colon Y \to X_j)_{j \in J}\) on the diagram,
    \begin{enumerate}
        \item there is a unique morphism \(g \colon Y \to X\) such that \[f_j g = g_j\] for all \(j \in J\), and
        \item the unique morphism \(g\) satisfies the equation
        \[g^\inv g = \sup_{j \in J} {g_j}^\inv g_j.\]
    \end{enumerate}
\end{itemize}
\end{definition}
            
The components of an \(\ell^2\)\nobreakdash-cone are not required to be contractions. Nevertheless, the components of an \(\ell^2\)\nobreakdash-limit are necessarily contractions, and, moreover, satisfy the universal property of an (ordinary) limit in the wide subcategory of contractions.

\begin{proposition}
\label{p:l2-lim-from-lim}
If \((f_j \colon X \to X_j)_{j \in J}\) is an \(\ell^2\)\nobreakdash-limit in a pre\nobreakdash-Hilbert \(\inv\)\nobreakdash-category \(\C\) of a codirected diagram \((f_{j \leq k} \colon X_k \to X_j)_{j,k \in J}\) in \(\ACon{\C}\), then all of the morphisms \(f_j\) are contractions, and \((f_j \colon X \to X_j)_{j \in J}\) is a limit in \(\ACon{\C}\) of the diagram.
\end{proposition}

\begin{proof}
Clearly \(f_j 1 = f_j\) for each \(j \in J\). By universality of the \(\ell^2\)\nobreakdash-limit,
\[1 = 1^\inv 1 = \sup_{j \in J} {f_j}^\inv f_j \geq {f_k}^\inv f_k\]
for each \(k \in J\). Hence each morphism \(f_k\) is a contraction.

Next, let \((g_j \colon Y \to X_j)_{j \in J}\) be a cone in \(\ACon{\C}\) on the diagram. As \({g_j}^\inv g_j \leq 1\) for each \(j \in J\), it is also an \(\ell^2\)\nobreakdash-cone in \(\C\) on the diagram. By universality of \(\ell^2\)\nobreakdash-limits, there is a unique morphism \(g \colon Y \to X\) in \(\C\) such that \(f_j g = g_j\) for all \(j \in J\). Also \(\sup_j {g_j}^\inv g_j\) exists, and \(g^\inv g = \sup_{j} {g_j}^\inv g_j \leq 1\). Hence \(g\) is a contraction.
\end{proof}

\cref{p:l2-lim-from-lim} makes mathematically precise the idea that the forgetful functor from \(\ACon{\C}\) to \(\C\) strictly creates limits from the codirected \(\ell^2\)\nobreakdash-limits that \(\C\) admits. The forgetful functor from the wide subcategory of coisometries to \(\ACon{\C}\) also strictly creates codirected limits. Since it is easier to think about isometries than coisometries in concrete categories like \(\Hilb\), we will instead prove the dual result about the wide subcategory \(\AIsom{\C}\) of isometries and creation of directed colimits.

\begin{proposition}
\label{p:dir-colim-create}
Let \(\C\) be a pre\nobreakdash-Hilbert \(\inv\)\nobreakdash-category. The forgetful functor from \(\AIsom{\C}\) to \(\ACon{\C}\) strictly creates all directed colimits that \(\ACon{\C}\) admits.
\end{proposition}

The proof below is adapted from the proof of an analogous result in a different setting~\cite[Lemma~20]{heunen:con}.

\begin{proof}
Let \((f_{j \leq k} \colon X_j \to X_k)_{j, k \in J}\) be a directed diagram in~\(\AIsom{\C}\) that admits a colimit \((f_j \colon X_j \to X)_{j \in J}\) in \(\ACon{\C}\).

Fix an index \(i \in J\). Observe that the restriction \((f_j \colon X_j \to X)_{j \in J_{\geq i}}\) of the colimit is a colimit of the restriction \((f_{j \leq k} \colon X_j \to X_k)_{j,k \in J_{\geq i}}\) of the diagram. Also, notice that \(({f_{i \leq j}}^\inv \colon X_j \to X_{i})_{j \in J_{\geq i}}\) is another cocone on the restricted diagram. Indeed,
\[
    {f_{i \leq j}}^\inv
    = {f_{i \leq j}}^\inv {f_{j \leq k}}^\inv f_{j \leq k}
    = {f_{i \leq k}}^\inv f_{j \leq k}
\]
for all \(j \leq k\) in \(J_{\geq i}\). As its components are contractions, there is a unique contraction \(r_{i} \colon X \to X_{i}\) such that \(r_{i} f_{j} = {f_{i \leq j}}^\inv\) for all \(j \in J_{\geq i}\). In particular,
\[r_{i} f_{i} = {f_{i \leq i}}^\inv = 1^\inv = 1,\]
so \(f_{i}\) is isometric~\cite[Proposition~8.4~(i)]{di-meglio:r-star-cats}.

Consider again the whole diagram \((f_{j \leq k} \colon X_j \to X_k)_{j, k \in J}\). Let \((g_j \colon X_j \to Y)_{j \in J}\) be a cocone in \(\AIsom{\C}\) on this diagram. By universality of the colimit \((f_j \colon X_j \to X)_{j \in J}\) in \(\ACon{\C}\), there is a unique contraction \(g \colon X \to Y\) such that \(gf_j = g_j\) for all \(j \in J\). Now, for all \(j, k \in J\), there exists \(l \in J\) with \(l \geq j\) and \(l \geq k\), and so
\[{f_k}^\inv g^\inv g f_j = {g_k}^\inv g_j = {f_{k \leq l}}^\inv {g_l}^\inv g_l f_{j \leq l} = {f_{k \leq l}}^\inv f_{j \leq l} = {f_{k \leq l}}^\inv {f_l}^\inv f_l f_{j \leq l} = {f_k}^\inv f_j.\]
As the morphisms \(f_j\) are jointly epic in \(\ACon{\C}\), it follows that \(g\) is isometric.
\end{proof}

We now turn our attention to Hilbert \(\inv\)\nobreakdash-categories.

\begin{definition}
\label{d:m-category}
A \textit{Hilbert \(\inv\)\nobreakdash-category} is a pre\nobreakdash-Hilbert \(\inv\)\nobreakdash-category that has \(\ell^2\)\nobreakdash-limits of all codirected diagrams of contractions.
\end{definition}

If a pre\nobreakdash-Hilbert \(\inv\)\nobreakdash-category is a Hilbert \(\inv\)\nobreakdash-category, then, by \cref{p:l2-lim-from-lim,p:dir-colim-create}, its wide subcategory of isometries has directed colimits. The converse of this statement is in fact also true, and will be proved in \cref{s:dir-colim-equivalence}.

We call a pre\nobreakdash-Hilbert \(\inv\)\nobreakdash-category \(\C\) \textit{monotone complete} if, for each object \(A\) of \(\C\), the ordered \(\inv\)\nobreakdash-ring \(\C(A, A)\) is monotone complete. 

\begin{theorem}
\label{p:monotone-completeness}
Every Hilbert \(\inv\)\nobreakdash-category is monotone complete.
\end{theorem}

The case \(\C = \Hilb\) is known as Vigier's theorem (see, e.g.,~\cite[Theorem~4.1.1]{murphy:cstaralgebras}).

\begin{proof}
Let \(A\) be an object of a Hilbert \(\inv\)\nobreakdash-category \(\C\), and let \((a_j)_{j \in J}\) be an increasing net in \(\C(A, A)\) with an upper bound~\(b\). As \(J\) is non-empty, it has an element \(i\). The net \((a_j - a_{i} + 1)_{j \in J_{\geq i}}\) is increasing and bounded above by \(b - a_{i} + 1\). For each \(j \in J_{\geq i}\), the morphism \(a_j - a_{i} + 1\) is invertible~\cite[Proposition~6.5]{di-meglio:r-star-cats}, so~\cite[Proposition~6.3]{di-meglio:r-star-cats} there is an isomorphism \(x_j \colon A \to X_j\) such that \({x_j}^\inv x_j = a_j - a_{i} + 1\). For all \(k \geq j\) in \(J_{\geq i}\), the morphism \(f_{j \leq k} \colon X_k \to X_j\) defined by \(f_{j \leq k} = x_{j} {x_{k}}^{-1}\) is a contraction because
\begin{multline*}
    {f_{j \leq k}}^\inv f_{j \leq k}
    = {x_k}^{-1\inv}{x_j}^\inv x_{j} {x_k}^{-1}
    = {x_k}^{-1\inv}(a_j - a_{i} + 1){x_k}^{-1}
    \\\leq {x_k}^{-1\inv}(a_k - a_{i} + 1){x_k}^{-1}
    = {x_k}^{-1\inv}{x_k}^\inv x_k {x_k}^{-1}
    = 1.
\end{multline*}
In fact, \((f_{j \leq k} \colon X_k \to X_j)_{j, k \in J_{\geq i}}\) is a codirected diagram of contractions, so it has an \(\ell^2\)\nobreakdash-limit \((f_j \colon X \to X_j)_{j \in J_{\geq i}}\). The net \(({x_j}^\inv x_j)_{j \in J_{\geq i}}\) is bounded above by \(b - a_{i} + 1\), so \((x_j \colon A \to X_j)_{j \in J_{\geq i}}\) is an \(\ell^2\)\nobreakdash-cone on the diagram. By universality of \(\ell^2\)\nobreakdash-limits, there is a unique morphism \(x \colon A \to X\) such that \(f_j x = x_j\) for each \(j \in J_{\geq i}\). Also,
\[x^\inv x = \sup_{j \in J_{\geq i}} {x_j}^\inv x_j = \sup_{j \in J_{\geq i}} (a_j - a_{i} + 1).\]
By the properties of suprema discussed in \cref{s:sup-properties},
\[\sup_{j \in J} a_j = \sup_{j \in J_{\geq i}} (a_j - a_{i} + 1) - 1 + a_{i} = x^\inv x - 1 + a_{i}. \qedhere
\]
\end{proof}

\subsection{The Hilbert \texorpdfstring{\(*\)}{*}-category of Hilbert spaces}
\label{s:hilbert-spaces}

Let \(\field{K}\) be \(\Reals\), \(\Comps\) or \(\Quats\). Our goal in this subsection is to show that \(\Hilb_{\field{K}}\) is a Hilbert \(\inv\)\nobreakdash-category. We already know that it is a pre\nobreakdash-Hilbert \(\inv\)\nobreakdash-category~\cite[Section~2]{di-meglio:r-star-cats}, so it suffices to construct \(\ell^2\)\nobreakdash-limits of codirected diagrams of contractions. A morphism \(f \colon X \to Y\) in \(\Hilb_{\field{K}}\) is a contraction in the pre\nobreakdash-Hilbert \(\inv\)\nobreakdash-category sense if and only if it is a contraction in the sense of functional analysis~\cite[Proposition~8.2]{di-meglio:r-star-cats}, that is, \(\norm{fx} \leq \norm{x}\) for all \(x \in X\).

Let \((f_{j \leq k} \colon X_k \to X_j)_{j, k \in J}\) be a codirected diagram of contractions in \(\Hilb_{\field{K}}\). The limit \((\bar{f}_j \colon \widebar{X} \to X_j)_{j \in J}\) of the underlying diagram in the category \(\Vect_{\field{K}}\) of vector spaces over~\(\field{K}\) can be computed using the standard construction of limits from products and equalisers. Its apex \(\widebar{X}\) is the subspace
\[\widebar{X} = \setb[\Big]{x \in \prod_{j \in J} X_j}{f_{j \leq k}x_k = x_j \text{ for all } j \leq k}\]
of the product \(\prod_{j \in J} X_j\), and the morphisms \(\bar{f}_j \colon \widebar{X} \to X_j\) are the restrictions of the product projections \(p_j \colon \prod_{j \in J} X_j \to X_j\) to this subspace.

The \(\ell^2\)\nobreakdash-limit \((f_j \colon X \to X_j)_{j \in J}\) of the original diagram in \(\Hilb_{\field{K}}\) has the following concrete description. Its apex \(X\) is the subspace
\[X = \setb[\Big]{x \in \widebar{X}}{\sup_{j \in J} \norm{x_j} < \infty}\]
of \(\widebar{X}\), equipped with the inner product defined by
\begin{equation}
    \label{e:dir-colim-inn-prod}
    \innerProd{x}{y} = \lim_{j \in J} \innerProd{x_j}{y_j}
\end{equation}
for all \(x, y \in X\), and the maps \(f_j\) are the restrictions of the morphisms \(\bar{f}_j\) to this subspace. We prove this fact in several steps.

\begin{lemma}
The set \(X\) is a subspace of \(\widebar{X}\).
\end{lemma}

\begin{proof}
For all \(x, y \in X\), all \(\lambda \in \field{K}\), and all \(k \in J\),
\[\norm{x_k \lambda + y_k} \leq \norm{x_k \lambda} + \norm{y_k} = \abs{\lambda}\,\norm{x_k} + \norm{y_k} \leq \abs{\lambda}\sup_{j \in J} \norm{x_j} + \sup_{j \in J} \norm{y_j} < \infty.\]
Hence \(x \lambda + y \in X\).
\end{proof}

\begin{lemma}
Equation \cref{e:dir-colim-inn-prod} defines an inner product on \(X\).
\end{lemma}

\begin{proof}
First, we check that the limit on the right-hand side of equation \cref{e:dir-colim-inn-prod} exists. By the polarisation identity (see, e.g.,~\cite[Proposition~2.2]{ghiloni:continuous-slice-functional} for the quaternionic case),
\[
    \lim_{j \in J} \innerProd{x_j}{y_j}
    = \lim_{j \in J} \frac{1}{4} \sum_{e \in B} e\norm{x_je + y_j}^2
    = \frac{1}{4} \sum_{e \in B} e \paren[\big]{\lim_{j \in J} \norm{x_je + y_j}}^2\text,
\]
where
\begin{align*}
    B &= \set{1, -1} &&\text{if \(\field{K} = \Reals\),}\\
    B &= \set{1, -1, \mathrm{i}, -\mathrm{i}} &&\text{if \(\field{K} = \Comps\), and}\\
    B &= \set{1, -1, \mathrm{i}, -\mathrm{i}, \mathrm{j}, -\mathrm{j}, \mathrm{k}, -\mathrm{k}} &&\text{if \(\field{K} = \Quats\),}
\end{align*}
provided that the limit of the net \(\paren[\big]{\norm{x_je + y_j}}_{j \in J}\) exists for each \(e \in B\). But this net is increasing because
\[
    \norm{x_j e + y_j} = \norm{f_{j \leq k}(x_k e + y_k)} \leq \norm{x_k e + y_k}
\]
for all \(j \leq k\) in \(J\), and it is bounded above because
\[\norm{x_j e + y_j} \leq \norm{x_j} + \norm{y_j} \leq \sup_{k \in J} \norm{x_k} + \sup_{k \in J} \norm{y_k} < \infty\]
for all \(j \in J\), so it has a supremum and a limit and these coincide.

We now check the inner product axioms. For all \(x, y, z \in X\) and all \(\lambda \in \field{K}\),
\[
    \innerProd{x}{y\lambda + z}
    = \lim_{j \in J} \innerProd{x_j}{y_j\lambda + z_j}
    = \lim_{j \in J} \innerProd{x_j}{y_j}\lambda + \lim_{j \in J} \innerProd{x_j}{z_j}
    = \innerProd{x}{y} \lambda + \innerProd{x}{z},
\]
and
\[
    \innerProd{x}{y}^\inv = \paren[\big]{\lim_{j \in J} \innerProd{x_j}{y_j}}^\inv = \lim_{j \in J} \innerProd{x_j}{y_j}^\inv = \lim_{j \in J} \innerProd{y_j}{x_j} = \innerProd{y}{x}\text.
\]
Next \(\innerProd{x}{x} = \lim_{j \in J} \norm{x_j}^2 \geq 0\). Finally, suppose that \(\innerProd{x}{x} = 0\). For each \(k \in J\),
\[0 \leq \norm{x_k}^2 \leq \sup_{j \in J} \norm{x_j}^2 = \innerProd{x}{x} = 0,\]
so \(\norm{x_k} = 0\). Hence \(x = 0\).
\end{proof}

\begin{lemma}
The inner product space \(X\) is a Hilbert space.
\end{lemma}

The proof below is a generalisation of the standard proof that the inner product space \(\ell^2(\Nats)\) of square-summable sequences is a Hilbert space.

\begin{proof}
We need to show that every Cauchy sequence in \(X\) converges. Let \((x^{(n)})_{n \in \Nats}\) be a Cauchy sequence in \(X\).

Fix \(j \in J\). The sequence \((x^{(n)}_j)_{n \in \Nats}\) in \(X_j\) is Cauchy because
\begin{equation}
    \label{e:piece-to-whole}
    \norm{x^{(n)}_j - x^{(m)}_j} \leq \sup_{k \in J} \norm{x^{(n)}_k - x^{(m)}_k} = \norm{x^{(n)} - x^{(m)}}
\end{equation}
for all \(m, n \in \Nats\). As \(X_j\) is complete, the sequence converges. Define \(x \in \prod_{j \in J} X_j\) by
\[x_j = \lim_{n \to \infty} x^{(n)}_j.\]
We will show that \(x \in X\) and that \(x^{(n)} \to x\) as \(n \to \infty\). First
\[
    f_{j \leq k} x_k
    = f_{j \leq k} \lim_{n \to \infty} x^{(n)}_k
    = \lim_{n \to \infty} f_{j \leq k} x^{(n)}_k
    = \lim_{n \to \infty} x^{(n)}_j = x_j\]
because \(f_{j \leq k}\) is continuous. Hence \(x \in \widebar{X}\).

Let \(\varepsilon \in \Reals_{>0}\). As \((x^{(n)})_{n \in \Nats}\) is Cauchy, there exists \(N \in \Nats\) such that \[\norm{x^{(n)} - x^{(m)}} < \varepsilon\] for all \(m, n \geq N\). For all \(j \in J\) and all \(m, n \geq N\),
\begin{align*}
    \norm{x^{(n)}_j - x_j}
    &= \norm{(x^{(n)}_j - x^{(m)}_j) + (x^{(m)}_j - x_j)}
    \\&\leq \norm{x^{(n)}_j - x^{(m)}_j} + \norm{x^{(m)}_j - x_j}
    \\&\leq \norm{x^{(n)} - x^{(m)}} + \norm{x^{(m)}_j - x_j}
    \\&< \varepsilon + \norm{x^{(m)}_j - x_j}
\end{align*}
by \cref{e:piece-to-whole}. Letting \(m \to \infty\), it follows that
\[\norm{x^{(n)}_j - x_j} \leq \varepsilon\]
for all \(j \in J\) and all \(n \geq N\). Hence, for all \(n \geq N\), 
\[\norm{x^{(n)} - x} = \sup_{j \in J}\norm{x^{(n)}_j - x_j} \leq \varepsilon,\]
so \(x^{(n)} - x \in X\), and thus \(x = x^{(n)} - (x^{(n)} - x) \in X\). Also,
\[\lim_{n \to \infty} x^{(n)} = x.\]
Thus \(X\) is complete.
\end{proof}

\begin{lemma}
The maps \(f_j\) form an \(\ell^2\)\nobreakdash-limit in \(\Hilb_{\field{K}}\) of the original diagram.
\end{lemma}

Recall~\cite[Proposition~5.12]{di-meglio:r-star-cats} that Hermitian maps \(f\) and \(g\) on a Hilbert space \(X\) satisfy \(f \leq g\) in the pre\nobreakdash-Hilbert \(\inv\)\nobreakdash-category order if and only if \(\innerProd{x}{fx} \leq \innerProd{x}{gx}\) for all \(x \in X\).

\begin{proof}
The maps \(f_j\) are the composites of the linear maps \(\bar{f}_j\) with the subspace inclusion \(X \hookrightarrow \widebar{X}\). Hence \((f_j \colon X \to X_j)_{j \in J}\) is a cone in \(\Vect_{\field{K}}\) on the diagram \((f_{j \leq k} \colon X_k \to X_j)_{j,k \in J}\). It is actually an \(\ell^2\)\nobreakdash-cone in \(\Hilb_{\field{K}}\) on this diagram because its components are contractions: for each \(j \in J\) and each \(x \in X\),
\[\norm{f_j x} = \norm{x_j} \leq \sup_{k \in J} \norm{x_k} = \norm{x}.\]

Let \((g_j \colon Y \to X_j)_{j \in J}\) be an \(\ell^2\)\nobreakdash-cone in \(\Hilb_{\field{K}}\) on the diagram. This means that there is a positive endomorphism \(b \colon Y \to Y\) in \(\Hilb_{\field{K}}\) such that
\[
    {g_j}^\inv g_j \leq b
\]
for all \(j \in J\). To see that the equation 
\[gy = (g_j y)_{j \in J}\]
defines a function \(g \colon Y \to X\), observe that
\[f_{j \leq k} g_k y = g_j y\]
for all \(j \leq k\) in \(J\), and also that
\[
    \norm{g_j y}^2
    = \innerProd{y}{{g_j}^\inv g_j y}
    \leq \innerProd{y}{by}
\]
for all \(j \in J\), where \(\innerProd{y}{by}\) is independent of \(j\). It is easy to check that \(g\) is linear. It is also bounded. Indeed, for all \(y \in Y\),
\[
    \norm{gy}^2 = \sup_{j \in J} \innerProd{y}{{g_j}^\inv g_j y} \leq \innerProd{y}{b y} = \norm{\sqrt{b}y}^2 \leq \norm{\sqrt{b}}^2 \norm{y}^2,
\]
where \(\sqrt{b}\) denotes the unique positive endomorphism on \(Y\) satisfying \(\sqrt{b} \sqrt{b} = b\). It is easy to check uniqueness of \(g\). Finally, for all \(y \in Y\),
\[
    \innerProd{y}{g^\inv g y}
    = \norm{gy}^2
    = \sup_{j \in J} \norm{g_j y}^2
    = \sup_{j \in J} \innerProd{y}{{g_j}^\inv g_j y},
\]
so \(g^\inv g = \sup_{j \in J} {g_j}^\inv g_j\).
\end{proof}

This completes the proof of the following proposition.

\begin{proposition}
\label{p:hilb}
All of \(\Hilb_{\Reals}\), \(\Hilb_{\Comps}\) and \(\Hilb_{\Quats}\) are Hilbert \(\inv\)\nobreakdash-categories.
\end{proposition}

\section{\texorpdfstring{\(\ell^2\)}{l2}-products and orthogonal completeness}

This section is about the second new universal construction---\textit{\(\ell^2\)\nobreakdash-products}---and the corresponding new analytic completeness property---\textit{orthogonal completeness}. \cref{s:sums} introduces \textit{order sums} and \textit{orthogonal completeness}. These notions are, to the best of the authors' knowledge, new, so time is taken to explain why order sums really are generalised sums, and why orthogonal completeness is a good alternative to norm completeness (for inner product \textit{spaces}, it is equivalent to norm completeness). \cref{s:l2products} introduces \textit{\(\ell^2\)\nobreakdash-products}, explains how they generalise finite orthonormal products, and links them to orthogonal completeness. It also shows that every Hilbert \(\inv\)\nobreakdash-category has all \(\ell^2\)\nobreakdash-products, and thus is orthogonally complete.

\subsection{Sums of orthogonal families}\label{s:sums}

For each set \(J\), let \(\FinSubs{J}\) denote the set of finite subsets of \(J\). It is a directed set under subset inclusion. For each family \((a_j)_{j \in J}\) of positive elements in an ordered \(*\)-ring, define
\[
    \sum_{j \in J} a_j = \sup_{F \in \FinSubs J} \sum_{j \in F} a_j.
\]

\begin{definition}
\label{d:order-sum}
An \textit{order sum} of an orthogonal family \((x_j)_{j \in J}\) of elements of an inner product module (over an ordered \(\inv\)\nobreakdash-ring) is an element \(s\) such that
\[\innerProd{s}{x_j} = \abs{x_j}^2\]
for each \(j \in J\), and
\[\abs{s}^2 = \sum_{j \in J} \abs{x_j}^2.\]
An orthogonal family is \textit{order summable} if it has an order sum.
\end{definition}

We begin by establishing some basic properties of order sums, such as uniqueness and compatibility with finite sums and scalar multiplication. This should reassure the reader that order sums are indeed generalised sums.

\begin{lemma}\label{r:order-sum-inf}
    If \((x_j)_{j \in J}\) is an orthogonal family with order sum \(s\), then
    \[
        \sup_{F \in \FinSubs{J}} \abs[\Big]{\sum_{j \in F} x_j}^2 = \abs{s}^2
        \qquad \text{ and } \qquad
        \inf_{F \in \FinSubs{J}} \abs[\Big]{s - \sum_{j \in F}x_j}^2 = 0\text. 
    \]
\end{lemma}
\begin{proof}
    For each \(F \in \FinSubs{J}\),
    \begin{align*}
        \abs[\Big]{\sum_{j \in F} x_j}^2
        = \sum_{j \in F} \abs{x_j}^2 + \sum_{j \in F} \sum_{\substack{k \in F\\k \neq j}} \innerProd{x_j}{x_k}
        = \sum_{j \in F} \abs{x_j}^2.
    \end{align*}
    Hence
    \[
        \sup_{F \in \FinSubs{J}}\abs[\Big]{\sum_{j \in F} x_j}^2
        = \sum_{j \in J} \abs{x_j}^2
        = \abs{s}^2.
    \]
    
    Additionally, for each \(F \in \FinSubs{J}\),
    \[
        \innerProd[\Big]{s}{\sum_{j \in F}x_j}
        = \sum_{j \in F}\innerProd{s}{x_j}
        = \sum_{j \in F}\abs{x_j}^2,
    \]
    so
    \begin{align*}
        \abs[\Big]{s - \sum_{j \in F}x_j}^2
        &= \abs{s}^2 - \innerProd[\Big]{s}{\sum_{j \in F}x_j} - \innerProd[\Big]{\sum_{j \in F}x_j}{s} + \abs[\Big]{\sum_{j \in F}x_j}^2
        = \abs{s}^2 - \sum_{j \in F}\abs{x_j}^2.
    \end{align*}
    Hence
    \[
        \inf_{F \in \FinSubs{J}}\abs[\Big]{s - \sum_{j \in F}x_j}^2
        = \abs{s}^2 - \sum_{j \in J}\abs{x_j}^2
        = \abs{s}^2 - \abs{s}^2
        = 0.
    \qedhere
    \]
\end{proof}

\begin{proposition}
\label{p:order-sum-unique}
Order sums of orthogonal families are unique.
\end{proposition}

\begin{proof}
Let \(s\) and \(t\) be order sums of an orthogonal family \((x_j)_{j \in J}\) of elements of an inner product module. For each \(F \in \FinSubs{J}\),
\[
    \abs{s - t}^2
    = \abs[\Big]{\paren[\Big]{s - \sum_{j \in F}x_j} - \paren[\Big]{t - \sum_{j \in F}x_j}}^2
    \leq 2 \abs[\Big]{s - \sum_{j \in F}x_j}^2 + 2 \abs[\Big]{t - \sum_{j \in F}x_j}^2
\]
by \cref{p:weak-tri-ineq}. Hence
\[
    0 \leq \abs{s - t}^2
    \leq 2 \inf_{F \in \FinSubs{J}}\abs[\Big]{s - \sum_{j \in F}x_j}^2 + 2 \inf_{F \in \FinSubs{J}}\abs[\Big]{t - \sum_{j \in F}x_j}^2
    = 0
\]
by \cref{r:order-sum-inf}. Therefore \(s=t\).
\end{proof}

\begin{proposition}\label{p:order-sum-finite}
    All finite orthogonal families in an inner product module are order summable: the order sum of such a family is merely its sum.
\end{proposition}
\begin{proof}
    Let \((x_j)_{j \in J}\) be a \textit{finite} orthogonal family of elements of an inner product module. As the family is orthogonal, 
    \[
        \innerProd[\Big]{\sum_{k \in J}x_k}{x_j}
        = \sum_{k \in J}\innerProd{x_k}{x_j}
        = \abs{x_j}^2
    \]
    for each \(j \in J\), and
    \[
        \innerProd[\Big]{\sum_{k \in J} x_k}{\sum_{j \in J} x_j}
        = \sum_{j \in J} \innerProd[\Big]{\sum_{k \in J} x_k}{x_j}
        = \sum_{j \in J} \abs{x_j}^2
        \text.
    \]
    Clearly \(\sum_{j \in J} \abs{x_j}^2\) is an upper bound of the increasing net \(\paren[\big]{\sum_{j \in F} \abs{x_j}^2}_{F \in \FinSubs{J}}\). It is actually the least upper bound because \(J \in \FinSubs{J}\).
\end{proof}

\cref{p:order-sum-unique,p:order-sum-finite} ensure that there is no ambiguity in using the notation
\[\sum_{j \in J} x_j\]
for the order sum of a possibly infinite orthogonal family \((x_j)_{j \in J}\).

Next, we show that orthogonal complements are closed under order sums.

\begin{proposition}
\label{p:orthogonal-complement-closed}
Let \(X\) be an inner product module over a symmetric ordered \(\inv\)\nobreakdash-ring. For each order-summable orthogonal family \((x_j)_{j \in J}\) in \(X\) and each \(y \in X\), if \(\innerProd{y}{x_j} = 0\) for all \(j \in J\), then
\[\innerProd[\Big]{y}{\sum_{j \in J}x_j} = 0.\]
\end{proposition}

\begin{proof}
By symmetry, the scalar \(1 + \abs{y}^2\) is invertible. Observe that
\[
    \paren[\big]{1 + \abs{y}^2}^{-1}\abs{y}^2 \paren[\big]{1 + \abs{y}^2}^{-1}
    \leq \paren[\big]{1 + \abs{y}^2}^{-1}\paren[\big]{1 + \abs{y}^2}\paren[\big]{1 + \abs{y}^2}^{-1}
    = \paren[\big]{1 + \abs{y}^2}^{-1}
    \leq 1.
\]
By the generalised Cauchy--Schwarz inequality (\cref{p:cauchy-schwarz-inverse}), it follows that
\[
    \abs[\Big]{\innerProd[\Big]{y\paren[\big]{1 + \abs{y}^2}^{-1}}{\sum_{j \in J}x_j}}^2
    = \abs[\Big]{\innerProd[\Big]{y\paren[\big]{1 + \abs{y}^2}^{-1}}{\sum_{j \in J}x_j - \sum_{j \in F}x_j}}^2
    \leq \abs[\Big]{\sum_{j \in J}x_j - \sum_{j \in F}x_j}^2
\]
for each \(F \in \FinSubs{J}\). Therefore
\[0 \leq \abs[\Big]{\innerProd[\Big]{y\paren[\big]{1 + \abs{y}^2}^{-1}}{\sum_{j \in J}x_j}}^2 \leq \inf_{F \in \FinSubs{J}} \abs[\Big]{\sum_{j \in J}x_j - \sum_{j \in F}x_j}^2 = 0\]
by \cref{r:order-sum-inf}. Hence
\[\innerProd[\Big]{y}{\sum_{j \in J}x_j} = \paren[\big]{1 + \abs{y}^2}\innerProd[\Big]{y\paren[\big]{1 + \abs{y}^2}^{-1}}{\sum_{j \in J}x_j}=0. \qedhere\]
\end{proof}

We now prove an analogue of the inclusion--exclusion principle for order sums, thus showing that they are compatible with addition and subtraction.

\begin{proposition}
Consider an orthogonal family \((x_j)_{j \in J}\) of elements of an inner product module over a symmetric ordered \(\inv\)\nobreakdash-ring. For all \(S, T \subseteq J\),
\[
    \sum_{j \in S \cup T} x_j
    = \sum_{j \in S} x_j + \sum_{j \in T} x_j - \sum_{j \in S \cap T} x_j.
\]
\end{proposition}

Equations like the one above should be interpreted as saying that the order sum (or supremum) on the left-hand side exists and is given by the right-hand side whenever the order sums (or suprema) on the right-hand side exist.

\begin{proof}
First, we show that
\[
    \sum_{j \in S \backslash T} x_j = \sum_{j \in S} x_j - \sum_{j \in S \cap T} x_j.
\]
For all \(k \in S \backslash T\),
\[
    \innerProd[\Big]{\sum_{j \in S} x_j - \sum_{j \in S \cap T} x_j}{x_k}
    =
    \innerProd[\Big]{\sum_{j \in S} x_j}{x_k}
    -
    \innerProd[\Big]{\sum_{j \in S \cap T} x_j}{x_k}
    = \abs{x_k}^2 - 0
    = \abs{x_k}^2
\]
by \cref{p:orthogonal-complement-closed}. Also, for all \(k \in S \cap T\),
\[
    \innerProd[\Big]{\sum_{j \in S} x_j - \sum_{j \in S \cap T} x_j}{x_k}
    =
    \innerProd[\Big]{\sum_{j \in S} x_j}{x_k}
    -
    \innerProd[\Big]{\sum_{j \in S \cap T} x_j}{x_k}
    = \abs{x_k}^2 - \abs{x_k}^2
    = 0.
\]
Applying \cref{p:orthogonal-complement-closed} a second time yields
\[
    \innerProd[\Big]{\sum_{j \in S} x_j - \sum_{j \in S \cap T} x_j}{\sum_{k \in S \cap T} x_k}
    = 0,
\]
so
\[
    \innerProd[\Big]{\sum_{j \in S} x_j}{\sum_{k \in S \cap T} x_k}
    = \abs[\Big]{\sum_{k \in S \cap T} x_k}^2.
\]
It follows that
\begin{align*}
    \sum_{j \in S \backslash T} \abs{x_j}^2
    &= \sum_{j \in S} \abs{x_j}^2 - \sum_{j \in S \cap T} \abs{x_j}^2\\
    &= \abs[\Big]{\sum_{j \in S} x_j}^2 - \abs[\Big]{\sum_{j \in S \cap T} x_j}^2 - \abs[\Big]{\sum_{j \in S \cap T} x_j}^2 + \abs[\Big]{\sum_{j \in S \cap T} x_j}^2\\
    &= \abs[\Big]{\sum_{j \in S} x_j - \sum_{j \in S \cap T} x_j}^2.
\end{align*}
This completes the proof of the first claim.

Next we show that
\[
    \sum_{j \in S \cup T} x_j = \sum_{j \in S \backslash T} x_j + \sum_{j \in T}x_j.
\]
For all \(k \in S \cup T\),
\[
    \innerProd[\Big]{\sum_{j \in S \backslash T} x_j + \sum_{j \in T}x_j}{x_k}
    =
    \innerProd[\Big]{\sum_{j \in S \backslash T} x_j}{x_k}
    +
    \innerProd[\Big]{\sum_{j \in T}x_j}{x_k}
    = \abs{x_k}^2;
\]
for the second equality, consider the disjoint cases \(k \in S \backslash T\) and \(k \in T\) separately, and apply \cref{p:orthogonal-complement-closed} again. As \(S \backslash T\) and \(T\) are disjoint, \cref{p:orthogonal-complement-closed} also yields
\[
    \innerProd[\Big]{\sum_{j \in S \backslash T}x_j}{x_k} = 0
\]
for all \(k \in T\), so also
\[
    \innerProd[\Big]{\sum_{j \in S \backslash T}x_j}{\sum_{k \in T} x_k} = 0.
\]
It follows that
\[
    \sum_{j \in S \cup T} \abs{x_j}^2
    = \sum_{j \in S \backslash T} \abs{x_j}^2 + \sum_{j \in T} \abs{x_j}^2
    = \abs[\Big]{\sum_{j \in S \backslash T} x_j + \sum_{j \in T}x_j}^2 \text.
\]
This completes the proof of the second claim. The overall result is obtained by combining both claims.
\end{proof}

Under stronger assumptions on the ordered \(\inv\)\nobreakdash-ring of scalars, order sums are also compatible with scalar multiplication.

\begin{proposition}
\label{p:scalar-mult}
Let \(R\) be an ordered \(\inv\)\nobreakdash-ring that is inverse closed and monotone complete, and let \(X\) be an inner product \(R\)\nobreakdash-module. For each orthogonal family \((x_j)_{j \in J}\) in \(X\) and each \(r \in R\), the family \((x_jr)_{j \in J}\) is orthogonal and
\[
    \sum_{j \in J} x_j r = \paren[\Big]{\sum_{j \in J}x_j}r.
\]
\end{proposition}

\begin{proof}
The family \((x_jr)_{j \in J}\) is orthogonal because
\(\innerProd{x_jr}{x_kr} = r^\inv \innerProd{x_j}{x_k} r = 0\)
for all \(j, k \in J\) with \(j \neq k\). Also
\[
    \innerProd[\Big]{\paren[\Big]{\sum_{j \in J}x_j}r}{x_k r}
    = r^\inv \innerProd[\big]{\sum_{j \in J}x_j}{x_k}r
    = r^\inv \abs{x_k}^2r
    = \abs{x_kr}^2
\]
for each \(k \in J\), and
\[
    \sum_{j \in J} \abs{x_j r}^2
    = \sum_{j \in J} r^\inv \abs{x_j}^2 r
    = r^\inv \paren[\Big]{\sum_{j \in J}\abs{x_j}^2} r
    = r^\inv \abs[\Big]{\sum_{j \in J}x_j}^2 r
    = \abs[\Big]{\paren[\Big]{\sum_{j \in J}x_j}r}^2
\]
by \cref{p:sup-mult}.
\end{proof}

Our next focus is the completeness property arising from order sums.

\begin{definition}
\label{d:l2-family}
An \textit{\(\ell^2\)\nobreakdash-family} in an inner product module is an orthogonal family \((x_j)_{j \in J}\) of elements of the module for which the increasing net
\[
    \paren[\Big]{\sum_{j \in F} \abs{x_j}^2}_{F \in \FinSubs{J}}
\]
is bounded above. An inner product module is \textit{orthogonally complete} if all \(\ell^2\)\nobreakdash-families are order summable.
\end{definition}

The following proposition should reassure the reader that orthogonal completeness could feasibly serve as an alternative to norm completeness. Intuitively it captures all of the completeness of a module that it does not inherit from its ring of scalars.

\begin{proposition}\label{p:hilb:orthogonallycomplete}
An inner product space over \(\Reals\), \(\Comps\) or \(\Quats\) is a Hilbert space if and only if it is orthogonally complete.
\end{proposition}

\begin{proof}
Let \(X\) be an inner product space over some \(\field{K}\) among \(\Reals\), \(\Comps\) and \(\Quats\).

For the \textit{if} direction, suppose that \(X\) is orthogonally complete. By Zorn's lemma, it has a maximal orthonormal family \((e_j)_{j \in J}\). To deduce that \(X\) is a Hilbert space, we will show that the equation
\[gy = \sum_{j \in J} e_j y_j\]
defines an isometric isomorphism to \(X\) from the Hilbert space
\[
    \ell^2(J) = \setb[\Big]{x \in \field{K}^J}{\sum_{j \in J}\abs{x_j}^2 < \infty}.
\]

We begin by showing that \(g\) is well defined. Let \(y \in \ell^2(J)\). Then \((e_jy_j)_{j \in J}\) is an \(\ell^2\)\nobreakdash-family in \(X\) because
\[
    \sum_{j \in J} \norm{e_j y_j}^2 = \sum_{j \in J} \abs{y_j}^2 < \infty.
\]
But \(X\) is orthogonally complete. Hence the series \(\sum_j e_j y_j\) is order summable, and, by \cref{r:order-sum-inf}, it also norm converges to its order sum.

Next, \(g\) is linear because addition and scalar multiplication preserve norm limits. Also, by the definition of order sum,
\[
    \norm{gy}^2 = \norm[\Big]{\sum_{j \in J} e_j y_j}^2 =
    \sum_{j \in J} \norm{e_j y_j}^2 = \sum_{j \in J} \abs{y_j}^2 = \norm{y}^2
\]
for each \(y \in \ell^2(J)\), so \(g\) is actually isometric.

It remains to show that \(g\) is surjective. Let \(x \in X\). First,
\[
    \sum_{j \in F} \abs{\innerProd{e_j}{x}}^2
    \leq \sum_{j \in F} \abs{\innerProd{e_j}{x}}^2 + \norm[\Big]{x - \sum_{j \in F}e_j \innerProd{e_j}{x}}^2 
    = \norm{x}^2
\]
for each \(F \in \FinSubs{J}\), so \((\innerProd{e_j}{x})_{j \in J} \in \ell^2(J)\). We will show that \(x = g (\innerProd{e_j}{x})_{j \in J}\). Let \(u = x - g(\innerProd{e_j}{x})_{j \in J}\). For each \(k \in J\),
\[
    \innerProd{e_k}{u}
    = \innerProd[\Big]{e_k}{x - \lim_{F \in \FinSubs{J}} \sum_{j \in F} e_j\innerProd{e_j}{x}}
    = \innerProd{e_k}{x} - \lim_{F \in \FinSubs{J}}\sum_{j \in F}\innerProd{e_k}{e_j}\innerProd{e_j}{x}
    = 0
\]
by norm continuity of the inner product. If \(u \neq 0\), then \((e_j)_{j \in J}\) can be extended to a strictly larger orthonormal family by adding \(u \frac{1}{\norm{u}}\), contradicting its maximality. Hence \(u = 0\), and thus \(x = g(\innerProd{e_j}{x})_{j \in J}\).

For the \textit{only if} direction, suppose that \(X\) is norm complete. Let \((x_j)_{j \in J}\) be an \(\ell^2\)\nobreakdash-family in \(X\). 
First, we show that the series \(\sum_{j \in J} x_j\) norm converges. Let \(\varepsilon \in \Reals_{>0}\). As the net
\(
    \paren[\big]{\sum_{j \in F} \norm{x_j}^2}_{F \in \FinSubs{J}}
\)
in \(\Reals\) is increasing and bounded above, it converges (to its supremum), and so is Cauchy. Hence, there exists \(E \in \FinSubs{J}\) such that
\[
    0 \leq \sum_{j \in G} \norm{x_j}^2 - \sum_{j \in F} \norm{x_j}^2 < \varepsilon
\]
for all \(F, G \in \FinSubs{J}\) such that \(E \subseteq F \subseteq G\). In particular, if \(U, V \in \FinSubs{J}\) and \(E \subseteq U\) and \(E \subseteq V\), then \(E \subseteq U \cap V \subseteq U \cup V\), so
\[
    \norm[\Big]{\sum_{j \in U} x_j - \sum_{j \in V} x_j}^2
    =
    \sum_{j \in U \cup V} \norm{x_j}^2 - \sum_{j \in U \cap V} \norm{x_j}^2 < \varepsilon.
\]
Thus the net in \(X\) is Cauchy. Its norm limit 
is an order sum of \((x_j)_{j \in J}\). Indeed
\[
    \innerProd[\Big]{\lim_{F \in \FinSubs{J}}\sum_{j\in F} x_j}{x_k}
    = \lim_{F \in \FinSubs{J}}\innerProd[\Big]{\sum_{j \in F} x_j}{x_k}
    = \lim_{F \in \FinSubs{J}}\innerProd{x_k}{x_k}
    = \innerProd{x_k}{x_k}
\]
for each \(k \in J\), and
\[
    \norm[\Big]{\lim_{F \in \FinSubs{J}}\sum_{j\in F} x_j}^2
    = \lim_{F \in \FinSubs{J}} \norm[\Big]{\sum_{j\in F} x_j}^2
    = \lim_{F \in \FinSubs{J}} \sum_{j \in F} \norm{x_j}^2.\qedhere
\]
\end{proof}

\subsection{\texorpdfstring{\(\ell^2\)}{l2}-products and \texorpdfstring{\(\ell^2\)}{l2}-biproducts}\label{s:l2products}

The following definitions of \(\ell^2\)\nobreakdash-product and \(\ell^2\)\nobreakdash-biproduct each parallel one of the equivalent characterisations~\cite[Theorem~5.1]{fritz:universal-property-infinite} of \textit{infinite direct sum} in a W*\nobreakdash-category.

\begin{definition}
\label{d:l2biproduct}
Let \((X_j)_{j \in J}\) be a family of objects of a pre\nobreakdash-Hilbert \(\inv\)\nobreakdash-category \(\C\).
\begin{itemize}
    \item An \textit{\(\ell^2\)\nobreakdash-span} on \((X_j)_{j \in J}\) is a span \((g_j \colon Y \to X_j)_{j \in J}\) such that the net
    \[
        \paren[\Big]{\sum_{j \in F} {g_j}^\inv g_j}_{F \in \FinSubs{J}}
    \]
    in \(\C(Y, Y)\) is bounded above.
    \item An \textit{\(\ell^2\)\nobreakdash-product} of \((X_j)_{j \in J}\) is an \(\ell^2\)\nobreakdash-span \((r_j \colon X \to X_j)_{j \in J}\) such that for all \(\ell^2\)\nobreakdash-spans \((g_j \colon Y \to X_j)_{j \in J}\),
    \begin{enumerate}
        \item there is a unique morphism \(g \colon Y \to X\) such that
            \[r_j g = g_j\]
        for all \(j \in J\); additionally,
        \item the morphism \(g\) satisfies the equation
        \[
            g^\inv g = \sum_{j \in J} {g_j}^\inv g_j.
        \]
    \end{enumerate}
    \item An \textit{\(\ell^2\)\nobreakdash-biproduct}\footnote{The definition of \(\ell^2\)\nobreakdash-biproduct, unlike the definition of (ordinary) biproduct, is not, strictly speaking, self-dual. This ostensible issue is easily remedied by incorporating into the definition the data of a cospan \((s_j \colon X_j \to X)_{j \in J}\) and the assumption that \(s_j = {r_j}^\inv\) for each \(j \in J\). While this change would strengthen the analogy with (ordinary) biproducts, it would also make the statements of several of the propositions below more awkward.} of \((X_j)_{j \in J}\) is a span \((r_j \colon X \to X_j)_{j \in J}\) such that
    \[
        r_j {r_k}^\inv = \begin{cases}
            1 &\text{if \(j = k\), and}\\
            0 &\text{otherwise}
        \end{cases}
    \]
    for all \(j, k \in J\), and
    \[
        \sum_{j \in J}{r_j}^\inv r_j = 1.
    \]
\end{itemize}
The terms \textit{\(\ell^2\)\nobreakdash-cospan} and \textit{\(\ell^2\)\nobreakdash-coproduct} are defined dually to \(\ell^2\)\nobreakdash-span and \(\ell^2\)\nobreakdash-product.
\end{definition}

In categories enriched in commutative monoids, finite product and finite biproduct are merely two different perspectives on a single concept. This raises the question of whether \(\ell^2\)\nobreakdash-product and \(\ell^2\)\nobreakdash-biproduct are also two different perspectives on a single concept. The answer is yes for sufficiently nice pre\nobreakdash-Hilbert \(\inv\)\nobreakdash-categories. Of the two directions of the equivalence, we will prove the one that holds in \textit{all} pre\nobreakdash-Hilbert \(\inv\)\nobreakdash-categories first.

\begin{proposition}
\label{p:l2-biproduct-eqs}
Every \(\ell^2\)\nobreakdash-product is an \(\ell^2\)\nobreakdash-biproduct.
\end{proposition}

The proof below recycles some ideas from the proofs of \cref{p:l2-lim-from-lim,p:dir-colim-create}.

\begin{proof}
Let \((r_j \colon X \to X_j)_{j \in J}\) be an \(\ell^2\)\nobreakdash-product in a pre\nobreakdash-Hilbert \(\inv\)\nobreakdash-category.

First, \(r_j 1 = r_j\) for all \(j \in J\), so
\[1 = 1^\inv 1 = \sum_{j \in J} {r_j}^\inv r_j.\]

Next, fix \(k \in J\). The morphism \(r_k\) is a contraction because
\[
    {r_k}^\inv r_k \leq \sum_{j \in J} {r_j}^\inv r_j = 1.
\]
Consider the span \((\delta_{jk} \colon X_k \to X_j)_{j \in J}\) defined by
\[
    \delta_{jk} =
    \begin{cases}
        1 &\text{if \(j = k\), and} \\
        0 &\text{otherwise.}
    \end{cases}
\]
It is an \(\ell^2\)\nobreakdash-span because
\[
    \sum_{j \in J} {\delta_{jk}}^\inv \delta_{jk}
    = 1.
\]
By universality of \(\ell^2\)\nobreakdash-products, there is a unique morphism \(s_k \colon X_k \to X\) such that
\[r_j s_k = \delta_{jk}\]
for each \(j \in J\). Also
\[
    {s_k}^\inv s_k = \sum_{j \in J} {\delta_{jk}}^\inv \delta_{jk}
    = 1,
\]
so \(s_k\) is isometric. Actually \(s_k = {r_k}^\inv\) because \(r_k s_k = 1\), and \(r_k\) and \(s_k\) are both contractions~\cite[Proposition~8.4~(i)]{di-meglio:r-star-cats}. Hence \(r_j {r_k}^\inv = \delta_{jk}\) for all \(j \in J\).
\end{proof}

For the other direction of the equivalence, it is helpful to recall the key idea behind the proof of the analogous fact about finite biproducts: for each span
\[
    \begin{tikzcd}[cramped]
        X_1
            \&
        Y
            \arrow[l, "g_1" swap]
            \arrow[r, "g_2"]
            \&
        X_2
    \end{tikzcd},
\]
the unique morphism \(g \colon Y \to X_1 \oplus X_2\) such that \(p_1 g = g_1\) and \(p_2g = g_2\) satisfies
\[g = i_1 g_1 + i_2g_2.\]
The comparison morphism from an \(\ell^2\)\nobreakdash-span to the corresponding \(\ell^2\)\nobreakdash-biproduct is defined by an analogous \textit{order sum}.

A pre\nobreakdash-Hilbert \(\inv\)\nobreakdash-category \(\C\) is \textit{orthogonally complete} if, for all objects \(A\) and \(X\) of \(\C\), the inner product \(\C(A, A)\)\nobreakdash-module \(\C(A, X)\) is orthogonally complete.

\begin{proposition}
\label{p:l2-biproduct-eqs-conv}
In a pre\nobreakdash-Hilbert \(\inv\)\nobreakdash-category that is monotone complete and orthogonally complete, every \(\ell^2\)\nobreakdash-biproduct is an \(\ell^2\)\nobreakdash-product.
\end{proposition}

\begin{proof}
Let \(\C\) be a pre\nobreakdash-Hilbert \(\inv\)\nobreakdash-category that is both monotone complete and orthogonally complete. Consider an \(\ell^2\)\nobreakdash-biproduct \((r_j \colon X \to X_j)_{j \in J}\) in \(\C\). Let \((g_j \colon Y \to X_j)_{j \in J}\) be an \(\ell^2\)\nobreakdash-span in \(\C\). Then \(({r_j}^\inv g_j)_{j \in J}\) is an orthogonal family in \(\C(Y, X)\). Also
\[
    \sum_{j \in F}{g_j}^\inv r_j{r_j}^\inv g_j = \sum_{j \in F}{g_j}^\inv g_j
\]
for each \(F \in \FinSubs{J}\), so the net
\[
    \paren[\Big]{\sum_{j \in F}{g_j}^\inv r_j{r_j}^\inv g_j}_{F \in \FinSubs{J}}
\]
in \(\C(Y, Y)\) is bounded above. By orthogonal completeness of \(\C(Y, X)\), the order sum
\[
    \sum_{j \in J} {r_j}^\inv g_j
\]
exists. Also, for each \(k \in J\),
\[
    r_k \sum_{j \in J} {r_j}^\inv g_j 
    = \sum_{j \in J}r_k{r_j}^\inv g_j
    = r_k {r_k}^\inv g_k 
    = g_k
\]
by the analogue of \cref{p:scalar-mult} for monotone complete pre\nobreakdash-Hilbert \(\inv\)\nobreakdash-categories.

For uniqueness, let \(g \in \C(Y, X)\), and suppose that \(r_j g = g_j\) for each \(j \in J\). Then \(g\) is the order sum of \(({r_j}^\inv g_j)_{j \in J}\). Indeed,
\[
    g^\inv {r_j}^\inv g_j
    = (r_j g)^\inv g_j
    = {g_j}^\inv g_j
\]
for each \(j \in J\), and
\[
    g^\inv g
    = g^\inv \paren[\Big]{\sum_{j \in J}{r_j}^\inv r_j}g
    = \sum_{j \in J}g^\inv {r_j}^\inv r_j g
    = \sum_{j \in J}{g_j}^\inv g_j,
\]
by the analogue of \cref{p:sup-mult} for monotone complete pre\nobreakdash-Hilbert \(\inv\)\nobreakdash-categories.
\end{proof}

Having proved \cref{p:l2-biproduct-eqs}, it should now be apparent that the notion of \(\ell^2\)\nobreakdash-product is a strict generalisation of the notion of finite orthonormal product.

\begin{proposition}
\label{p:fin-l2-is-prod}
A finite span in a pre\nobreakdash-Hilbert \(\inv\)\nobreakdash-category is an \(\ell^2\)\nobreakdash-product if and only if it is an orthonormal product.
\end{proposition}

\begin{proof}
Let \((X_j)_{j \in J}\) be a finite family of objects. Each span 
\((g_j \colon Y \to X_j)_{j \in J}\) on \((X_j)_{j \in J}\) satisfies
\begin{equation}
    \label{e:finite-sup}
    \sup_{F \in \FinSubs{J}}\sum_{j \in F} {g_j}^\inv g_j = \sum_{j \in J} {g_j}^\inv g_j
\end{equation}
because \(J\) is finite. Hence every orthonormal product of \((X_j)_{j \in J}\) is an \(\ell^2\)\nobreakdash-product.

Conversely, suppose that \((r_j \colon X \to X_j)_{j \in J}\) is an \(\ell^2\)\nobreakdash-product of \((X_j)_{j \in J}\). It is actually an (ordinary) product because, by equation \cref{e:finite-sup}, every span on \((X_j)_{j \in J}\) is an \(\ell^2\)\nobreakdash-span. Additionally, by \cref{p:l2-biproduct-eqs}, its components \(r_j\) are orthogonal coisometries. Together, this means that it is an orthonormal product.
\end{proof}

In fact, \(\ell^2\)\nobreakdash-products are just \(\ell^2\)\nobreakdash-limits of finite orthonormal products.

\begin{proposition}
\label{p:l2-prod-from-lim}
For each family \((X_j)_{j \in J}\) of objects in a pre\nobreakdash-Hilbert \(\inv\)\nobreakdash-category, under the canonical isomorphism from the category of spans on \((X_j)_{j \in J}\) to the category of cones on the codirected diagram
\[
    \paren[\big]{p_F \colon \bigoplus_{j \in G} X_j \to \bigoplus_{j \in F}X_j}_{F \subseteq G \in \FinSubs{J}}
\]
of coisometries, \(\ell^2\)\nobreakdash-spans correspond to \(\ell^2\)\nobreakdash-cones, and \(\ell^2\)\nobreakdash-products to \(\ell^2\)\nobreakdash-limits.
\end{proposition}

The canonical isomorphism maps a span
\begin{equation}
    \label{e:corr-span}
    \paren[\big]{g_j \colon Y \to X_j}_{j \in J}
\end{equation}
to the cone
\begin{equation}
    \label{e:corr-cone}
    \paren[\big]{g_F \colon Y \to \bigoplus_{j \in F} X_j}_{F \in \FinSubs{J}}
\end{equation}
whose component \(g_F\) is the unique morphism satisfying
\begin{equation*}
    \label{e:span-to-cone}
    p_j g_F = g_j
\end{equation*}
for each \(j \in F\). Also, it maps a morphism of spans to the same morphism regarded as a morphism of cones. Its inverse maps a cone \cref{e:corr-cone} to the span \cref{e:corr-span} defined by
\begin{equation*}
    \label{e:cone-to-span}
    g_j = g_{\set{j}},
\end{equation*}
and it maps a morphism of cones to the same morphism regarded as a morphism of spans. Checking that these functors are well defined and mutually inverse involves repeated application of universality of (finite) products.

\begin{proof}
For each span \cref{e:corr-span} and cone \cref{e:corr-cone} that correspond to each other,
\[
    {g_F}^\inv g_F
    = {g_F}^\inv \paren[\Big]{\sum_{j \in F}{p_j}^\inv p_j} g_F
    = \sum_{j \in F} {g_F}^\inv {p_j}^\inv p_j g_F
    = \sum_{j \in F} {g_j}^\inv g_j
\]
for each \(F \in \FinSubs{J}\). Hence \cref{e:corr-span} is an \(\ell^2\)\nobreakdash-span if and only if \cref{e:corr-cone} is an \(\ell^2\)\nobreakdash-cone, and so~\cref{e:corr-span} is an \(\ell^2\)\nobreakdash-product if and only if \cref{e:corr-cone} is an \(\ell^2\)\nobreakdash-limit.
\end{proof}

\begin{corollary}
\label{p:m-l2-products}
Every Hilbert \(\inv\)\nobreakdash-category has all \(\ell^2\)\nobreakdash-products.
\end{corollary}

Now that we know the connection between \(\ell^2\)\nobreakdash-products and \(\ell^2\)\nobreakdash-limits, it is easy to extract from \cref{s:hilbert-spaces} a description of \(\ell^2\)\nobreakdash-products in \(\Hilb_{\field{K}}\) where \(\field{K}\) is \(\Reals\), \(\Comps\) or \(\Quats\). First, the Hilbert space \(\ell^2(\Nats)\) of square-summable sequences is merely the \(\ell^2\)\nobreakdash-product in \(\Hilb_{\field{K}}\) of a copy of~\(\field{K}\) for each natural number. More generally, the \(\ell^2\)\nobreakdash-product in \(\Hilb_{\field{K}}\) of a family \((X_j)_{j \in J}\) of Hilbert spaces is the Hilbert space
\[
    \setb[\Big]{x \in \prod_{j \in J} X_j}{\sum_{j \in J} \norm{x_j}^2 < \infty}
\]
equipped with the inner product
\[
    \innerProd{x}{y}
    = \lim_{F \in \FinSubs{J}} \sum_{j \in F} \innerProd{x_j}{y_j}.
\]

\cref{s:elementary-concepts} introduced standard notation for chosen finite orthonormal products. Extending this notation, we write
\[
    \paren[\Big]{p_k \colon \bigoplus_{j \in J} X_j \to X_k}_{k \in J}
\]
for a chosen \(\ell^2\)\nobreakdash-product of a family \((X_j)_{j \in J}\) of objects, and let
\(
    i_j = {p_j}^\inv
\)
for each \(j \in J\). As the involution is a self-duality, the tuple
\[
    \paren[\Big]{i_k \colon X_k \to \bigoplus_{j \in J} X_j}_{k \in J}
\]
is then an \(\ell^2\)\nobreakdash-coproduct of \((X_j)_{j \in J}\). We also write \(X^{\oplus J}\) for the \(\ell^2\)\nobreakdash-product of a copy of \(X\) for each element of a set \(J\).

The following proposition is the \(\ell^2\)\nobreakdash-product analogue of \cref{p:l2-lim-from-lim,p:dir-colim-create}.

\begin{proposition}
\label{p:l2-coproduct-cone-isometries}
Consider an \(\ell^2\)\nobreakdash-coproduct \((s_j \colon X_j \to X)_{j \in J}\) in a pre\nobreakdash-Hilbert \(\inv\)\nobreakdash-category. For each orthogonal cospan \((m_j \colon X_j \to Y)_{j \in J}\) of isometries, there is a unique isometry \(m \colon X \to Y\) such that \(ms_j = m_j\) for each \(j \in J\).
\end{proposition}

\begin{proof}
First, we show that \((m_j \colon X_j \to Y)_{j \in J}\) is an \(\ell^2\)\nobreakdash-cospan. For each \(F \in \FinSubs{J}\),
\[\sum_{j \in F} m_j {m_j}^\inv\]
is a projection (i.e., Hermitian idempotent) because
\[
    \paren[\Big]{\sum_{j \in F} m_j {m_j}^\inv}\paren[\Big]{\sum_{k \in F} m_k {m_k}^\inv}
    = \sum_{j \in F}\sum_{k \in F} m_j \delta_{jk} {m_k}^\inv
    = \sum_{j \in F} m_j {m_j}^\inv.
\]
It follows that
\[
    1 - \sum_{j \in F} m_j {m_j}^\inv
\]
is also a projection, and so is positive. Hence the net
\[
    \paren[\Big]{\sum_{j \in F} m_j {m_j}^\inv}_{F \in \FinSubs{J}}
\]
in \(\C(Y, Y)\) is bounded above by \(1\), as required.

By the dual of \cref{p:l2-prod-from-lim}, the \(\ell^2\)\nobreakdash-cospan \((m_j \colon X_j \to Y)_{j \in J}\) corresponds to an \(\ell^2\)\nobreakdash-cocone
\[
    \paren[\big]{m_F \colon \bigoplus_{j \in F} X_j \to Y}_{F \in \FinSubs{J}}
\]
on the directed diagram
\[
    \paren[\big]{i_{F} \colon \bigoplus_{j \in F} X_j \to \bigoplus_{j \in G}X_j}_{F \subseteq G \in \FinSubs{J}}
\]
of isometries, and the \(\ell^2\)\nobreakdash-coproduct \((s_j \colon X_j \to X)_{j \in J}\) corresponds to an \(\ell^2\)\nobreakdash-colimit
\[
    \paren[\big]{s_F \colon \bigoplus_{j \in F} X_j \to X}_{F \in \FinSubs{J}}
\]
of the same directed diagram. By universality of \(\ell^2\)\nobreakdash-colimits, there is a unique morphism \(m \colon X \to Y\) such that
\[
    ms_F = m_F
\]
for each \(F \in \FinSubs{J}\), or, equivalently (by the dual of \cref{p:l2-prod-from-lim}), such that
\[
    ms_j = m_j
\]
for each \(j \in J\). The morphism \(m\) is isometric by \cref{p:l2-lim-from-lim,p:dir-colim-create}.
\end{proof}

\begin{proposition}
\label{p:orthogonal-completeness}
Every pre\nobreakdash-Hilbert \(\inv\)\nobreakdash-category with all \(\ell^2\)\nobreakdash-products is orthogonally complete.
\end{proposition}

\begin{proof}
Let \(\C\) be a pre\nobreakdash-Hilbert \(\inv\)\nobreakdash-category that has \(\ell^2\)\nobreakdash-products. Fix objects \(A\) and \(X\) of \(\C\), and consider an \(\ell^2\)\nobreakdash-family \((x_j)_{j \in J}\) in the inner product \(\C(A,A)\)\nobreakdash-module \(\C(A, X)\). For each \(j \in J\), there is~\cite[Section~4.2]{di-meglio:r-star-cats} an epimorphism \(y_j \colon A \to Y_j\) and an isometry \(m_j \colon Y_j \to X\) such that \(x_j = m_j y_j\).

First, the span \((y_j \colon A \to Y_j)_{j \in J}\) is an \(\ell^2\)\nobreakdash-span.
Indeed,
\[
    \sum_{j \in F} {y_j}^\inv y_j
    = \sum_{j \in F} {y_j}^\inv {m_j}^\inv m_j y_j
    = \sum_{j \in F} {x_j}^\inv x_j
\]
for all \(F \in \FinSubs{J}\) because each \(m_j\) is isometric, and the net
\[
    \paren[\Big]{\sum_{j \in F} {x_j}^\inv x_j}_{F \in \FinSubs{J}}
\]
in \(\C(A, A)\) is bounded above because \((x_j)_{j \in J}\) is an \(\ell^2\)\nobreakdash-family in \(\C(A, X)\). By universality of \(\ell^2\)\nobreakdash-products, there is a unique morphism
\[
    y \colon A \to \bigoplus_{j \in J} Y_j
\]
such that \(p_k y = y_k\) for all \(k \in J\). Additionally,
\[
    y^\inv y = \sum_{j \in J} {y_j}^\inv y_j = \sum_{j \in J}{x_j}^\inv x_j.
\]

Next, we show that the cospan \((m_j \colon Y_j \to X)_{j \in J}\) of isometries is orthogonal. Fix indices \(j, k \in J\), and suppose that \(j \neq k\). By definition, \(\ell^2\)\nobreakdash-families are orthogonal. Hence \({x_j}^\inv x_k = 0\), and so
\[
    {y_j}^\inv {m_j}^\inv m_{k} y_k = {x_j}^\inv x_k = 0 = {y_j}^\inv 0 y_k.
\]
As \(y_k\) is epic and \({y_j}^\inv\) is monic, it follows that \({m_j}^\inv m_{k} = 0\). By \cref{p:l2-coproduct-cone-isometries}, there is a unique isometry
\[
    m \colon \bigoplus_{j \in J} Y_j \to X
\]
such that \(mi_k = m_k\) for each \(k \in J\).

Let \(x = my\). Then
\[
    x = \sum_{j \in J} x_j
\]
because
\[
    x^\inv x_j
    = y^\inv m^\inv m_j y_j
    = y^\inv m^\inv m i_j y_j
    = y^\inv i_j y_j
    = (p_j y)^\inv y_j
    = {y_j}^\inv y_j
    = {x_j}^\inv x_j\]
for each \(j \in J\), and
\[
    x^\inv x
    = y^\inv m^\inv m y
    = y^\inv y
    = \sum_{j \in J}{x_j}^\inv x_j. \qedhere
\]
\end{proof}

\begin{corollary}
Every Hilbert \(\inv\)\nobreakdash-category is orthogonally complete.
\end{corollary}

\begin{proof}
Combine \cref{p:orthogonal-completeness} and \cref{p:m-l2-products}.
\end{proof}

In light of \cref{p:l2-biproduct-eqs-conv,p:l2-biproduct-eqs}, for Hilbert \(\inv\)\nobreakdash-categories, the concepts of \(\ell^2\)\nobreakdash-product and \(\ell^2\)\nobreakdash-biproduct coincide.

\section{Equivalent definitions of Hilbert \texorpdfstring{\(*\)}{*}-category}
\label{s:dir-colim-equivalence}

The goal of this section is to prove the following theorem, which gives two alternative characterisations of Hilbert \(\inv\)\nobreakdash-categories.

\begin{theorem}\label{p:dir-colim-equivalence}
    The following statements about a pre\nobreakdash-Hilbert \(\inv\)\nobreakdash-category \(\C\) are equivalent:
    \begin{enumerate}
        \item \(\C\) is a Hilbert \(\inv\)\nobreakdash-category;
        \item \(\ACon{\C}\) has directed colimits;
        \item \(\AIsom{\C}\) has directed colimits.
    \end{enumerate}
\end{theorem}

We already showed that (i) implies (ii) in \cref{p:l2-lim-from-lim}, and that (ii) implies (iii) in \cref{p:dir-colim-create}. Proving the converses of these implications is much more challenging. Our proof that (iii) implies (ii) is broken into three steps. First, in \cref{s:dir-colim-isoms}, the wide subcategory of contractions is shown to have colimits of directed diagrams of isometries. Second, in \cref{s:dir-colim-coisom}, it is also shown to have limits of codirected diagrams of isometries. Third, in \cref{s:dir-colim-con}, the general case is reduced to these two special cases. Proving that (ii) implies (i) is the focus of \cref{s:con-l2}.

\subsection{Directed colimits of isometries}
\label{s:dir-colim-isoms}

Recall~\cite[Definition~8.6]{di-meglio:r-star-cats} that a \emph{codilation} of a morphism \(f \colon X \to Y\) in a \(\inv\)\nobreakdash-category is a triple \((T,t_1,t_2)\) consisting of isometries \(t_1 \colon X \to T\) and \(t_2 \colon Y \to T\) such that \({t_2}^\inv t_1 = f\); and a \emph{codilator} of \(f\) is a codilation \((F,f_1,f_2)\) of \(f\) such that, for all codilations \((T,t_1,t_2)\) of \(f\), there is a unique isometry \(t \colon F \to T\) such that \(t_1 = tf_1\) and \(t_2 = tf_2\).

\begin{proposition}
\label{p:dir-colim-preserve}
Let \(\D\) be a \(\inv\)\nobreakdash-category in which every morphism has a codilator. If \(\AIsom{\D}\) has directed colimits, then the forgetful functor from \(\AIsom{\D}\) to \(\D\) preserves them.
\end{proposition}

Our proof resembles the proof~\cite[Theorem~3.5]{frei:limits-categories-relations} that the embedding of an AB5 abelian category into its category of corelations preserves directed colimits.

\begin{proof}
\label{c:dir-colim-codilator-comm}
Let \((f_{j \leq k} \colon X_j \to X_k)_{j,k \in J}\) be a directed diagram in \(\AIsom{\D}\), let \((f_j \colon X_j \to X)_{j \in J}\) be its colimit, and let \((a_j \colon X_j \to A)_{j \in J}\) be a cocone on the diagram in \(\D\). By assumption, each morphism \(a_j\) has a codilator.
Let us construct the morphisms in~\(\D\) that are depicted in \cref{d:codil-colim-1,d:codil-colim-2}.
\begin{diagram}
\centering
\begin{minipage}{0.5\textwidth}
    \centering
    \begin{tikzpicture}[x=7.5em, y=4.8em, commutative diagrams/.cd, every diagram]
        \begin{scope}[commutative diagrams/.cd, every node]
            \node (Xj) at (0,2) {\(X_j\)};
            \node (Xk) at (1,2) {\(X_k\)};
            \node (X) at (2,1.65) {\(X\)};
            \node (Al) at (0,0) {\(A\)};
            \node (Am) at (1,0) {\(A\)};
            \node (Ar) at (2,-0.35) {\(A\)};
        \end{scope}
        
        \path[
            commutative diagrams/.cd,
            every arrow,
            every label
        ]
        (Al)
            edge[
                commutative diagrams/.cd, mono
            ] node{\(1\)} (Am)
            edge[
                in=180, out=-45, looseness=.4,
                commutative diagrams/.cd, mono
            ] node[swap]{\(1\)} ([yshift=-0.8ex]Ar.west)
        (Xk)
            edge node[swap]{\(a_k\)} (Am)
        (Xj)
            edge node[swap]{\(a_j\)} (Al)
            edge[
                commutative diagrams/.cd, mono
            ] node{\(f_{j \leq k}\)} (Xk)
            edge[
                in=180, out=-45, looseness=.4,
                commutative diagrams/.cd, crossing over, mono
            ] node[pos=.325,swap]{\(f_j\)} ([yshift=-0.8ex]X.west)
        (Xk)
            edge[
                in=180, out=-45, looseness=.5,
                commutative diagrams/.cd, mono
            ] node{\(f_k\)} ([yshift=0.6ex]X.west)
        (X)
            edge node[swap]{\(a\)} (Ar)
        (Am)
            edge[in=180, out=-45, looseness=.5, commutative diagrams/.cd, mono] node{\(1\)} ([yshift=0.6ex]Ar.west)
        ;
    \end{tikzpicture}
    \caption{}
    \label{d:codil-colim-1}
\end{minipage}%
\begin{minipage}{0.5\textwidth}
    \centering
    \begin{tikzpicture}[x=7.5em, y=4.8em, commutative diagrams/.cd, every diagram]
        \begin{scope}[commutative diagrams/.cd, every node]
            \node (Xj) at (0,2) {\(X_j\)};
            \node (Xk) at (1,2) {\(X_k\)};
            \node (X) at (2,1.65) {\(X\)};
            \node (Yj) at (0,1) {\(Y_j\)};
            \node (Yk) at (1,1) {\(Y_k\)};
            \node (Y) at (2,0.65) {\(Y\)};
            \node (Al) at (0,0) {\(A\)};
            \node (Am) at (1,0) {\(A\)};
            \node (Ar) at (2,-0.35) {\(A\)};
        \end{scope}
        
        \path[
            commutative diagrams/.cd,
            every arrow,
            every label,
            mono
        ]
        (Al)
            edge node{\(1\)} (Am)
            edge node[swap, pos=.3]{\(t_j\)} (Yj)
            edge[
                in=180, out=-45, looseness=.4
            ] node[swap]{\(1\)} ([yshift=-0.8ex]Ar.west)
        (Am)
            edge node[pos=.3, swap]{\(t_k\)} (Yk)
            edge[
                in=180, out=-45, looseness=.5
            ] node{\(1\)} ([yshift=0.6ex]Ar.west)
        (Ar)
            edge node[swap]{\(t\)} (Y)
        (Yj)
            edge node{\(g_{j \leq k}\)} (Yk)
            edge[
                in=180, out=-45, looseness=.4,
                commutative diagrams/.cd, crossing over
            ] node[swap,pos=.325]{\(g_j\)} ([yshift=-0.8ex]Y.west)
        (Yk)
            edge[
                in=180, out=-45, looseness=.5
            ] node[pos=.4]{\(g_k\)} ([yshift=0.6ex]Y.west)
        (Xk)
            edge node[pos=.7]{\(s_k\)} (Yk)
        (Xj)
            edge node[pos=.7]{\(s_j\)} (Yj)
            edge node{\(f_{j \leq k}\)} (Xk)
            edge[
                in=180, out=-45, looseness=.4,
                commutative diagrams/.cd, crossing over
            ] node[pos=.325,swap]{\(f_j\)} ([yshift=-0.8ex]X.west)
        (Xk)
            edge[
                in=180, out=-45, looseness=.5
            ] node{\(f_k\)} ([yshift=0.6ex]X.west)
        (X) edge node{\(s\)} (Y)
        ;
    \end{tikzpicture}
    \caption{}
    \label{d:codil-colim-2}
\end{minipage}
\end{diagram}
For each \(j \in J\), the morphism \(a_j\) has a codilator \((Y_j, s_j, t_j)\).

Let \(j \leq k\) in~\(J\). Then
\[
    {t_k}^\inv s_k f_{j \leq k}
    = a_k f_{j \leq k}
    = a_j.
\]
Hence \((Y_k, s_kf_{j \leq k}, t_k)\) is a codilation of \(a_j\). Universality of the codilator \((Y_j, s_j, t_j)\) of \(a_j\) yields a unique isometry \(g_{j \leq k} \colon Y_j \to Y_k\) such that
\[
    g_{j \leq k} s_j = s_k f_{j \leq k} \qquad\text{and}\qquad
    g_{j \leq k}t_j = t_k.
\]

Universality of codilators also implies that the morphisms \((g_{j \leq k} \colon Y_j \to Y_k)_{j,k \in J}\) form a directed diagram. By assumption, this diagram has a colimit \((g_j \colon Y_j \to Y)_{j \in J}\) in \(\AIsom{\D}\).
For all \(j \leq k\) in \(J\),
\[
    g_ks_kf_{j \leq k} = g_k g_{j \leq k} s_j = g_js_j
    \qquad\text{and}\qquad
    g_kt_k = g_k g_{j \leq k}t_j = g_j t_j.
\]
Hence \((g_js_j \colon X_j \to Y)_{j \in J}\) and \((g_j t_j \colon A \to Y)_{j \in J}\) are cocones on the directed diagrams \((f_{j \leq k} \colon X_j \to X_k)_{j,k \in J}\) and \((1 \colon A \to A)_{j,k \in J}\), respectively. The respective universal properties of the colimits \((f_j \colon X_j \to X)_{j \in J}\) and \((1 \colon A \to A)_{j \in J}\) yield unique isometries \(s \colon X \to Y\) and \(t \colon A \to Y\) satisfying
\[
    sf_j = g_j s_j
    \qquad\text{and}\qquad
    t = g_jt_j
\]
for each \(j \in J\).
Let \(a = t^\inv s\). Then \((Y, s, t)\) is a codilation of \(a\). Also, for each \(j \in J\),
\[
af_j
= t^\inv s f_j
= t^\inv g_j s_j
= {t_j}^\inv {g_j}^\inv g_j s_j
= {t_j}^\inv s_j
= a_j \text.
\]
This completes the proof of existence.

For uniqueness, let \(a' \colon X \to A\) be a morphism in \(\D\), and suppose that
\[a' f_j = a_j\]
for all \(j \in J\). By assumption, the morphism \(a'\) has a codilator \((Y',s',t')\). For each \(j \in J\), universality of the codilator \((Y_j, s_j, t_j)\) yields a unique isometry \(g_j' \colon Y_j \to Y'\) such that
\[
    g_j's_j = s'f_j \qquad\text{and}\qquad g_j' t_j = t'.
\]
In fact \((g_j' \colon Y_j \to Y')_{j \in J}\) is a cocone on \((g_{j \leq k} \colon Y_j \to Y_k)_{j, k \in J}\), also by universality of codilators. By universality of the colimit \((g_j \colon Y_j \to Y)_{j \in J}\), there is a unique isometry \(m \colon Y \to Y'\) such that
\[g'_j = m g_j\]
for each \(j \in J\), as depicted in \cref{d:dir-colim-isom-3}.
\begin{diagram}
    \centering
    \begin{tikzpicture}[x=7.5em, y=4.8em, commutative diagrams/.cd, every diagram]
        \begin{scope}[
            commutative diagrams/.cd,
            every node
        ]
            \node (Xj) at (0,2) {\(X_j\)};
            \node (Yj) at (0,1) {\(Y_j\)};
            \node (A) at (0,0) {\(A\)};
            \node (X) at (1,2) {\(X\)};
            \node (Y) at (1,1) {\(Y\)};
            \node (Y') at (1.5,0.65) {\(Y'\)};
            \node (Ar) at (1,0) {\(A\)};

        \end{scope}
    
        \path[
            commutative diagrams/.cd,
            every arrow,
            every label,
            mono
        ]
        (A)
            edge node {\(t_j\)} (Yj)
            edge node[swap] {\(1\)} (Ar)
        (Ar)
            edge node[pos=0.4] {\(t\)} (Y)
            edge[swap] node {\(t'\)} (Y')

        (Xj)
            edge node[swap] {\(s_j\)} (Yj)
            edge node {\(f_j\)} (X)

        (Yj)
            edge[out=-45, in=180, looseness=0.4,commutative diagrams/.cd, crossing over, pos=0.425] node[swap] {\(g_j'\)} (Y')    
            edge node {\(g_j\)} (Y)
        
        (X)
            edge node[swap] {\(s\)} (Y)
            edge[shorten <={0.5ex}] node {\(s'\)} (Y')
        
        (Y)
            edge[shorten <={-0.5ex}] node[pos=0.2,inner sep={0.3ex}]{\(m\)} (Y')
        
        ;
    \end{tikzpicture}
    \caption{}
    \label{d:dir-colim-isom-3}
\end{diagram}
There exists \(j \in J\), so
\[mt = m g_j t_j = g'_j t_j = t'.\]
Also, for each \(j \in J\),
\[msf_j = mg_j s_j = g'_j s_j = s' f_j.\]
Hence \(ms = s'\) by universality of the colimit \((f_j \colon X_j \to X)_{j \in J}\).
It follows that
\[a' = {t'}^\inv s' = t^\inv m^\inv m s = t^\inv s = a\text. \qedhere\]
\end{proof}

\begin{corollary}
\label{p:dir-colim-preserve-con}
Let \(\C\) be a pre\nobreakdash-Hilbert \(\inv\)\nobreakdash-category. If \(\AIsom{\C}\) has directed colimits, then the forgetful functor from \(\AIsom{\C}\) to \(\ACon{\C}\) preserves them.
\end{corollary}

\begin{proof}
Every morphism in \(\ACon{\C}\) has a codilator~\cite[Theorem~8.9]{di-meglio:r-star-cats}, so the result is a special case of \cref{p:dir-colim-preserve}.
\end{proof}

\subsection{Directed colimits of coisometries}
\label{s:dir-colim-coisom}

Next we consider directed colimits of coisometries. The goal of this subsection is to prove the following proposition.

\begin{proposition}\label{p:codir-lim-preserve}
    Let \(\C\) be a pre\nobreakdash-Hilbert \(\inv\)\nobreakdash-category. If \(\AIsom{\C}\) has directed colimits, then it also has codirected limits, and the forgetful functor \(U \colon \AIsom{\C} \to \ACon{\C}\) preserves them.
\end{proposition}

The key idea is to exploit the contravariance of orthonormal complements. An \textit{orthonormal complement} of a morphism \(x\) is an isometric kernel of \(x^\inv\). As \(\C\) is a pre\nobreakdash-Hilbert \(\inv\)\nobreakdash-category, every morphism has an isometric kernel, so every morphism has an orthonormal complement. For each object \(X\), the contravariant functor
\[
    {\perp} \colon \CommaCat{\ACon{\C}}{X} \to \CommaCat{\ACon{\C}}{X}
\]
maps each object \(x\) of \(\CommaCat{\ACon{\C}}{X}\) to a chosen orthonormal complement \(x^\perp\), and its action on morphisms is uniquely determined by universality of kernels. This contravariant functor is self-adjoint on the right (see~\cite[1.818]{freyd:allegories}, \cite[3.G]{freyd:abelian-categories}, or~\cite{isbell:functorial}). The components \(\eta_x \colon x \to x^{\perp\perp}\) of the unit of the adjunction are also uniquely determined by universality of kernels. An object \(x\) of \(\CommaCat{\ACon{\C}}{X}\) is an isometry if and only if \(\eta_x\) is unitary, so the full subcategory
\[
    \CommaCat{U}{X} \colon \CommaCat{\AIsom{\C}}{X} \hookrightarrow \CommaCat{\ACon{\C}}{X}
\]
is (essentially) the fixed point of the adjunction~\cite[Proposition~3.4.3]{borceux:handbook-categorical-algebra-1}. In particular, the self-adjunction restricts along the inclusion functor, yielding a self-duality
\[
    {\perp} \colon \CommaCat{\AIsom{\C}}{X} \to \CommaCat{\AIsom{\C}}{X}\text.
\]
This story is merely a different perspective on the \(\inv\)\nobreakdash-category specialisation of the well-known adjunction between kernels and cokernels (see~\cite[Proposition~56]{di-meglio:fcon}).

The following lemma is a specialisation of a standard result about colimits and slice categories. The proof of the part about colimit preservation, which does not appear in any of the standard textbooks on category theory, is included here for completeness.

\begin{lemma}\label{l:codir-lim-preserve:1}
    Let \(\C\) be a pre\nobreakdash-Hilbert \(\inv\)\nobreakdash-category. If \(\AIsom{\C}\) has directed colimits, then, for each object \(X\) of \(\C\), the comma category \(\CommaCat{\AIsom{\C}}{X}\) has directed colimits and the forgetful functor \(\CommaCat{U}{X} \colon \CommaCat{\AIsom{\C}}{X} \to \CommaCat{\ACon{\C}}{X}\) preserves them.
\end{lemma}
\begin{proof}
    Consider the following commutative diagram of functors.
    \[
    \begin{tikzcd}[column sep={between origins, 8em}]
        \CommaCat{\AIsom{\C}}{X} 
        \ar[d, "\Pi",swap]
        \ar[hook, r, "\CommaCat{U}{X}"]
        \& 
        \CommaCat{\ACon{\C}}{X} 
        \ar[d, "\Pi"]
        \\
        \AIsom{\C}
        \ar[r,hook,"U",swap]
        \& 
        \ACon{\C} 
    \end{tikzcd}
    \]
    For existence, the category \(\AIsom{\C}\) has directed colimits and the functor \(\Pi \colon \CommaCat{\AIsom{\C}}{X} \to \AIsom{\C}\) strictly creates colimits~\cite[Proposition~3.3.8]{riehl:category-theory-context}, so \(\CommaCat{\AIsom{\C}}{X}\) also has directed colimits.

    For preservation, let \(\alpha\) be a colimit of a directed diagram \(F\) in \(\CommaCat{\AIsom{\C}}{X}\). The functor \(\Pi \colon \CommaCat{\AIsom{\C}}{X} \to \AIsom{\C}\) strictly creates colimits, so it preserves them. Also \(U\) preserves directed colimits by~\cref{p:dir-colim-preserve-con}. Thus \(U \circ \Pi \circ \alpha\) is a colimit of \(\Pi \circ (\CommaCat{U}{X}) \circ F = U \circ \Pi \circ F\). But \(\Pi \colon \CommaCat{\ACon{\C}}{X} \to \ACon{\C}\) strictly creates colimits, so
    \begin{enumerate}
        \item there is a unique cocone \(\beta\) on \((\CommaCat{U}{X}) \circ F\) such that \(\Pi \circ \beta = U \circ \Pi \circ \alpha\), and
        \item the cocone \(\beta\) is a colimit.
    \end{enumerate}
    But \(\Pi \circ (\CommaCat{U}{X}) \circ \alpha = U \circ \Pi \circ \alpha\), so \(\beta = (\CommaCat{U}{X}) \circ \alpha\).
    Hence \((\CommaCat{U}{X}) \circ \alpha\) is a colimit.
\end{proof}

\begin{lemma}\label{l:codir-lim-preserve:2}
    Let \(\C\) be a pre\nobreakdash-Hilbert \(\inv\)\nobreakdash-category. If \(\AIsom{\C}\) has directed colimits, then, for each object \(X\) of \(\C\), the comma category \(\CommaCat{\AIsom{\C}}{X}\) has codirected limits and the forgetful functor \(\CommaCat{U}{X} \colon \CommaCat{\AIsom{\C}}{X} \to \CommaCat{\ACon{\C}}{X}\) preserves them.
\end{lemma}
\begin{proof}
    Consider the commutative diagram 
    \[
    \begin{tikzcd}[column sep={between origins, 8em}]
        \CommaCat{\AIsom{\C}}{X_{\phantom{1}}}
            \ar[d, "\perp", swap]
            \ar[r,hook,"\CommaCat{U}{X}", shorten <=-1ex]
        \&
        \CommaCat{\ACon{\C}}{X}_{\phantom{\leq 1}}
            \ar[d, "\perp"]
        \\
        \CommaCat{\AIsom{\C}}{X}_{\phantom{1}}
            \ar[shorten <=-1ex, hook, swap, r, "\CommaCat{U}{X}"]
        \& 
        \CommaCat{\ACon{\C}}{X}_{\phantom{\leq 1}}
    \end{tikzcd}
    \]
    of covariant and contravariant functors.

    For existence, the category \(\CommaCat{\AIsom{\C}}{X}\) has directed colimits and \({\perp} \colon \CommaCat{\AIsom{\C}}{X} \to \CommaCat{\AIsom{\C}}{X}\) is a self-duality, so \(\CommaCat{\AIsom{\C}}{X}\) also has codirected limits.

    For preservation, let \(\alpha\) be a limit of a codirected diagram \(F\) in \(\CommaCat{\AIsom{\C}}{X}\). As \({\perp} \colon \CommaCat{\AIsom{\C}}{X} \to \CommaCat{\AIsom{\C}}{X}\) is a self-duality, it maps limits to colimits.
    We already know that \(\CommaCat{U}{X}\) preserves directed colimits by~\cref{l:codir-lim-preserve:1}. As \({\perp} \colon \CommaCat{\ACon{\C}}{X} \to \CommaCat{\ACon{\C}}{X}\) is self-adjoint on the right, it maps colimits to limits. Hence the (covariant) composite
    \[{\perp} \circ (\CommaCat{U}{X}) \circ {\perp} = {\perp} \circ {\perp} \circ (\CommaCat{U}{X})\]
    preserves codirected limits. But the functor \({\perp} \colon \CommaCat{\AIsom{\C}}{X} \to \CommaCat{\AIsom{\C}}{X}\) is a self-duality, so its unit \(\eta \colon 1 \Rightarrow {\perp} \circ {\perp}\) is a natural isomorphism.
    It follows 
    that the functor \(\CommaCat{U}{X} \colon \CommaCat{\AIsom{\C}}{X} \to \CommaCat{\ACon{\C}}{X}\) is naturally isomorphic to \({\perp} \circ {\perp} \circ (\CommaCat{U}{X})\) and hence also preserves codirected limits.
\end{proof}

We are now ready for the third step, proving the main result in this subsection.

\begin{proof}[Proof of \cref{p:codir-lim-preserve}]
    Let \(\J\) be a codirected set regarded as a category, and let \(F\) be a \(\J\)\nobreakdash-shaped diagram in \(\AIsom{\C}\). As \(\J\) is codirected, it has an object \(J\). Consider the following commutative diagram of (covariant) functors.
    \[
        \begin{tikzcd}[column sep={between origins, 8em}]
            \CommaCat{\cat{J}}{J}
            \ar[r, "\CommaCat{F}{J}"]
            \ar[d, "\Pi",swap]
            \&
            \CommaCat{\AIsom{\C}}{FJ} 
            \ar[d, "\Pi",swap]
            \ar[hook, r, "\CommaCat{U}{FJ}"]
            \& 
            \CommaCat{\ACon{\C}}{FJ} 
            \ar[d, "\Pi"]
            \\
            \cat{J}
            \ar[r,"F",swap]
            \& 
            \AIsom{\C} 
            \ar[r,hook,"U",swap]
            \& 
            \ACon{\C} 
        \end{tikzcd}
    \]

    We begin with existence. As \(\J\) is a codirected set regarded as a category, so is the comma category \(\CommaCat{\J}{J}\). Now \(\CommaCat{\AIsom{\C}}{FJ}\) has codirected limits by~\cref{l:codir-lim-preserve:2}. Hence the codirected diagram \(\CommaCat{F}{J}\) has a limit \(\alpha\). The functor \(\Pi \colon \CommaCat{\AIsom{\C}}{FJ} \to \AIsom{\C}\) preserves codirected limits~\cite[Proposition~43]{vicary:conn-lims}, so \(\Pi \circ \alpha\) is a limit of the diagram \(\Pi \circ (\CommaCat{F}{J}) = F \circ \Pi\). Since \(\J\) is a codirected set, the comma category \(\CommaCat{\Pi}{K}\) is connected for each object \(K\) of \(\J\), so the functor \(\Pi \colon \CommaCat{\J}{J} \to \J\) is initial~\cite[Lemma~8.3.4]{riehl:homotopy}. Hence~\cite[Definition~8.3.2]{riehl:homotopy} the diagram \(F\) has a limit \(\beta\) such that \(\beta \circ \Pi = \Pi \circ \alpha\).

    For preservation, let \(\beta\) be a limit of \(F\). Since the functor \(\Pi \colon \CommaCat{\J}{J} \to \J\) is initial, \(\beta \circ \Pi\) is a limit of the diagram \(F \circ \Pi = \Pi \circ (\CommaCat{F}{J})\).
    As \(\Pi \colon \CommaCat{\AIsom{\C}}{FJ} \to \AIsom{\C}\) creates codirected limits~\cite[Proposition~3.3.8]{riehl:category-theory-context}, the diagram \(\CommaCat{F}{J}\) has a limit \(\alpha\) such that \(\Pi \circ \alpha = \beta \circ \Pi\).
    By~\cref{l:codir-lim-preserve:2}, the functor \(\CommaCat{U}{FJ}\) preserves codirected limits. The functor \(\Pi \colon \CommaCat{\ACon{\C}}{FJ} \to \ACon{\C}\) also preserves them~\cite[Proposition~43]{vicary:conn-lims}.
    Hence 
    \[
        \Pi \circ (\CommaCat{U}{FJ}) \circ \alpha = U \circ \Pi \circ \alpha = U \circ \beta \circ \Pi
    \]
    is a limit of 
    \[
        \Pi \circ (\CommaCat{U}{FJ}) \circ (\CommaCat{F}{J}) = U \circ F \circ \Pi\text.
    \]
    Using the fact that \(\Pi \colon \CommaCat{\J}{J} \to \J\) is an initial functor again, we conclude that \(U \circ \beta\) is a limit of \(U \circ F\).
\end{proof}

\subsection{Directed colimits of contractions}
\label{s:dir-colim-con}

Finally, in \Cref{p:colim-isometries} below, we combine \Cref{p:dir-colim-preserve-con,p:codir-lim-preserve} to prove that \(\ACon{\C}\) has directed colimits. Our construction is similar to \citeauthor{milius:colimits-categories-relations}'s construction of lax adjoint cooplimits of \(\omega\)\nobreakdash-chains in certain 2\nobreakdash-categories of relations~\cite[Theorem~4.2]{milius:colimits-categories-relations}.

\begin{lemma}
\label{p:dir_limit_comm_coisom}
Let \(\C\) be a pre\nobreakdash-Hilbert \(\inv\)\nobreakdash-category, and suppose that \(\AIsom{\C}\) has directed colimits. Let \((f_{j \leq j'} \colon X_j \to X_{j'})_{j,j'\in J}\)
and
\((f'_{j \leq j'} \colon X'_j \to X'_{j'})_{j,j'\in J}\) be directed diagrams in~\(\AIsom{\C}\), and let \((f_j \colon X_j \to X)_{j \in J}\) and \((f'_j \colon X'_j \to X')_{j \in J}\) be their respective colimits. Let \((g_j \colon X_j \to X'_{j})_{j \in J}\) be a natural transformation in \(\ACon{\C}\) between these diagrams, and let \(g \colon X \to X'\) be the induced contraction between their colimits (see \cref{p:dir-colim-preserve-con}). If each of the morphisms \(g_j\) is coisometric, then the morphism \(g\) is also coisometric.
\end{lemma}

\begin{proof}
Throughout this proof, refer to \cref{d:dir-lim-comm-coisom}.
\begin{diagram}
\begin{tikzpicture}[x=7.5em, y=4.8em, commutative diagrams/.cd, every diagram]
    \begin{scope}[commutative diagrams/.cd, every node]
        \node (Xj) at (0,2) {\(K_j\)};
        \node (Xk) at (1,2) {\(K_{j'}\)};
        \node (X) at (2,1.65) {\(K\)};
        \node (Yj) at (0,1) {\(X_j\)};
        \node (Yk) at (1,1) {\(X_{j'}\)};
        \node (Y) at (2,0.65) {\(X\)};
        \node (Al) at (0,0) {\(X_j'\)};
        \node (Am) at (1,0) {\(X'_{j'}\)};
        \node (Ar) at (2,-0.35) {\(X'\)};
    \end{scope}
    
    \path[
        commutative diagrams/.cd,
        every arrow,
        every label,
        mono
    ]
    (Al)
        edge node{\(f'_{j \leq j'}\)} (Am)
        edge[commutative diagrams/.cd, twoheadleftarrow] node{\(g_j\)} (Yj)
        edge[
            in=180, out=-45, looseness=.4
        ] node[swap]{\(f'_j\)} ([yshift=-0.8ex]Ar.west)
    (Am)
        edge[commutative diagrams/.cd, twoheadleftarrow] node[pos=.3, swap]{\(g_{j'}\)} (Yk)
        edge[
            in=180, out=-45, looseness=.45
        ] node{\(f'_{j'}\)} ([yshift=0.6ex]Ar.west)
    (Ar)
        edge[commutative diagrams/.cd, twoheadleftarrow] node[swap]{\(g\)} (Y)
    (Yj)
        edge node{\(f_{j \leq j'}\)} (Yk)
        edge[
            in=180, out=-45, looseness=.4,
            commutative diagrams/.cd, crossing over
        ] node[swap,pos=.325]{\(f_j\)} ([yshift=-0.8ex]Y.west)
    (Yk)
        edge[
            in=180, out=-45, looseness=.45
        ] node[pos=.4]{\(f_{j'}\)} ([yshift=0.6ex]Y.west)
    (Xk)
        edge node[pos=.7]{\(m_{j'}\)} (Yk)
    (Xj)
        edge node[swap]{\(m_j\)} (Yj)
        edge node{\(k_{j \leq j'}\)} (Xk)
        edge[
            in=180, out=-45, looseness=.4,
            commutative diagrams/.cd, crossing over
        ] node[pos=.325,swap]{\(k_j\)} ([yshift=-0.8ex]X.west)
    (Xk)
        edge[
            in=180, out=-45, looseness=.45
        ] node{\(k_{j'}\)} ([yshift=0.6ex]X.west)
    (X) edge node{\(m\)} (Y)
    ;
\end{tikzpicture}
\caption{}
\label{d:dir-lim-comm-coisom}
\end{diagram}
For each \(j \in J\), the morphism \(g_j\) has an isometric kernel \(m_j \colon K_j \to X_j\) in \(\C\). For all \(j \leq j'\) in \(J\),
\[
    g_{j'}f_{j \leq j'}m_j = f'_{j \leq j'}g_j m_j = f'_{j \leq j'}0 = 0,
\]
so universality of the kernel \(m_{j'}\) yields a unique morphism \(k_{j \leq j'} \colon K_j \to K_{j'}\) in \(\C\) such that \(m_{j'} k_{j \leq j'} = f_{j \leq j'}m_{j}\). Universality also ensures that \((k_{j \leq j'} \colon K_j \to K_{j'})_{j,j' \in J}\) is a directed diagram in \(\C\). It is actually a directed diagram in \(\AIsom{\C}\) because
\[
    {k_{j \leq j'}}^\inv k_{j \leq j'}
    = {k_{j \leq j'}}^\inv {m_{j'}}^\inv m_{j'}k_{j \leq j'}
    = {m_j}^\inv {f_{j \leq j'}}^\inv f_{j \leq j'}m_j
    = {m_j}^\inv m_j
    = 1
\]
for all \(j \leq j'\) in \(J\). By assumption, it has a colimit \((k_j \colon K_j \to K)_{j \in J}\) in \(\AIsom{\C}\); this cocone remains a colimit in \(\ACon{\C}\) by \cref{p:dir-colim-preserve-con}. Universality in \(\AIsom{\C}\) of the colimit yields a unique isometry \(m \colon K \to X\) such that \(mk_j = f_j m_j\).

Let \(j \in J\). As \(g_j\) is coisometric, it is a cokernel in \(\C\) of \(1 - {g_j}^\inv g_j\). But normal epimorphisms are always cokernels of their kernels, so it is also a cokernel of \(m_j\) in \(\C\). The forgetful functor \(\ACon{\C} \hookrightarrow \C\) creates coisometric coequalisers~\cite[Proposition~8.4]{di-meglio:r-star-cats}, so \(g_j\) remains a cokernel of \(m_j\) in \(\ACon{\C}\).

As colimits commute with colimits, \(g\) is a cokernel of \(m\) in \(\ACon{\C}\), and thus is coisometric~\cite[Proposition~8.4]{di-meglio:r-star-cats}.
\end{proof}

\begin{construction}
\label{c:colim_construction}
Read this construction in consultation with \cref{d:codilator_construction}.
\begin{diagram}
\centering
\begin{tikzcd}[row sep={3.2em,between origins}, column sep={3.4em,between origins}, labels={inner sep=-0.5ex}, cramped, cells={inner sep=0.3ex}]
    \&
    \&
    \&
    \&
    \&
X_{j}
    \arrow[dr, "g_{j \leq j'}"{inner sep=0.25ex}, epi]
    \&
    \&
    \&
    \&
    \&
\\
    \&
    \&
    \&
    \&
    \&
    \&
X_{j'}
    \arrow[dr, "g_{j' \leq k}"{inner sep=0.3ex}, epi]
\\
    \&
    \&
    \&
X_{j}^{k'}
    \arrow[dr, "g_{j \leq j'}^{k'}" swap, epi]
    \arrow[uurr, "f_{j}^{k'}"{inner sep=0ex}, mono, dashed]
    \&
    \&
    \&
    \&
X_{k}
    \arrow[dr, "g_{k \leq k'}"{inner sep=0.3ex}, epi]
\\
    \&
    \&
X_{j}^{k}
    \arrow[ur, "f_{j}^{k \leq k'}", mono]
    \arrow[dr, "g_{j \leq j'}^{k}" swap, epi]
    \&
    \&
X_{j'}^{k'}
    \arrow[uurr, "f_{j'}^{k'}"{inner sep=0ex}, mono, dashed]
    \arrow[dr, "g_{j' \leq k}^{k'}" swap, epi]
    \&
    \&
    \&
    \&
X_{k'}
    \arrow[ddrr, "g_{k'}"{inner sep=0.3ex}, epi, dashed]
\\
    \&
X_{j}^{j'}
    \arrow[dr, "g_{j \leq j'}^{j'}" swap, epi]
    \arrow[ur, "f_{j}^{j' \leq k}", mono]
    \&
    \&
X_{j'}^{k}
    \arrow[ur, "f_{j'}^{k \leq k'}", mono]
    \arrow[dr, "g_{j' \leq k}^{k}" swap, epi]
    \&
    \&
X_{k}^{k'}
    \arrow[dr, "g_{k \leq k'}^{k'}" swap, epi]
    \arrow[uurr, "f_{k}^{k'}"{inner sep=0ex}, mono, dashed]
\\
X_{j}^{j}
    \arrow[ur, "f_{j}^{j \leq j'}", mono]
    \&
    \&
X_{j'}^{j'}
    \arrow[ur, "f_{j'}^{j' \leq k}", mono]
    \&
    \&
X_{k}^{k}
    \arrow[ur, "f_{k}^{k \leq k'}", mono]
    \&
    \&
X_{k'}^{k'}
    \arrow[uurr, "f_{k'}^{k'}"{inner sep=0ex}, mono, dashed]
    \&
    \&
    \&
    \&
X
\end{tikzcd}
\caption{}
\label{d:codilator_construction}
\end{diagram}
Let \(\C\) be a pre\nobreakdash-Hilbert \(\inv\)\nobreakdash-category, and suppose that \(\AIsom{\C}\) has directed colimits. Fix a directed set \(J\). For each \(j \in J\), let
\((f_{j}^{k \leq k'} \colon X_{j}^{k} \to X_{j}^{k'})_{k,k' \in J_{\geq j}}\)
be a diagram of isometries indexed by the directed set \(J_{\geq j}\), and let
\((f_{j}^{k} \colon X_{j}^{k} \to X_j)_{k \in J_{\geq j}}\)
be its colimit in \(\AIsom{\C}\). 
For each \(k \in J\), let \((g_{j \leq j'}^{k} \colon X_{j}^{k} \to X_{j'}^{k})_{j,j' \in J_{\leq k}}\)
be a diagram of coisometries indexed by the directed subset \(J_{\leq k} = \setb{j \in J}{j \leq k}\) of \(J\). Suppose that the square
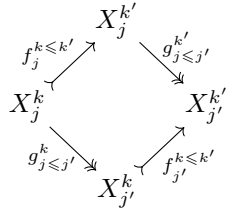
\[
    \begin{tikzcd}[row sep={3.2em,between origins}, column sep={3.4em,between origins}, labels={inner sep=-0.5ex}, cramped, cells={inner sep=0.3ex}]
        \&
    X_{j}^{k'}
        \arrow[dr, "g_{j \leq j'}^{k'}"{inner sep=0ex}, epi]
        \&
    \\
    X_{j}^{k}
        \arrow[ur, "f_{j}^{k \leq k'}", mono]
        \arrow[dr, "g_{j \leq j'}^{k}" swap, epi]
        \&
        \&
    X_{j'}^{k'}
    \\
        \&
    X_{j'}^{k}
        \arrow[ur, "f_{j'}^{k \leq k'}"{swap, inner sep=0ex}, mono]
        \&
    \end{tikzcd}
\]
is commutative for all \(j \leq j' \leq k \leq k'\) in \(J\).
\Cref{p:dir-colim-preserve-con} yields a unique contraction
\(g_{j \leq j'} \colon X_j \to X_{j'}\)
such that the square
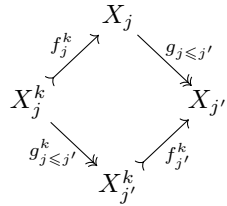
\[
    \begin{tikzcd}[row sep={3.2em,between origins}, column sep={3.4em,between origins}, labels={inner sep=-0.5ex}, cramped, cells={inner sep=0.3ex}]
        \&
    X_{j}
        \arrow[dr, "g_{j \leq j'}"{inner sep=0ex}, epi]
        \&
    \\
    X_{j}^{k}
        \arrow[ur, "f_{j}^{k}"{inner sep=0ex}, mono]
        \arrow[dr, "g_{j \leq j'}^{k}" swap, epi]
        \&
        \&
    X_{j'}
    \\
        \&
    X_{j'}^{k}
        \arrow[ur, "f_{j'}^{k}"{inner sep=0ex,swap}, mono]
        \&
    \end{tikzcd}
\]
is commutative; also \(g_{j \leq j'}\) is coisometric by \cref{p:dir_limit_comm_coisom}. In fact, \((g_{j \leq j'} \colon X_j \to X_{j'})_{j,j' \in J}\) is a directed diagram of coisometries. By the dual of \cref{p:codir-lim-preserve}, it has a colimit \((g_j \colon X_j \to X)_{j \in J}\) in the wide subcategory of coisometries.
\end{construction}
    
\begin{lemma}
    \label{l:colim_construction}
    In \Cref{c:colim_construction}, the contractions \((g_j f_{j}^{j} \colon X_{j}^{j} \to X)_{j \in J}\) form a colimit in \(\ACon{\C}\) of the directed diagram \((g_{j \leq k}^{k} f_{j}^{j \leq k} \colon X_{j}^{j} \to X_{k}^{k})_{j,k \in J}\) of contractions. 
\end{lemma}

\begin{proof}
    First, \((g_j f_{j}^{j} \colon X_{j}^{j} \to X)_{j \in J}\) is a cocone on the diagram. Indeed,
    \[
        g_{k}f^{k}_{k}g^{k}_{j \leq k} f^{j \leq k}_{j}
        = g_{k} g_{j \leq k} f^{k}_{j} f^{j \leq k}_{j} = g_j f^{j}_{j}
    \]
    for all \(j \leq k\) in \(J\).

    Let \((s_{j} \colon X^{j}_{j} \to Y)_{j \in J}\) be another cocone of contractions. We must show that there is a unique morphism \(t \colon X \to Y\) such that \(tg_jf^j_j = s_j\) for all \(j \in J\).
    
    For each \(j \in J\), the contractions \((s_k g^{k}_{j \leq k} \colon X_{j}^k \to Y)_{k \in J_{\geq j}}\) form a cocone on the diagram \((f_{j}^{k \leq k'} \colon X_{j}^k \to X_{j}^{k'})_{k,k' \in J_{\geq j}}\) of isometries; indeed, for all \(k' \geq k\) in \(J_{\geq j}\),
    \[
        s_k g^{k}_{j \leq k}
        = s_{k'}g^{k'}_{k \leq k'}f^{k \leq k'}_{k}g^{k}_{j \leq k}
        = s_{k'}g^{k'}_{k \leq k'} g^{k'}_{j \leq k} f_{j}^{k \leq k'}
        = s_{k'}g^{k'}_{j\leq k'} f_{j}^{k \leq k'}.
    \]
    Hence, for each \(j \in J\), there is a unique contraction \(t_j \colon X_j \to Y\) such that
    \[t_j f^{k}_j = s_k g^{k}_{j \leq k}\]
    for all \(k \in J_{\geq j}\), by \cref{p:dir-colim-preserve-con}.
    
    Let us show that the contractions \((t_j \colon X_j \to Y)_{j \in J}\) form a cocone on the directed diagram \((g_{j \leq j'} \colon X_{j} \to X_{j'})_{j,j' \in J}\) of coisometries. Let \(j' \geq j\) in \(J\). For all \(k \in J_{\geq j'}\),
    \[t_{j'}g_{j \leq j'} f^{k}_{j} = t_{j'}f^{k}_{j'} g^k_{j \leq j'} = s_{k}g^{k}_{j' \leq k} g^k_{j \leq j'} = s_{k} g^{k}_{j \leq k}.\]
    The right-hand side is independent of \(j'\). Also \((f^k_j \colon X^k_j \to X_j)_{k \in J_{\geq j'}}\) is the tail of a colimit in \(\AIsom{\C}\), so it is itself a colimit in \(\AIsom{\C}\). Hence
    \(t_{j'} g_{j \leq j'} = t_j g_{j \leq j} = t_j\) by \cref{p:dir-colim-preserve-con},
    as claimed.
    
    By \cref{p:codir-lim-preserve}, there is a unique contraction \(t \colon X \to Y\) such that \(t g_j = t_j\) for all \(j \in J\). As desired,
    \(tg_j f_j^{j} = t_j f_j^j = s_j g^{j}_{j \leq j} = s_j\) for all \(j \in J\).
    
    For uniqueness, let \(t' \colon X \to Y\) be another contraction satisfying \(t'g_j f_j^{j} = s_j\) for all \(j \in J\). For each \(j \in J\),
    \[t'g_j f^{k}_{j} = t'g_k g_{j \leq k} f^{k}_{j} = t'g_k f^{k}_{k}g^k_{j \leq k} = s_k g^k_{j \leq k}\]
    for each \(k \in J_{\geq j}\), so \(t'g_j = t_j\) by \cref{p:dir-colim-preserve-con}. Hence \(t = t'\) by \cref{p:codir-lim-preserve}.
\end{proof}

\begin{proposition}
    \label{p:colim-isometries}
    Let \(\C\) be a pre\nobreakdash-Hilbert \(\inv\)\nobreakdash-category. If \(\AIsom{\C}\) has directed colimits, then so does \(\ACon{\C}\).
\end{proposition}

To obtain a colimit of a directed diagram in \(\ACon{\C}\) from \cref{l:colim_construction}, we need to be able to choose codilations of all of the morphisms in the diagram in a compositional way. If the diagram is a \textit{smooth chain}, then we can do this by transfinite recursion. A \textit{chain} is a directed diagram indexed by an ordinal. A chain \((h_{j \leq k} \colon Z_j \to Z_k)_{j,k \in \lambda}\) is \textit{smooth} if, for each limit ordinal \(\mu < \lambda\), the morphisms \((h_{j \leq \mu} \colon Z_j \to Z_\mu)_{j \in \mu}\) form a colimit of the restricted diagram \((h_{j \leq k} \colon Z_j \to Z_k)_{j,k \in \mu}\). A category has directed colimits if and only if it has colimits of smooth chains~\cite[Corollary~1.7]{adamek:locally-presentable-accessible}.

\begin{proof}
As discussed above, it suffices to show that \(\ACon{\C}\) has colimits of smooth chains.

Let \((h_{j \leq k} \colon Z_j \to Z_k)_{j,k \in \lambda}\) be a smooth chain in \(\ACon{\C}\). We will recursively define isometries \(f_j^{k \leq k'} \colon X^k_j \to X^{k'}_j\) and coisometries \(g^k_{j \leq j'} \colon X^k_j \to X^k_{j'}\) that satisfy the hypotheses of \cref{c:colim_construction}, such that
\[X^{j}_{j} = Z_j \qquad\text{and}\qquad g^k_{j \leq k} f^{j \leq k}_{j} = h_{j \leq k}\]
for all \(j \leq k\) in \(\lambda\). The result then follows from \cref{l:colim_construction}.

Let \(\mu < \lambda\). Suppose that \(f_{j}^{k \leq k'} \colon X^k_j \to X^{k'}_j\) and \(g^{k}_{j \leq j'} \colon X^{k}_{j} \to X^k_{j'}\) are already defined for all \(j \leq j' \leq k \leq k'\) in \(\mu\).

First, for the case where \(\mu = 0\), let
\[X_0^0 = Z_0, \qquad f_0^{0 \leq 0} = 1, \qquad\text{and}\qquad g^0_{0 \leq 0} = 1.\]

Next, suppose that \(\mu = \nu + 1\) for some ordinal \(\nu\). As \(h_{\nu \leq \mu} \colon Z_\nu \to Z_\mu\) is a contraction, there is an object \(A\) and a coisometry \(e \colon Z_\nu \oplus A \to Z_{\mu}\) such that \(h_{\nu \leq \mu} = e i_1\)~\cite[Theorem~8.9]{di-meglio:r-star-cats}. For all \(j \leq k \leq \nu\), let
\[
\begin{tikzcd}[row sep={3.6em,between origins}, column sep={3.8em,between origins}, labels={inner sep=-0.5ex}, cramped, cells={inner sep=0.5ex}]
    \&
X_j^\mu
    \arrow[dr, epi, "g^\mu_{j \leq k}"{inner sep=0ex}]
    \&
    \&
    \&
\\
|[color=Gray]| X^\nu_j
    \arrow[ur, mono, "f^{\nu \leq \mu}_j"]
    \arrow[dr, epi, "g^\nu_{j \leq k}" swap, Gray]
    \&
    \&
X^\mu_k
    \arrow[dr, "g^\mu_{k \leq \nu}"{inner sep=0}, epi]
    \&
    \&
\\
    \&
|[color=Gray]| X^\nu_k
    \arrow[ur, "f^{\nu \leq \mu}_k", mono]
    \arrow[dr, "g^\nu_{k \leq \nu}" swap, epi, Gray]
    \&
    \&
X^\mu_\nu
    \arrow[dr, "g^\mu_{\nu \leq \mu}"{inner sep=0ex}, epi]
    \&
\\
    \&
    \&
|[color=Gray]| X^\nu_\nu
    \arrow[ur, mono, "f^{\nu \leq \mu}_\nu"]
    \arrow[rr, "h_{\nu \leq \mu}"{inner sep=0.2ex, swap}, Gray]
    \&
    \&
X^\mu_\mu
\end{tikzcd}
\nqquad=\quad
\begin{tikzcd}[row sep={3.6em,between origins}, column sep={3.8em,between origins}, labels={inner sep=0ex}, cramped, cells={inner sep=0.5ex}]
    \&
X^\nu_j \oplus A
    \arrow[dr, "g^\nu_{j \leq k} \oplus 1", epi]
    \&
    \&
    \&
\\
|[color=Gray]| X^\nu_j
    \arrow[ur, mono, shorten >=0.7ex, "i_1"{inner sep=0.2ex}]
    \arrow[dr, epi, "g^\nu_{j \leq k}"{inner sep=-0.5ex, swap}, Gray]
    \&
    \&
X^\nu_k \oplus A
    \arrow[dr, epi, "g^\nu_{k \leq \nu} \oplus 1"]
    \&
    \&
\\
    \&
|[color=Gray]| X^\nu_k
    \arrow[ur, mono, shorten >=0.7ex, "i_1"{inner sep=0.2ex}]
    \arrow[dr, epi, "g^\nu_{k \leq \nu}"{inner sep=-0.5ex, swap}, Gray]
    \&
    \&
Z_\nu \oplus A
    \arrow[dr, epi, "e"{inner sep=0.3ex}]
    \&
\\
    \&
    \&
|[color=Gray]| Z_\nu
    \arrow[ur, mono, shorten >=0.7ex, "i_1"{inner sep=0.2ex}]
    \arrow[rr, "h_{\nu \leq \mu}"{inner sep=0.2ex, swap}, Gray]
    \&
    \&
Z_\mu
\end{tikzcd}
\nqquad,\qquad
\]
and then let
\[
    f^{k \leq \mu}_{j} = f^{\nu \leq \mu}_{j}f^{k \leq \nu}_{j}, \qquad
    f^{\mu \leq \mu}_{j} = 1,
    \qquad\text{and}\qquad
    g^{\mu}_{\mu \leq \mu} = 1.
\]

Finally, suppose that \(\mu\) is a limit ordinal. From \cref{c:colim_construction}, the morphisms 
\[
    f_j^k \colon X_j^k \to X_j,
    \qquad
    g_{j \leq k} \colon X_j \to X_k,
    \qquad\text{and}\qquad
    g_k \colon X_k \to X
\]
are defined for all \(j \leq k\) in \(\mu\). By \cref{l:colim_construction}, \((g_j f_j^j \colon Z_j \to X)_{j \in \mu}\) is a colimit in \(\ACon{\C}\) of the restricted diagram \((h_{j \leq k} \colon Z_j \to Z_k)_{j,k \in \mu}\).
But \((h_{j \leq k} \colon Z_j \to Z_k)_{j,k \in \lambda}\) is a smooth chain, so \((h_{j \leq \mu} \colon Z_j \to Z_\mu)_{j \in \mu}\) is another colimit in \(\ACon{\C}\) of the restricted diagram. Universality of these colimits yields an isomorphism \(u \colon X \to Z_\mu\) in \(\ACon{\C}\) such that \(ug_jf^j_j = h_{j \leq \mu}\) for all \(j \in \mu\). 
For all \(j \leq k \leq k'\) in \(\mu\), we may let
\[
\begin{tikzcd}[row sep={3.6em,between origins}, column sep={3.8em,between origins}, labels={inner sep=-0.5ex}, cramped, cells={inner sep=0.5ex}]
    \&
X_j^\mu
    \arrow[dr, epi, "g^\mu_{j \leq k}"{inner sep=0ex}]
    \&
    \&
    \&
\\
|[color=Gray]| X^{k'}_j
    \arrow[ur, mono, "f^{k' \leq \mu}_j"]
    \arrow[dr, epi, "g^{k'}_{j \leq k}" swap, Gray]
    \&
    \&
X^\mu_k
    \arrow[dr, "g^\mu_{k \leq k'}"{inner sep=0}, epi]
    \&
    \&
\\
    \&
|[color=Gray]| X^{k'}_k
    \arrow[ur, "f^{k' \leq \mu}_k", mono]
    \arrow[dr, "g^{k'}_{k \leq k'}" swap, epi, Gray]
    \&
    \&
X^\mu_{k'}
    \arrow[dr, "g^\mu_{k' \leq \mu}"{inner sep=0ex}, epi]
    \&
\\
    \&
    \&
|[color=Gray]| X^{k'}_{k'}
    \arrow[ur, mono, "f^{k' \leq \mu}_{k'}"]
    \arrow[rr, "h_{k' \leq \mu}"{inner sep=0.5ex, swap}, Gray]
    \&
    \&
X^\mu_\mu
\end{tikzcd}
\nqquad=\quad
\begin{tikzcd}[row sep={3.6em,between origins}, column sep={3.8em,between origins}, labels={inner sep=-0.5ex}, cramped, cells={inner sep=0.5ex}]
    \&
X_j
    \arrow[dr, "g_{j \leq k}"{inner sep=0.2ex}, epi]
    \&
    \&
    \&
\\
|[color=Gray]| X^{k'}_j
    \arrow[ur, mono, "f^{k'}_j"{inner sep=-0.2ex}]
    \arrow[dr, epi, "g^{k'}_{j \leq k}" swap, Gray]
    \&
    \&
X_k
    \arrow[dr, epi, "g_{k \leq k'}"{inner sep=0.2ex}]
    \&
    \&
\\
    \&
|[color=Gray]| X^{k'}_k
    \arrow[ur, mono, "f^{k'}_k"{inner sep=-0.2ex}]
    \arrow[dr, epi, "g^{k'}_{k \leq k'}" swap, Gray]
    \&
    \&
X_{k'}
    \arrow[dr, epi, "ug_{k'}"{inner sep=0.2ex}]
    \&
\\
    \&
    \&
|[color=Gray]| Z_{k'}
    \arrow[ur, mono, "f_{k'}^{k'}"{inner sep=-0.2ex}]
    \arrow[rr, "h_{k' \leq \mu}"{swap, inner sep=0.5ex}, Gray]
    \&
    \&
Z_\mu
\end{tikzcd}
\nqquad.\qquad
\]
because \(u\) is unitary~\cite[Proposition~8.4~(i)]{di-meglio:r-star-cats}.
\end{proof}

This finishes the proof that (iii) implies (ii) in \cref{p:dir-colim-equivalence}.

\subsection{Contractions and \texorpdfstring{\(\ell^2\)}{l2}-limits}
\label{s:con-l2}

To complete our proof of \cref{p:dir-colim-equivalence}, it remains to verify that (ii) implies (i). We need the following lemma, whose proof is inspired by a proof of an analogous fact in a different setting~\cite[Proposition~19]{di-meglio:fcon}.

\begin{lemma}
\label{p:dir-colim-con-oplus}
Let \(\C\) be a pre\nobreakdash-Hilbert \(\inv\)\nobreakdash-category. If \(\ACon{\C}\) has directed colimits, then, for each object \(A\), the functor \(\blank \oplus A \colon \ACon{\C} \to \ACon{\C}\) preserves directed colimits.
\end{lemma}

\begin{proof}
Read this proof in conjunction with \cref{d:dir-colim-con-oplus}.
\begin{diagram}
    \centering
    \begin{tikzpicture}[x=7.5em, y=4.8em, commutative diagrams/.cd, every diagram]
        \begin{scope}[commutative diagrams/.cd, every node]
            \node (Xj) at (0,2) {\(X_j\)};
            \node (Xk) at (1,2) {\(X_k\)};
            \node (X) at (2,1.65) {\(X\)};
            \node (Yj) at (0,1) {\(X_j \oplus A\)};
            \node (Yk) at (1,1) {\(X_k \oplus A\)};
            \node (Y) at (2,0.65) {\(Y\)};
            \node (Al) at (0,0) {\(A\)};
            \node (Am) at (1,0) {\(A\)};
            \node (Ar) at (2,-0.35) {\(A\)};
        \end{scope}
        
        \path[
            commutative diagrams/.cd,
            every arrow,
            every label
        ]
        (Al)
            edge node{\(1\)} (Am)
            edge[commutative diagrams/.cd, twoheadleftarrow] node {\(p_2\)} (Yj)
            edge[
                in=180, out=-45, looseness=.4
            ] node[swap]{\(1\)} ([yshift=-0.8ex]Ar.west)
        (Am)
            edge[commutative diagrams/.cd, twoheadleftarrow] node[swap, pos=0.3] {\(p_2\)} (Yk)
        (Am)
            edge[
                in=180, out=-45, looseness=.5
            ] node{\(1\)} ([yshift=0.6ex]Ar.west)
        (Ar)
            edge[commutative diagrams/.cd, leftarrow] node[swap] {\(r_2\)} (Y)
        (Yj)
            edge node{\(f_{j \leq k} \oplus 1\)} (Yk)
            edge[
                in=180, out=-45, looseness=.4,
                commutative diagrams/.cd, crossing over
            ] node[swap,pos=.325]{\(g_j\)} ([yshift=-0.8ex]Y.west)
        (Yk)
            edge[
                in=180, out=-45, looseness=.5
            ] node[pos=.4]{\(g_k\)} ([yshift=0.6ex]Y.west)
        (Xk)
            edge[transform canvas={xshift=-0.75ex}, commutative diagrams/.cd, mono] node[pos=.7, swap] {\(i_1\)} (Yk)
            edge[transform canvas={xshift=0.75ex}, commutative diagrams/.cd, twoheadleftarrow] node[pos=.7]{\(p_1\)} (Yk)
        (Xj)
            edge[transform canvas={xshift=-0.75ex}, commutative diagrams/.cd, mono] node[pos=.6, swap] {\(i_1\)} (Yj)
            edge[transform canvas={xshift=0.75ex}, commutative diagrams/.cd, twoheadleftarrow] node[pos=.6]{\(p_1\)} (Yj)
            edge node{\(f_{j \leq k}\)} (Xk)
            edge[
                in=180, out=-45, looseness=.4,
                commutative diagrams/.cd, crossing over
            ] node[pos=.325,swap]{\(f_j\)} ([yshift=-0.8ex]X.west)
        (Xk)
            edge[
                in=180, out=-45, looseness=.5
            ] node{\(f_k\)} ([yshift=0.6ex]X.west)
        (X) 
            edge[transform canvas={xshift=-0.75ex}] node[swap] {\(s_1\)} (Y)
            edge[transform canvas={xshift=0.75ex}, commutative diagrams/.cd, leftarrow] node{\(r_1\)} (Y)
        ;
    \end{tikzpicture}
    \caption{}
    \label{d:dir-colim-con-oplus}
\end{diagram}

Let \((f_j \colon X_j \to X)_{j \in J}\) be a colimit of a directed diagram \((f_{j \leq k} \colon X_j \to X_k)_{j,k \in J}\) in \(\ACon{\C}\). By assumption, the directed diagram \((f_{j \leq k} \oplus 1 \colon X_j \oplus A \to X_k \oplus A)_{j,k \in J}\) in \(\ACon{\C}\) has a colimit \((g_j \colon X_j \oplus A \to Y)_{j \in J}\).

Also \((f_j p_1 \colon X_j \oplus A \to X)_{j \in J}\) and \((p_2 \colon X_j \oplus A \to A)_{j \in J}\) are cocones on this diagram. Indeed, for all \(j \leq k\) in \(J\),
\[
    f_k p_1(f_{j \leq k} \oplus 1) = f_k f_{j \leq k}p_1 = f_jp_1
    \qquad\text{and}\qquad
    p_2(f_{j \leq k} \oplus 1) = p_2.
\]
Hence, universality of its colimit \((g_j \colon X_j \oplus A \to Y)_{j \in J}\) yields unique contractions \(r_1 \colon Y \to X\) and \(r_2 \colon Y \to A\) such that
\[r_1 g_j = f_j p_1 \qquad\text{and}\qquad r_2 g_j = p_2\]
for all \(j \in J\).

Furthermore, \((g_ji_1 \colon X_j \to Y)_{j \in J}\) is a cocone on \((f_{j \leq k} \colon X_j \to X_k)_{j,k \in J}\). Indeed,
\[
    g_ji_1 = g_k(f_{j \leq k} \oplus 1)i_1 = g_k i_1 f_{j \leq k}
\]
for all \(j \leq k\) in \(J\). Universality of the colimit \((f_j \colon X_j \to X)_{j \in J}\) of this diagram yields a unique contraction \(s_1 \colon X \to Y\) such that, for all \(j \in J\),
\[s_1f_j = g_ji_1.\]

Also, for all \(j \in J\),
\[
    r_1 s_1f_j = r_1 g_j i_1 = f_j p_1i_1 = f_j.
\]
Hence \(r_1 s_1 = 1\) by universality of the colimit \((f_j \colon X_j \to X)_{j \in J}\).

For each \(j \in J\), the morphism \(p_2 \colon X_j \oplus A \to A\) is a cokernel of \(i_1 \colon X_j \to X_j \oplus A\) in \(\C\), and thus~\cite[Proposition~8.4~(iii)]{di-meglio:r-star-cats} also in \(\ACon{\C}\). But colimits commute with colimits, so \(r_2\) is a cokernel of \(s_1\) in \(\ACon{\C}\), and thus~\cite[Proposition~8.4~(iii)]{di-meglio:r-star-cats} also in~\(\C\). Thus \(r_2\) and \(s_1\) form a split exact sequence in \(\C\). But \(\C\) is additive, so \((Y, r_1, r_2)\) is a product in \(\C\) (see, e.g., \cite[Lemma~20]{di-meglio:fcon} and \cite[Proposition~1.8.7]{borceux:handbook-categorical-algebra-2}),
and in fact an orthonormal product~\cite[Proposition~8.4~(iv)]{di-meglio:r-star-cats}.
\end{proof}

\begin{proposition}
\label{p:dir-colim-bounded}
Let \(\C\) be a pre\nobreakdash-Hilbert \(\inv\)\nobreakdash-category. If \(\ACon{\C}\) has codirected limits, then \(\C\) is a Hilbert \(\inv\)\nobreakdash-category. In particular, every codirected limit in \(\ACon{\C}\) is an \(\ell^2\)\nobreakdash-limit in \(\C\).
\end{proposition}

\begin{proof}
Read this proof in conjunction with \cref{d:dir-colim-bounded}.
\begin{diagram}
    \centering
    \begin{tikzpicture}[x=7.5em, y=4.8em, commutative diagrams/.cd, every diagram]
        \begin{scope}[commutative diagrams/.cd, every node]
            \node (X1) at (0,0) {\(X_j \oplus A\)};
            \node (X2) at (1,0) {\(X_k \oplus A\)};
            \node (X) at (2,0.5) {\(X \oplus A\)};
            \node (A) at (2,1.5) {\(A\)};
            \node (Y) at (2.75,1.5) {\(Y\)};
        \end{scope}
        
        \path[
            commutative diagrams/.cd,
            every arrow,
            every label
        ]
        (X1)
            edge[
                commutative diagrams/.cd, leftarrow
            ] node[
                swap
            ] {\(f_{j \leq k} \oplus 1\)} (X2)

            edge[
                in=180, out=70, looseness=0.9,
                commutative diagrams/.cd, crossing over,
                leftarrow
            ] node[
                pos=.25,
            ] {\(\pair{x_j}{1}\)} (A.west)

            edge[
                in=180, out=45, looseness=.5,
                commutative diagrams/.cd, crossing over,
                leftarrow
            ] node[
                pos=.325,
                swap,
                inner sep={0.4ex}
            ] {\(f_j \oplus 1\)} (X.west)

        (X2)

            edge[
                in=180, out=75, looseness=1,
                commutative diagrams/.cd, crossing over,
                leftarrow
            ] node[
                pos=.35,
                inner sep={0.4ex}
            ] {\(\pair{x_k}{1}\)} ([yshift=-0.8ex]A.west)

            edge[
                in=180, out=45, looseness=.8,
                commutative diagrams/.cd, leftarrow
            ] node[
                swap
            ] {\(f_k \oplus 1\)} ([yshift=-0.8ex]X.west)

        (X)
            edge[
                commutative diagrams/.cd, leftarrow
            ] node {\(\pair{x}{1}\)} (A)

            edge[
                commutative diagrams/.cd, leftarrow
            ] node[swap] {\(g\)} (Y)
        (A)
            edge[
                commutative diagrams/.cd, leftarrow
            ] node {\(y^{-1}\)} (Y)
        ;
    \end{tikzpicture}
    \caption{}
    \label{d:dir-colim-bounded}
\end{diagram}

Let \((f_{j \leq k} \colon X_k \to X_j)_{j,k \in J}\) be a codirected diagram in \(\ACon{\C}\), let \((f_j \colon X \to X_j)_{j \in J}\) be a limit of the diagram in \(\ACon{\C}\), and let \((x_j \colon A \to X_j)_{j \in J}\) be an \(\ell^2\)\nobreakdash-cone on the diagram in \(\C\). Then there exists \(b \in \SelfAdj{\C(A, A)}\) such that \({x_j}^\inv x_j \leq b\) for all \(j \in J\). Now \(b + 1 \sgt 0\) by symmetry~\cite[Proposition~6.5]{di-meglio:r-star-cats}. Hence~\cite[Proposition~6.3]{di-meglio:r-star-cats} there is an isomorphism \(y \colon A \to Y\) in \(\C\) such that \(y^\inv y = b + 1\).

For each \(j \in J\), the morphism \(\pair{x_j}{1} y^{-1} \colon Y \to X_j \oplus A\) is a contraction because
\[
    y^{-1\inv}\pairBig{x_j}{1}^\inv \pairBig{x_j}{1} y^{-1}
    = y^{-1\inv} ({x_j}^\inv x_j + 1) y^{-1}
    \leq y^{-1\inv} y^\inv yy^{-1}
    = 1.
\]
Also, for all \(j \leq k\) in \(J\),
\[
    (f_{j \leq k} \oplus 1) \pairBig{x_k}{1} y^{-1}
    = \pairBig{f_{j \leq k} x_k}{1} y^{-1}
    = \pairBig{x_j}{1} y^{-1}.
\]
Hence \(\paren[\big]{\pair{x_j}{1} y^{-1} \colon Y \to X_j \oplus A}_{j \in J}\) is a cone on \((f_{j \leq k} \oplus 1 \colon X_k \oplus A \to X_j \oplus A)_{j,k \in J}\). But \((f_j \oplus 1 \colon X \oplus A \to X_j \oplus A)_{j \in J}\) is a limit in \(\ACon{\C}\) of this diagram, by the dual of \cref{p:dir-colim-con-oplus}. Hence there is a unique contraction \(g \colon Y \to X \oplus A\) such that
\[(f_j \oplus 1) g = \pairBig{x_j}{1} y^{-1}\]
for all \(j \in J\).

Let \(x = p_1 gy\). Then
\[gy = \pairBig{x}{1}.\]
Indeed, as \(J\) is directed, it has an element \(l\), and so
\[p_2gy = p_2(f_l \oplus 1)gy = p_2 \pairBig{x_l}{1} = 1.\]
Hence
\[f_j x = p_1 \pairBig{f_j x}{1} = p_1 (f_j \oplus 1) gy = p_1 \pairBig{x_j}{1} y^{-1} y = x_j\]
for all \(j \in J\). This equation uniquely determines \(x\) because the morphisms \(f_j\), when viewed as morphisms in \(\C\), remain jointly monic~\cite[Proposition~8.4~(ii)]{di-meglio:r-star-cats}.

Also, for each \(j \in J\),
\[
    {x_j}^\inv x_j
    = x^\inv {f_j}^\inv f_j x
    \leq x^\inv x.
\]
Similarly,
\[
    x^\inv x + 1
    = \pairBig{x}{1}^\inv \pairBig{x}{1}
    = y^\inv g^\inv g y
    \leq y^\inv y
    = b + 1,
\]
so \(x^\inv x \leq b\). Since the same argument works for every upper bound \(b\), it follows that \[x^\inv x = \sup_{j \in J} {x_j}^\inv x_j.\qedhere\]
\end{proof}

This completes the proof of \cref{p:dir-colim-equivalence}.

\section{More examples of Hilbert \texorpdfstring{\(*\)}{*}-categories}
\label{s:examples}

\subsection{Unitary representations}
\label{s:unitary-representations}

If \(\C\) is a \(\inv\)\nobreakdash-category and \(\D\) is a pre\nobreakdash-Hilbert \(\inv\)\nobreakdash-category, then the \(\inv\)\nobreakdash-category \([\C, \D]_\inv\) of \(\inv\)\nobreakdash-functors \(\C \to \D\) and natural transformations is again a pre\nobreakdash-Hilbert \(\inv\)\nobreakdash-category~\cite[Proposition~2.5]{di-meglio:r-star-cats}. As a special case, the \(\inv\)\nobreakdash-category \(\URep_{G}\) of unitary representations of a groupoid \(G\) is a pre\nobreakdash-Hilbert \(\inv\)\nobreakdash-category~\cite[Corollary~2.6]{di-meglio:r-star-cats}. Indeed, thinking of groupoids as small \(\inv\)\nobreakdash-categories in which every morphism is unitary, \textit{unitary representations} of \(G\) are precisely \(\inv\)\nobreakdash-functors \(G \to \Hilb\) and \textit{intertwiners} are precisely natural transformations, that is, \(\URep_{G} = [G, \Hilb]_\inv\).

\begin{proposition}
\label{p:unitary-representations}
For all groupoids \(G\), the \(\inv\)\nobreakdash-category \(\URep_{G}\) is a Hilbert \(\inv\)\nobreakdash-category.
\end{proposition}

\begin{proof}
    As \(\URep_G\) is a pre\nobreakdash-Hilbert \(\inv\)\nobreakdash-category~\cite[Corollary~2.6]{di-meglio:r-star-cats}, it suffices to show that its wide subcategory \(\URep_{G,1}\) of isometries has directed colimits. Consider a directed diagram \((f_{j \leq k} \colon F_j \Rightarrow F_k)_{j,k \in J}\) in \(\URep_{G,1}\).
    
    The isometries in \(\URep_{G,1}\) are precisely the natural transformations whose components are isometries in \(\Hilb\). Hence, for all objects \(X\) of \(G\), the morphisms \((f_{j \leq k, X} \colon F_j X \to F_k X)_{j,k \in J}\) form a directed diagram in \(\AIsom{\Hilb}\),
    which has a colimit \((f_{j, X} \colon F_j X \to FX)_{j \in J}\) because \(\Hilb\) is a Hilbert \(\inv\)\nobreakdash-category (see \cref{p:hilb}).
    
    Consider a morphism \(g \colon X \to Y\) in \(G\). As \(g\) is unitary, so is \(F_j g \colon F_j X \to F_j Y\) for each \(j \in J\). For all \(j \leq k\) in \(J\), the square
    \[
    \begin{tikzcd}[sep=large]
        F_jX
            \ar[r, "f_{j \leq k, X}"]
            \ar[d, "F_j g" swap] 
            \&
        F_k X
            \ar[d, "F_k g"]
        \\
        F_jY
            \ar[r, "f_{j \leq k, Y}" swap] 
            \&
        F_k Y
    \end{tikzcd}
    \]
    in \(\AIsom{\Hilb}\) is commutative. Hence, there is a unique isometry \(Fg\) such that the square
    \[
    \begin{tikzcd}[sep=large]
        F_j X
            \ar[d, "F_j g" swap]
            \ar[r, "f_{j, X}"]
            \&
        F X
            \ar[d, "Fg"]
        \\
        F_j Y
            \ar[r, "f_{j,Y}" swap]
            \&
        FY
    \end{tikzcd}
    \]
    in \(\AIsom{\Hilb}\) is commutative for all \(j \in J\).
    
    By universality of directed colimits, \(F\) is a \(\inv\)-functor \(G \to \Hilb\) and the original diagram in \(\URep_{G,1}\) has colimit \((f_j \colon F_j \Rightarrow F)_{j \in J}\).
\end{proof}

More generally, for each groupoid \(G\) and each Hilbert \(\inv\)\nobreakdash-category \(\D\), the \(\inv\)\nobreakdash-category \([G, \D]_\inv\) is again a Hilbert \(\inv\)\nobreakdash-category; the proof of \cref{p:unitary-representations} works \textit{mutatis mutandis}. The condition that every morphism in the \(\inv\)\nobreakdash-category \(G\) be unitary can also be relaxed somewhat. For example, it suffices that every morphism in \(G\) be the composite of a sequence of isometries and coisometries. However, as shown in the following proposition, such a condition on \(G\) cannot be dropped altogether.

Regard the additive monoid \(\Nats\) of natural numbers as a one-object \(\inv\)\nobreakdash-category with the trivial involution.

\begin{proposition}
The \(\inv\)\nobreakdash-category \([\Nats, \Hilb]_\inv\) is not a Hilbert \(\inv\)\nobreakdash-category.
\end{proposition}

To prove this proposition, it will be helpful to have a simpler description of the objects and morphisms of the \(\inv\)\nobreakdash-category \([\Nats, \Hilb]_\inv\). Its objects are completely determined by their value at the single object of \(\Nats\) and at the generating morphism of~\(\Nats\) (the morphism corresponding to the natural number \(1\)). In fact, its objects are in bijection with the pairs \((X, s)\) consisting of a Hilbert space \(X\) and a Hermitian map \(s \colon X \to X\). Under this bijection, the morphisms \((X, s) \to (Y, t)\) are the bounded linear maps \(f \colon X \to Y\) such that \(tf = fs\).

\begin{proof}
For each \(n \in \Nats\), define \(s_n \colon \Comps^n \to \Comps^n\) by
\[s_n = \begin{bmatrix} 1 & 0 & \cdots & 0 \\ 0 & 2 & \cdots & 0 \\ \vdots & \smash{\vdots} & \smash{\ddots} & \smash{\vdots} \\ 0 & 0 & \cdots & n \end{bmatrix}.\]
Assume, for a contradiction, that the directed diagram
\[
    \begin{tikzcd}[column sep=large]
        (\Comps^1, s_1)
            \arrow[r, "{\bsmallmat{ 1 \\ 0}}"]
            \&
        (\Comps^2, s_2)
            \arrow[r, "{\bsmallmat{ 1 & 0 \\0 &  1 \\ 0 & 0}}"]
            \&
        (\Comps^3, s_3)
            \arrow[r]
            \&[-2em]
        \cdots
    \end{tikzcd}
\]
in \([\Nats, \Hilb]_{\inv,1}\) has a colimit \(f_n \colon (\Comps^n, s_n) \to (X, s)\). Then, for each \(n \in \Nats\), as \(f_n\) is isometric, we have
\[
    \norm{f_n e_n} = \norm{e_n} = 1
\]
and
\[
    \norm{sf_ne_n} = \norm{f_ns_ne_n} = \norm{s_n e_n} = n,
\]
where \(e_n\) is the \(n\)\nobreakdash-th standard basis vector of \(\Comps^n\). Hence 
\[
    \norm{s} = \sup_{\substack{x \in X\\\norm{x} = 1}} \norm{sx} \geq \norm{sf_ne_n} = n
\]
for all \(n \in \Nats\), which is a contradiction.
\end{proof}

\subsection{Hilbert W*-modules}
\label{s:hilbert-w-modules}

Hilbert C*-modules over a complex C*\nobreakdash-algebra generalise Hilbert spaces over the complex numbers. Unfortunately, the \(\inv\)\nobreakdash-category of Hilbert C*-modules and adjointable maps over a fixed C*\nobreakdash-algebra is not necessarily even a pre\nobreakdash-Hilbert \(\inv\)\nobreakdash-category~\cite[Remark~5.9]{di-meglio:r-star-cats}, let alone a Hilbert \(\inv\)\nobreakdash-category. To obtain a Hilbert \(\inv\)\nobreakdash-category, one must look to well-behaved Hilbert C*-modules over well-behaved C*\nobreakdash-algebras. In this subsection, we show that over a W*\nobreakdash-algebra \(A\), the Hilbert W*\nobreakdash-modules and bounded \(A\)-linear maps form a Hilbert \(\inv\)\nobreakdash-category~\(\Hilb_A\). The proof is an adaptation of the proof in \cref{s:hilbert-spaces} that the \(\inv\)\nobreakdash-category \(\Hilb = \Hilb_\Comps\) of complex Hilbert spaces is a Hilbert \(\inv\)\nobreakdash-category.

\subsubsection{W*-algebras}

A (\textit{complex unital}) \textit{C*\nobreakdash-algebra} is a \(\inv\)\nobreakdash-ring \(A\) equipped with a \(\inv\)\nobreakdash-ring homomorphism from \(\Comps\) to the centre of \(A\) and a complete norm \(\norm{\blank} \colon A \to \PosReals\) such that \(\norm{ab} \leq \norm{a}\norm{b}\) and \(\norm{a^\inv a} = \norm{a}^2\) for all \(a, b \in A\).

The \textit{Hermitian} and \textit{skew-Hermitian parts} of an element \(a \in A\) are, respectively, the Hermitian and skew-Hermitian elements
\[\He a = \frac{1}{2}(a + a^\inv) \qquad\text{and}\qquad \Sk a = \frac{1}{2}(a - a^\inv)\]
of \(A\). Of course \(a = \He a + \Sk a\) and \(\He \mathrm{i}a = \mathrm{i} \Sk a\).
A C*\nobreakdash-algebra \(A\) is canonically an ordered \(\inv\)\nobreakdash-ring with \(a \leq b\) if and only if \(b - a = c^\inv c\) for some \(c \in A\). 

A net \((a_j)_{j \in J}\) in a C*\nobreakdash-algebra \(A\) \textit{order converges}\footnote{Order convergence was introduced by \textcite[310]{widom:embedding-algebras} for commutative C*\nobreakdash-algebras and by \textcite[206]{kadison:equivalence-operator-algebras} for arbitrary C*\nobreakdash-algebras. Our definition is equivalent to the one given by~\textcite[260]{hamana:tensor-products-monotone}, which itself is a slight generalisation of the original one given by \citeauthor{kadison:equivalence-operator-algebras}.} to an element \(a \in A\) if there is a decreasing net \((\epsilon_j)_{j \in J}\) in \(\SelfAdj{A}\) with infimum \(0\) such that
\[-\epsilon_j \leq \He(a - a_j) \leq \epsilon_j\qquad\text{and}\qquad -\epsilon_j \leq \mathrm{i}\Sk(a - a_j) \leq \epsilon_j\]
for all \(j \in J\). There is at most one such element \(a \in A\); it is called the \textit{order limit} of the net \((a_j)_{j \in J}\), and we denote it by
\[\olim_{j \in J} a_j.\]
An increasing net order converges if and only if it has a supremum, in which case its order limit and its supremum coincide. Similarly, a decreasing net order converges if and only if it has an infimum, in which case its order limit and its infimum coincide. In this way, order limits generalise suprema and infima.

Let \(A\) be a C*-algebra. A linear functional \(\omega \colon A \to \Comps\) is \textit{positive}~\cite[10~II]{westerbaan:2019:category} if \(a \geq 0\) implies \(\omega a \geq 0\). A positive linear functional \(\omega \colon A \to \Comps\) is \textit{normal}~\cite[42~II]{westerbaan:2019:category} if, for all increasing nets \((a_j)_{j \in J}\) in \(\SelfAdj{A}\) with a supremum,
    \[\omega \sup_{j \in J} a_j = \sup_{j \in J} \omega a_j.\]
If a linear functional on a C*\nobreakdash-algebra is positive and normal, then it preserves the involution~\cite[10~IV]{westerbaan:2019:category} and is bounded~\cite[20~II]{westerbaan:2019:category}.

A \textit{W*\nobreakdash-algebra} is a C*\nobreakdash-algebra with the following properties~\cite[42~I]{westerbaan:2019:category}:
\begin{itemize}
    \item \textit{monotone completeness}: every increasing net of Hermitian elements with an upper bound has a supremum, and
    \item \textit{faithfulness of the set of normal positive linear functionals}: for each \(a \geq 0\), if \(\omega a = 0\) for all normal positive linear functionals \(\omega\), then \(a = 0\).
\end{itemize}
Of course \(\Comps\) is itself a W*\nobreakdash-algebra.

\begin{lemma}
\label{l:positive-boundedness}
Let \(A\) be a W*-algebra, and let \(B \subseteq A_{\geq 0}\). If \(\omega(B)\) has an upper bound for each normal positive linear functional \(\omega\), then \(B\) has an upper bound.
\end{lemma}

\begin{proof}
Suppose that \(B\) does not have an upper bound. Then, for each \(n \geq 1\), there exists \(b_n \in B\) such that \(b_n \nleq 4^n\). The normal positive linear functionals are order separating~\cite[44~XI]{westerbaan:2019:category}, so, for each integer \(n \geq 1\), there is a normal positive linear functional \(\omega_n\) such that \(\omega_n(b_n) \nleq \omega_n(4^n)\), that is, such that \(\omega_n(b_n) > \omega_n(4^n)\); in particular, \(\omega_n\) is nonzero, so \(\norm{\omega_n} > 0\).

For each \(a \in A\), the series
\[
    \omega(a) = \sum_{n = 1}^\infty 2^{-n} \frac{\omega_n(a)}{\norm{\omega_n}}
\]
in \(\Comps\) is absolutely convergent. Indeed, for each integer \(N \geq 1\),
\[
    \sum_{n = 1}^N 2^{-n} \frac{\abs{\omega_n(a)}}{\norm{\omega_n}}
    \leq \sum_{n = 1}^N 2^{-n} \frac{\norm{\omega_n}\norm{a}}{\norm{\omega_n}}
    = \norm{a}\sum_{n = 1}^N 2^{-n}
    \leq \norm{a};
\]
since the sums on the left-hand side are increasing and bounded above, they have a supremum and converge to this supremum. The computation above also implies that \(\abs{\omega(a)} \leq \norm{a}\) for all \(a \in A\), so the function \(\omega \colon A \to \Comps\) is bounded. It is linear because each \(\omega_n\) is linear and limits in \(\Comps\) preserve addition and scalar multiplication. It is positive because each \(\omega_n\) is positive and limits in \(\Comps\) preserve positivity.

Next, we show that it is normal. Consider an increasing net \((a_j)_{j \in J}\) in \(\SelfAdj{A}\) with supremum \(a\). Since \(J\) is directed, there exists \(i \in J\). For each \(j \in J_{\geq i}\), let \(b_j = a_j - a_i\). Then \((b_j)_{j \in J_{\geq i}}\) is an increasing net in \(A_{\geq 0}\) with supremum \(b = a - a_i\). Now
\[
    \omega(b)
    = \sum_{n = 1}^\infty 2^{-n}\frac{\omega_n(b)}{\norm{\omega_n}}
    = \sum_{n = 1}^\infty \sup_{j \in J_{\geq i}} 2^{-n}\frac{\omega_n(b_j)}{\norm{\omega_n}}
    = \sup_{j \in J_{\geq i}}\sum_{n = 1}^\infty 2^{-n}\frac{\omega_n(b_j)}{\norm{\omega_n}}
    = \sup_{j \in J_{\geq i}} \omega(b_j)
\]
because the doubly-indexed net
\[
    \paren[\Big]{\sum_{n = 1}^N 2^{-n}\frac{\omega_n(b_j)}{\norm{\omega_n}}}_{N \geq 1, j \geq i}
\]
is increasing and bounded above by \(\omega(b)\). It follows that
\[
    \omega(a) - \omega(a_i)
    = \omega(b)
    = \sup_{j \in J_{\geq i}} \omega(b_j)
    = \sup_{j \in J_{\geq i}}\omega(a_j) - \omega(a_i).
\]
Adding \(\omega(a_i)\) and using cofinality of \(J_{\geq i}\) in \(J\), we obtain
\[
    \omega(a) = \sup_{j \in J} \omega(a_j).
\]

Finally, for all integers \(n \geq 1\),
\[
    \omega(b_n)
    \geq 2^{-n} \frac{\omega_n(b_n)}{\norm{\omega_n}}
    > 2^{-n}\frac{\omega_n(4^n)}{\norm{\omega_n}}
    = 2^n \frac{\omega_n(1)}{\norm{\omega_n}}
    \geq 2^n
\]
by~\cite[86~II]{westerbaan:2019:category}. Hence \(\omega(B)\) does not have an upper bound.
\end{proof}

\subsubsection{Hilbert W*-modules}

Let \(A\) be a W*\nobreakdash-algebra, and \(X\) be an inner product \(A\)\nobreakdash-module. For each \(x \in X\), the unique positive square root of \(\innerProd{x}{x}\) is called the \textit{absolute value} of \(x\) and denoted \(\abs{x}\). While it does not, in general, satisfy the triangle inequality, the following weaker inequality will suffice:
\begin{equation}
    \label{e:weak-tri-inequality}
    \abs{x + y}^2 \leq \abs{x + y}^2 + \abs{x - y}^2 = 2 \abs{x}^2 + 2 \abs{y}^2.
\end{equation}

For each positive linear functional \(\omega\) on \(A\), the equation
\[
    \norm{x}_{\omega} = \sqrt{\omega \abs{x}^2}
\]
defines~\cite[142~II]{westerbaan:dagger-dilation-category} a seminorm on \(X\), that is, it satisfies all of the norm axioms except possibly anisotropy. A net \((x_j)_{j \in J}\) in \(X\) \textit{ultranorm converges} to \(x \in X\) if
\[
    \lim_{j \in J} {\norm{x_j - x}_{\omega}} = 0
\]
for each \textit{normal} positive linear functional \(\omega\) on \(A\), that is, if it converges to \(x\) with respect to each of the seminorms \(\norm{\blank}_\omega\) where \(\omega\) is normal.\footnote{\label{f:uniformity} A net is ultranorm Cauchy or ultranorm convergent exactly when it is, respectively, Cauchy or convergent in the \textit{ultranorm uniformity}~\cite[146~VII]{westerbaan:dagger-dilation-category}. The topology generated by the ultranorm uniformity is known as the \textit{\(\tau_1\)\nobreakdash-topology}~\cite[Theorem~3.5.1]{manuilov:hilbert-w-modules-their} and the \textit{s\nobreakdash-topology}~\cite[Corollary~2.9]{ghez:w-categories}.} As the set of normal positive linear functionals on \(A\) is faithful, there is at most one such element \(x\); it will be referred to as the \textit{ultranorm limit} of \((x_j)_{j \in J}\), and denoted by
\[\unlim_{j \in J}x_j,\]
if it exists.\footnote{In contrast, limits with respect to a seminorm are not necessarily unique.} Similarly, a net \((x_j)_{j \in J}\) in \(X\) is \textit{ultranorm Cauchy} if, for each normal positive linear functional \(\omega\) on \(A\), it is Cauchy with respect to the seminorm \(\norm{\blank}_\omega\).
An inner product \(A\)\nobreakdash-module is \textit{ultranorm complete} if every ultranorm-Cauchy net ultranorm converges. A \textit{Hilbert W*\nobreakdash-module} over \(A\) is an inner product \(A\)\nobreakdash-module that is ultranorm complete.

\subsubsection{The category of Hilbert W*-modules}

Let \(A\) be a complex W*\nobreakdash-algebra. A function \(f \colon X \to Y\) between inner product \(A\)\nobreakdash-modules is \textit{adjointable} if there is a function \(g \colon Y \to X\) such that \(\innerProd{y}{fx} = \innerProd{gy}{x}\) for all \(x \in X\) and \(y \in Y\). There is at most one such function \(g\); it is called the \textit{adjoint} of \(f\), and denoted by \(f^\inv\), if it exists. Hilbert W*\nobreakdash-modules over \(A\) and adjointable maps form a \(\inv\)\nobreakdash-category \(\Hilb_A\).\footnote{The definition of \(\Hilb_A\) given here is equivalent to the one in~\cite[Section~5.3]{di-meglio:r-star-cats} because an inner product \(A\)\nobreakdash-module is ultranorm complete if and only if it is self-dual~\cite[149~V]{westerbaan:dagger-dilation-category}. Also, a function between self-dual inner product \(A\)-modules is bounded and \(A\)-linear if and only if it is adjointable (see~\cite[144~III and 152~VIII]{westerbaan:dagger-dilation-category} and \cite[35~VI]{westerbaan:2019:category}).}

\begin{proposition}
\label{p:hilb-mod}
The \(\inv\)\nobreakdash-category \(\Hilb_A\) is a Hilbert \(\inv\)\nobreakdash-category.
\end{proposition}

We already know that \(\Hilb_A\) is a pre\nobreakdash-Hilbert \(\inv\)\nobreakdash-category~\cite[Proposition~5.10]{di-meglio:r-star-cats}, so, to show that it is a Hilbert \(\inv\)\nobreakdash-category, it suffices, by \cref{p:dir-colim-equivalence}, to check that its wide subcategory \(\Hilb_{A,\leq 1}\) of contractions has directed colimits. A morphism \(f \colon X \to Y\) in \(\Hilb_A\) satisfies the relation \(\abs{fx}^2 \leq \abs{x}^2\) for all \(x \in X\) if and only if it is a contraction~\cite[Proposition~8.2]{di-meglio:r-star-cats}.

Let \((f_{j \leq k} \colon X_k \to X_j)_{j,k \in J}\) be a codirected diagram in \(\Hilb_{A,\leq 1}\).
The limit of the underlying diagram in the category \(\Mod_A\) of right \(A\)\nobreakdash-modules can be computed using the standard construction of limits from products and equalisers. Its apex is the submodule
\[\widebar{X} = \setb[\Big]{x \in \prod_{j \in J} X_j}{f_{j \leq k}x_k = x_j \text{ for all } j \leq k \text{ in } J}\]
of the product \(\prod_{j \in J} X_j\), and its projections are the restrictions \(\bar{f}_j \colon \widebar{X} \to X_j\) of the product projections \(p_j \colon \prod_{j \in J} X_j \to X_j\) to this submodule.

The limit \((f_j \colon X \to X_j)_{j \in J}\) in \(\Hilb_{A,\leq 1}\) of the diagram \((f_{j \leq k} \colon X_k \to X_j)_{j,k \in J}\) has the following concrete description. Its apex is the submodule
\[X = \setb[\Big]{x \in \widebar{X}}{\sup_{j \in J} \abs{x_j}^2 < \infty}\]
of \(\widebar{X}\), equipped with the \(A\)\nobreakdash-valued inner product defined by
\begin{equation}
    \label{e:dir-colim-inn-prod-mod}
    \innerProd{x}{y} = \olim_{j \in J} \innerProd{x_j}{y_j}
\end{equation}
for all \(x, y \in X\). Its projections are the restrictions \(f_j\) of the morphisms \(\bar{f}_j\) to this submodule. We prove this fact in several steps.

\begin{lemma}
The set \(X\) is a submodule of \(\widebar{X}\).
\end{lemma}

\begin{proof}
For all \(x, y \in X\), all \(a \in A\), and all \(k \in J\),
\[\abs{x_k a + y_k}^2 \leq 2\abs{x_k a}^2 + 2\abs{y_k}^2 = 2 a^\inv \abs{x_k}^2 a + 2 \abs{y_k}^2 \leq 2 a^\inv \paren[\big]{\sup_{j \in J} \abs{x_j}^2}a + 2 \sup_{j \in J} \abs{y_j}^2\]
by equation~\cref{e:weak-tri-inequality}, and so \(x a + y \in X\).
\end{proof}

\begin{lemma}
Equation \cref{e:dir-colim-inn-prod-mod} defines an \(A\)\nobreakdash-valued inner product on \(X\).
\end{lemma}

\begin{proof}
First, we check that the order limit in equation~\cref{e:dir-colim-inn-prod-mod} exists. Standard properties of order limits~\cite[Lemma~1.2]{hamana:tensor-products-monotone} and the complex polarisation identity yield
\[
    \olim_{j \in J} \innerProd{x_j}{y_j}
    = \olim_{j \in J} \frac{1}{4} \sum_{n = 0}^3 \mathrm{i}^n \abs{x_j\mathrm{i}^n + y_j}^2
    = \frac{1}{4} \sum_{n = 0}^3 \mathrm{i}^n \olim_{j \in J} \abs{x_j\mathrm{i}^n + y_j}^2,
\]
provided that the net \(\paren[\big]{\abs{x_j\mathrm{i}^n + y_j}^2}_{j \in J}\) order converges for each \(n \in \set{0, 1, 2, 3}\). But this net is increasing because
\[
    \abs{x_j \mathrm{i}^n + y_j}^2 = \abs{f_{j \leq k}(x_k \mathrm{i}^n + y_k)}^2 \leq \abs{x_k \mathrm{i}^n + y_k}^2
\]
for all \(j \leq k\) in \(J\), and it is bounded above because, by equation \cref{e:weak-tri-inequality},
\[\abs{x_j \mathrm{i}^n + y_j}^2 \leq 2\abs{x_j}^2 + 2\abs{y_j}^2 \leq 2 \sup_{k \in J} \abs{x_k}^2 + 2 \sup_{l \in J} \abs{y_l}^2\]
for all \(j \in J\), so it has a supremum and thus order converges to its supremum.

We now check the inner product axioms. For all \(x, y, z \in X\) and all \(a \in A\),
\begin{align*}
    \innerProd{x}{ya + z}
    &= \olim_{j \in J} \innerProd{x_j}{(ya + z)_j}
    \\&= \olim_{j \in J} \innerProd{x_j}{y_ja + z_j}
    \\&= \olim_{j \in J} \paren[\big]{\innerProd{x_j}{y_j}a + \innerProd{x_j}{z_j}}
    \\&= \paren[\big]{\olim_{j \in J} \innerProd{x_j}{y_j}}a + \olim_{j \in J} \innerProd{x_j}{z_j}
    \\&= \innerProd{x}{y}a + \innerProd{x}{z}
\end{align*}
by standard properties of order limits~\cite[Lemma~1.2]{hamana:tensor-products-monotone}.
Also, for all \(x, y \in X\), we have
\[\innerProd{x}{y}^\inv = \paren[\big]{\olim_{j \in J} \innerProd{x_j}{y_j}}^\inv = \olim_{j \in J} \innerProd{x_j}{y_j}^\inv = \olim_{j \in J} \innerProd{y_j}{x_j} = \innerProd{y}{x}.\]
As \(J\) is directed, it has an element \(k\), and so, for each \(x \in X\),
\[\abs{x}^2 = \sup_{j \in J} \abs{x_j}^2 \geq \abs{x_k}^2 \geq 0.\]
Finally, let \(x \in X\) and suppose that \(\abs{x}^2 = 0\). For each \(k \in J\),
\[0 \leq \abs{x_k}^2 \leq \sup_{j \in J} \abs{x_j}^2 = \abs{x}^2 = 0,\]
so \(\abs{x_k}^2 = 0\), and thus \(x_k = 0\). Hence \(x = 0\).
\end{proof}

\begin{proposition}
\label{p:un-complete}
The inner product \(A\)\nobreakdash-module \(X\) is ultranorm complete.
\end{proposition}

The proof below is a generalisation of the standard proof that the inner product space \(\ell^2(\Nats)\) of square-summable sequences is a Hilbert space (see \cref{s:hilbert-spaces}).

\begin{proof}
Let \((x^{(s)})_{s \in S}\) be a net in \(X\) that is ultranorm Cauchy.

Fix \(j \in J\). For all positive linear functionals \(\omega\) on \(A\) and all \(s, t \in S\),
\begin{equation}
    \label{e:piece-to-whole-mod}
    \begin{aligned}
        \norm[\big]{x^{(s)}_j - x^{(t)}_j}_\omega^2
        &= \omega \abs[\big]{x^{(s)}_j - x^{(t)}_j}^2
        \\&\leq \omega \sup_{k \in J} \abs[\big]{x^{(s)}_k - x^{(t)}_k}^2
        = \omega \abs{{x^{(s)} - x^{(t)}}}^2
        = \norm{x^{(s)} - x^{(t)}}_\omega^2.
    \end{aligned}
\end{equation}
Hence the net \((x^{(s)}_j)_{s \in S}\) in \(X_j\) is ultranorm Cauchy. It ultranorm converges because \(X_j\) is ultranorm complete. Let
\[x_j = \unlim_{s \in S} x^{(s)}_j.\]
We will show that \(x = (x_j)_{j \in J}\) is an element of \(X\) and that \(\unlim_{s \in S} x^{(s)} = x\). First
\[f_{j \leq k} x_k = f_{j \leq k} \unlim_{s \in S} x^{(s)}_k = \unlim_{s \in S} f_{j \leq k} x^{(s)}_k = \unlim_{s \in S} x^{(s)}_j = x_j\]
because contractions are ultranorm continuous~\cite[148~I]{westerbaan:dagger-dilation-category}. Hence \(x \in \widebar{X}\).

Let \(\varepsilon \in \Reals_{>0}\), and let \(\omega\) be a normal positive linear functional on \(A\). As \((x^{(s)})_{s \in S}\) is ultranorm Cauchy, there exists \(s_{\omega, \varepsilon} \in S\) such that \[\norm{x^{(s)} - x^{(t)}}_\omega < \varepsilon\] for all \(s, t \geq s_{\omega, \varepsilon}\). Using \cref{e:piece-to-whole-mod}, it follows that
\begin{align*}
    \norm[\big]{x^{(s)}_j - x_j}_\omega
    &= \norm[\big]{\paren[\big]{x^{(s)}_j - x^{(t)}_j} + \paren[\big]{x^{(t)}_j - x_j}}_\omega
    \\&\leq \norm[\big]{x^{(s)}_j - x^{(t)}_j}_\omega + \norm[\big]{x^{(t)}_j - x_j}_\omega
    \\&\leq \norm{x^{(s)} - x^{(t)}}_\omega + \norm[\big]{x^{(t)}_j - x_j}_\omega
    \\&< \varepsilon + \norm[\big]{x^{(t)}_j - x_j}_\omega
\end{align*}
for all \(j \in J\) and all \(s, t \geq s_{\omega, \varepsilon}\). Thus
\[\norm[\big]{x^{(s)}_j - x_j}_\omega \leq \varepsilon + \lim_{t \geq s_{\omega, \varepsilon}}\norm[\big]{x^{(t)}_j - x_j}_\omega = \varepsilon + 0 = \varepsilon\]
for all \(j \in J\) and all \(s \geq s_{\omega, \varepsilon}\). Hence, for all \(s \geq s_{\omega, \varepsilon}\),
\[
    \sup_{j \in J}\norm[\big]{x^{(s)}_j - x_j}_\omega \leq \varepsilon.
\]
and thus
\[
    \sup_{j \in J} \norm{x_j}_\omega
    \leq \sup_{j \in J} \norm{x_j - x^{(s)}_j}_\omega + \sup_{j \in J} \norm{x_j^{(s)}}_\omega
    \leq \varepsilon + \norm{x^{(s)}}_\omega
    < \infty,
\]
that is,
\[
    \sup_{j \in J} \omega \abs{x_j}^2 < \infty.
\]
Since \(\omega\) is arbitrary, it follows by \cref{l:positive-boundedness} that 
\[
    \sup_{j \in J} \abs{x_j}^2 < \infty,
\]
that is, that \(x \in X\). Thus
\[
    \norm{x^{(s)} - x}_\omega = \sup_{j \in J}\norm{x^{(s)}_j - x_j}_\omega \leq \varepsilon
\]
for all \(s \geq s_{\omega, \varepsilon}\). Hence
\[\unlim_{s \in S} x^{(s)} = x.\]
Thus \(X\) is ultranorm complete.
\end{proof}

\begin{lemma}
The morphisms \((f_j \colon X \to X_j)_{j \in J}\) form a limit in \(\Hilb_{A, \leq 1}\) of the diagram \((f_{j \leq k} \colon X_k \to X_j)_{j,k \in J}\).
\end{lemma}

\begin{proof}
For each \(j \in J\) and each \(x \in X\),
\[\abs{f_j x}^2 = \abs{x_j}^2 \leq \sup_{k \in J} \abs{x_k}^2 = \abs{x}^2,\]
so \(f_j\) is a contraction. The maps \(f_j\) are the composites of the submodule inclusion \(X \hookrightarrow \widebar{X}\) with the limit projections \(\bar{f}_j \colon \widebar{X} \to X_j\) in \(\Mod_A\). It follows that they are \(A\)-linear and thus adjointable~\cite[152~VIII]{westerbaan:dagger-dilation-category}, and hence also that they form a cone on the diagram \(f_{j \leq k}\) in \(\Hilb_{A, \leq 1}\).

Let \((g_j \colon Y \to X_j)_{j \in J}\) be another cone in \(\Hilb_{A, \leq 1}\) on this diagram. For each \(y \in Y\), as
\[f_{j \leq k} g_k y = g_j y\]
for all \(j \leq k\) in \(J\) and
\[\abs{g_j y}^2 \leq \abs{y}^2\]
for all \(j \in J\), we may define the function \(g \colon Y \to X\) by
\[gy = (g_j y)_{j \in J}.\]
It is easy to check that \(g\) is \(A\)-linear. It is a contraction because
\[\abs{gy}^2 = \sup_{j \in J} \abs{g_j y}^2 \leq \abs{y}^2.\]
Since \(g\) is bounded and \(A\)-linear, it is adjointable~\cite[152~VIII]{westerbaan:dagger-dilation-category}. It is also easy to check uniqueness.
\end{proof}

This completes the proof of \cref{p:hilb-mod}.
\newpage
\printbibliography

@article{handelman:rings-involution-partially,
    Author = {David Handelman},
    Title = {Rings with involution as partially ordered abelian groups},
    Journal = {Rocky Mountain Journal of Mathematics},
    Year = {1981},
    Volume = {11},
    Number = {3},
    Pages = {337--382},
    DOI = {10.1216/RMJ-1981-11-3-337}
}

@unpublished{di-meglio:r-star-cats,
    Author = {Di Meglio, Matthew},
    Title = {Pre-Hilbert \(*\)-categories: The Hilbert-space analogue of abelian categories},
    Year = {2023},
    eprint = {2312.02883},
    archivePrefix = {arXiv}
}

@article{di-meglio:fcon,
    Author = {Di Meglio, Matthew and Heunen, Chris},
    Title = {Dagger categories and the complex numbers: Axioms for the category of finite-dimensional Hilbert spaces and linear contractions},
    Year = {2025},
    journal = {Applied Categorical Structures},
    volume = {33},
    pages = {18},
    doi = {10.1007/s10485-025-09803-5},
    eprint = {2401.06584},
    archivePrefix = {arXiv}
}

@article{heunen:quantum-logic-dagger-order,
    Author = {Heunen, Chris and Jacobs, Bart},
    Title = {Quantum logic in dagger kernel categories},
    Journal = {Order},
    Year = {2010},
    Volume = {27},
    Number = {2},
    Pages = {177--212},
    DOI = {10.1007/s11083-010-9145-5}
}

@article{demarr:partially-ordered-fields,
    Author = {R. {DeMarr}},
    Title = {Partially ordered fields},
    Journal = {American Mathematical Monthly},
    Year = 1967,
    Volume = 74,
    Number = 4,
    Pages = {418--420},
    DOI = {10.2307/2314578}
}

@article{prestel:solers-characterization-hilbert,
    Author = {Prestel, Alexander},
    Title = {On Solèr's characterization of Hilbert spaces},
    Journal = {Manuscripta Mathematica},
    Year = {1995},
    Volume = {86},
    Number = {1},
    Pages = {225--238},
    DOI = {10.1007/BF02567991}
}

@book{adamek:locally-presentable-accessible,
    Title = {Locally Presentable and Accessible Categories},
    Publisher = {Cambridge University Press},
    Year = {1994},
    Author = {Adamek, J. and Rosicky, J.},
    Series = {London Mathematical Society Lecture Note Series},
    DOI = {10.1017/CBO9780511600579}
}

@article{manuilov:hilbert-w-modules-their,
    Author = {Manuilov, V. M. and Troitsky, E. V.},
    Title = {Hilbert C*- and W*-modules and their morphisms},
    Journal = {Journal of Mathematical Sciences},
    Year = {2000},
    Volume = {98},
    Number = {2},
    Pages = {137--201},
    DOI = {10.1007/BF02355447}
}

@book{borceux:handbook-categorical-algebra-1,
    Title = {Handbook of Categorical Algebra},
    Publisher = {Cambridge University Press},
    Year = {1994},
    Author = {Borceux, Francis},
    Volume = {1},
    Series = {Encyclopedia of Mathematics and its Applications},
    DOI = {10.1017/CBO9780511525858}
}

@book{borceux:handbook-categorical-algebra-2,
    Title = {Handbook of Categorical Algebra},
    Publisher = {Cambridge University Press},
    Year = {1994},
    Author = {Borceux, Francis},
    Volume = {2},
    Series = {Encyclopedia of Mathematics and its Applications},
    DOI = {10.1017/CBO9780511525865}
}

@phdthesis{westerbaan:dagger-dilation-category,
    Author = {Westerbaan, B. E.},
    Title = {Dagger and dilation in the category of von Neumann algebras},
    Institution = {Radboud University},
    Year = {2019},
    eprint={1803.01911},
    archivePrefix={arXiv},
    pagination={none}
}

@article{vicary:completeness--categories-complex,
    Author = {Vicary,Jamie},
    Title = {Completeness of $\dagger$-categories and the complex numbers},
    Journal = {Journal of Mathematical Physics},
    Year = {2011},
    Volume = {52},
    Number = {8},
    Pages = {082104},
    DOI = {10.1063/1.3549117}
}

@article{soler:characterization-hilbert-spaces,
    Author = {Solèr, M.\ P.},
    Title = {Characterization of Hilbert spaces by orthomodular spaces},
    Journal = {Communications in Algebra},
    Year = {1995},
    Volume = {23},
    Number = {1},
    Pages = {219-243},
    DOI = {10.1080/00927879508825218}
}

@book{riehl:category-theory-context,
    Title = {Category Theory in Context},
    Publisher = {Dover Publications Inc.},
    Year = 2016,
    Author = {E. Riehl},
}

@article{heunen:con,
    Author = {Heunen, Chris and Kornell, Andre and van der Schaaf, Nesta},
    Title = {Axioms for the category of Hilbert spaces and linear contractions},
    Journal = {Bulletin of the London Mathematical Society},
    Year = {2024},
    Pages = {1532-1549},
    DOI={10.1112/blms.13010}
}

@article{heunen:hilb,
    Author = {Chris Heunen and Andre Kornell},
    Title = {Axioms for the category of Hilbert spaces},
    Journal = {Proceedings of the National Academy of Sciences},
    Year = {2022},
    Volume = {119},
    Number = {9},
    Pages = {e2117024119},
    DOI = {10.1073/pnas.2117024119}
}

@article{freyd:abelian-categories,
    Author = {P. J. Freyd},
    Title = {Abelian categories},
    Journal = {Reprints in Theory and Applications of Categories},
    Year = {2003},
    Number = {3},
    Pages = {23--164},
}

@inproceedings{abramsky:categorical-semantics-quantum,
    Author = {Abramsky, S. and Coecke, B.},
    Title = {A categorical semantics of quantum protocols},
    Booktitle = {Proceedings of the 19th Annual IEEE Symposium on Logic in Computer Science},
    Year = {2004},
    Pages = {415-425},
    Booktitle = {Proceedings of the 19th Annual IEEE Symposium on Logic in Computer Science},
    DOI = {10.1109/LICS.2004.1319636}
}

@article{ghez:w-categories,
    Author = {P. Ghez and R. Lima and J. E. Roberts},
    Title = {W*-categories},
    Journal = {Pacific Journal of Mathematics},
    Year = {1985},
    Volume = {120},
    Number = {1},
    Pages = {79-109},
    DOI = {10.2140/pjm.1985.120.79}
}

@article{schotz:equivalence-order-algebraic,
    Author = {Schötz, Matthias},
    Title = {Equivalence of order and algebraic properties in ordered $*$-algebras},
    Journal = {Positivity},
    Year = {2021},
    Volume = {25},
    Number = {3},
    Pages = {883--909},
    DOI = {10.1007/s11117-020-00792-4}
}

@book{saitowright:monotone-complete-generic,
    Title = {Monotone Complete C*-algebras and Generic Dynamics},
    Publisher = {Springer London},
    Year = {2015},
    Author = {Kazuyuki Saitô and {J. D. Maitland} Wright},
    DOI = {10.1007/978-1-4471-6775-4}
}

@book{berberian:baer-star-rings,
  author         = {S. K. Berberian},
  title          = {Baer $*$-Rings},
  year           = {1972},
  publisher      = {Springer},
  doi            = {10.1007/978-3-642-15071-5}
}

@phdthesis{westerbaan:2019:category,
      title={The category of von Neumann algebras}, 
      author={Abraham A. Westerbaan},
      year={2019},
      institution={Radboud University},
      eprint={1804.02203},
      archivePrefix={arXiv},
      pagination={none}
}

@article{tobin-lack:hilb,
      title={A characterisation for the category of Hilbert spaces}, 
      author={Lack, Stephen and Tobin, Shay},
      year={2025},
	  journal={Applied Categorical Structures},
	  volume={33},
	  number={13},
	  doi={10.1007/s10485-025-09805-3}
}

@article{hamana:tensor-products-monotone,
    Author = {Masamichi Hamana},
    Title = {Tensor products for monotone complete C*-algebras, I},
    Journal = {Japanese Journal of Mathematics},
    Year = {1982},
    Volume = {8},
    Number = {2},
    Pages = {259-283},
    DOI = {10.4099/math1924.8.259}
}

@article{kadison:equivalence-operator-algebras,
    Author = {Richard V. Kadison and Gert Kjærgård Pedersen},
    Title = {Equivalence in operator algebras},
    Journal = {Mathematica Scandinavica},
    Year = {1970},
    Volume = {27},
    Number = {2},
    Pages = {205--222},
    DOI = {10.7146/math.scand.a-10999}
}

@article{widom:embedding-algebras,
    Author = {Harold Widom},
    Title = {Embedding in algebras of type I},
    Journal = {Duke Mathematical Journal},
    Year = {1956},
    Volume = {23},
    Number = {2},
    Pages = {309 -- 324},
    DOI = {10.1215/S0012-7094-56-02330-4}
}

@article{kaplansky:modules-operator-algebras,
    Author = {Irving Kaplansky},
    Title = {Modules over operator algebras},
    Journal = {American Journal of Mathematics},
    Year = {1953},
    Volume = {75},
    Number = {4},
    Pages = {839--858},
    DOI={10.2307/2372552}
}

@article{arens:space-lomega-convex,
    Author = {Richard Arens},
    Title = {The space $L^\omega$ and convex topological rings},
    Journal = {Bulletin of the American Mathematical Society},
    Year = {1946},
    Volume = {52},
    Number = {10},
    Pages = {931-935},
    DOI={10.1090/S0002-9904-1946-08681-4}
}

@article{kadison:operator-algebras-faithful,
    Author = {Richard V. Kadison},
    Title = {Operator algebras with a faithful weakly-closed representation},
    Journal = {Annals of Mathematics},
    Year = {1956},
    Volume = {64},
    Number = {1},
    Pages = {175--181},
    DOI = {10.2307/1969954}
}

@book{ murphy:cstaralgebras,
  author = {G. J. Murphy},
  publisher = {Academic Press},
  title = {C*-Algebras and Operator Algebras},
  year = {1990}
}

@article{ghiloni:continuous-slice-functional,
    Author = {Ghiloni, Riccardo and Moretti, Valter and Perotti, Alessandro},
    Title = {Continuous slice functional calculus in quaternionic Hilbert spaces},
    Journal = {Reviews in Mathematical Physics},
    Year = {2013},
    Volume = {25},
    Number = {04},
    Pages = {1350006},
    DOI = {10.1142/S0129055X13500062}
}

@book{popescu:abelian-categories,
    author={N. Popescu},
    title={Abelian Categories with Applications to Rings and Modules},
    publisher={Academic Press},
    year={1973}
}

@article{grothendieck:quelques-points-dalgebre,
    Author = {Alexander Grothendieck},
    Title = {Sur quelques points d'algèbre homologique, I},
    Journal = {Tohoku Mathematical Journal},
    Year = {1957},
    Volume = {9},
    Number = {2},
    Pages = {119-221},
    DOI = {10.2748/tmj/1178244839}
}

@BOOK{heunenvicary:quantum,
    author = {C. Heunen and J. Vicary},
    title = {Categories for Quantum Theory},
    publisher = {Oxford University Press},
    year = {2019}
}

@article{milius:colimits-categories-relations,
    Author = {Milius, Stefan},
    Title = {On colimits in categories of relations},
    Journal = {Applied Categorical Structures},
    Year = {2003},
    Volume = {11},
    Number = {3},
    Pages = {287--312},
    DOI = {10.1023/A:1024201720423}
}

@article{frei:limits-categories-relations,
    Author = {Frei, Armin and MacDonald, John L.},
    Title = {Limits in categories of relations and limit-colimit commutation},
    Journal = {Journal of Pure and Applied Algebra},
    Year = {1971},
    Volume = {1},
    Number = {2},
    Pages = {179-197},
    DOI = {10.1016/0022-4049(71)90017-X}
}

@article{fritz:universal-property-infinite,
    Author = {Fritz, Tobias and Westerbaan, Bas},
    Title = {The universal property of infinite direct sums in C*-categories and W*-categories},
    Journal = {Applied Categorical Structures},
    Year = {2020},
    Volume = {28},
    Number = {2},
    Pages = {355--365},
    DOI = {10.1007/s10485-019-09583-9}
}

@article{cimpric:2009:quadratic-module,
    title={A representation theorem for Archimedean quadratic modules on \(*\)-rings},
    volume={52},
    DOI={10.4153/CMB-2009-005-4},
    number={1},
    journal={Canadian Mathematical Bulletin},
    author={Cimprič, Jakob},
    year={2009},
    pages={39–52}
}

@article{albeverio:partially-ordered-involutory,
    Author = {Albeverio, S. and Ayupov, Sh. A. and Dadakhodjayev, R. A.},
    Title = {On partially ordered real involutory algebras},
    Journal = {Acta Applicandae Mathematica},
    Year = {2006},
    Volume = {94},
    Number = {3},
    Pages = {195--214},
    DOI = {10.1007/s10440-006-9074-x}
}

@article{saito:defining-aw-algebras-rickart,
    Author = {Saitô, Kazuyuki; Wright, J.D. Maitland},
    Title = {On defining AW*-algebras and Rickart C*-algebras},
    Journal = {The Quarterly Journal of Mathematics},
    Year = {2015},
    Volume = {66},
    Number = {3},
    Pages = {979-989},
    DOI = {10.1093/qmath/hav015}
}

@book{halmos:introduction-hilbert,
  author = {Halmos, Paul R.},
  publisher = {Dover Publications},
  title = {Introduction to Hilbert Space and the Theory of Spectral Multiplicity},
  edition = {2},
  year = {2017},
  isbn = {9780486817330}
}

@article{abian:direct-product,
    Author = {Alexander Abian},
    Title = {Direct product decomposition of commutative semisimple rings},
    Journal = {Proceedings of the American Mathematical Society},
    Year = {1970},
    Volume = {24},
    Number = {3},
    Pages = {502--507},
    DOI = {10.2307/2037396}
}

@article{gabriel-popescu,
    Author = {Popesco, Nicolae and Gabriel, Pierre},
    Title = {Caractérisation des catégories abéliennes avec générateurs et limites inductives exactes},
    Journal = {Comptes Rendus Hebdomadaires des Séances de l'Académie des Sciences},
    Volume={258},
    Year={1964},
    Pages={4188-4190},
    URL={https://gallica.bnf.fr/ark:/12148/bpt6k4011c/f1826.item},
    note={French}}

@book{freyd:allegories,
  author = {Freyd, Peter J. and Scedrov, Andre},
  publisher = {Elsevier Science Publishers},
  title = {Categories, Allegories},
  year = {1990},
  ISBN={0-444-70367-5},
  pagination={section}
}

@article{isbell:functorial,
    Author = {J. R. Isbell},
    Title = {General functorial semantics, I},
    Journal = {American Journal of Mathematics},
    Year = {1972},
    Volume = {94},
    Number = {2},
    Pages = {535--596},
    DOI = {10.2307/2374638}
}

@book{riehl:homotopy,
    author={Riehl, Emily},
    title={Categorical Homotopy Theory},
    publisher={Cambridge University Press},
    year={2014},
    doi={10.1017/CBO9781107261457}
}

@inproceedings{vicary:conn-lims,
    author = {Reutter, David and Vicary, Jamie},
    title = {High-level methods for homotopy construction in associative \(n\)-categories},
    year = {2021},
    publisher = {IEEE Press},
    booktitle = {Proceedings of the 34th Annual ACM/IEEE Symposium on Logic in Computer Science},
    number = {62},
    DOI={10.1109/lics.2019.8785895}}
\end{document}